\documentclass[a4paper, 12pt]{article}
\usepackage{enumerate, theorem}
\usepackage{amsmath, amsfonts, amssymb}
\usepackage[height=22.5cm, width=15cm]{geometry}
\usepackage[all]{xy}
\usepackage{mathrsfs, yfonts}

\DeclareMathOperator{\id}{id}

\DeclareMathOperator{\C}{\mathbb{C}}
\DeclareMathOperator{\ptc}{\widehat{\otimes}}

\newcommand{\A}{\tilde{\mathcal{A}}}
\newcommand{\p}{\Lambda}

\newcommand{\parag}[1]{\paragraph{\sc{#1.}}}

\newtheorem{thm}{Th\'eor\`eme}[subsection]
\newtheorem{defn}[thm]{D\'efinition}
\newtheorem{cor}[thm]{Corollaire}
\newtheorem{prop}[thm]{Proposition}
\newtheorem{lemma}[thm]{Lemme}

\begin{document}

\title{Le th\`eme d'une p\'eriode \'evanescente.}

\author{Daniel Barlet\footnote{Barlet Daniel, Institut Elie Cartan UMR 7502  \newline
Nancy-Universit\'e, CNRS, INRIA  et  Institut Universitaire de France, \newline
BP 239 - F - 54506 Vandoeuvre-l\`es-Nancy Cedex.France. r\newline
e-mail : barlet@iecn.u-nancy.fr}.}

\date{28/11/09 corrig\'ee.}

\maketitle


 

 \markright{Le th\`eme d'une p\'eriode \'evanescente.}

\section*{Abstract.}

In this article we study holomorphic deformations of the filtered Gauss-Manin systems associated to a vanishing period integral. For that purpose we introduce a new sub-class of the class of monogenic (a,b)-modules (Brieskorn modules) which was studied in our previous article [B. 09]. We show that these new objects, called "themes", have good functorial properties and that there exists a canonical order on the roots of the corresponding Bernstein polynomial.\\
We construct, for given fundamental invariants, a finite dimensional versal holomorphic family and we show that, when all themes with these fundamental invariants are "stable", this versal family is in fact universal. We also give a sufficient condition on the roots of the Bernstein polynomial in order that the previous condition is satisfied. We show with  an example that a universal family may not exist for some values of the fundamental invariants.

\parag{Key words} Vanishing period, Berstein polynomial, filtered Gauss-Manin system, (a,b)-module, Brieskorn module.

\parag{AMS Classification 2000} 32S05, 32S25, 32S40.

\tableofcontents

\section*{Introduction.}

Mon article pr\'ec\'edent [B.09]  est focalis\'e sur les (a,b)-modules monog\`enes dans l'id\'ee d'\'etudier, pour un \'el\'ement \ $x$ \ donn\'e dans un (a,b)-module r\'egulier, ou de mani\`ere plus g\'eom\'etrique,  pour une forme holomorphe donn\'ee dont on veut \'etudier la p\'eriode \'evanescente (voir [A-G-V], [M. 75], [S. 89]), le sous-(a,b)-module engendr\'e par \ $x$. Concr\`etement cela signifie que l'on se concentre sur l'\'equation diff\'erentielle (filtr\'ee) minimale satisfaite par les fonctions obtenues par int\'egration de cette forme sur les cycles \'evanescents.\\
Le pr\'esent article se consacre \`a l'\'etude plus pr\'ecise d'une int\'egrale de p\'eriode \'evanescente, ce qui revient cette fois \`a fixer la classe d'homologie \'evanescente sur laquelle on int\`egre. Ceci conduit \`a une sous-classe int\'eressante des (a,b)-modules monog\`enes r\'eguliers \'etudi\'es dans  [B.09] qui est caract\'eris\'ee par une propri\'et\'e alg\'ebrique remarquable et simple dans le cas \ $[\lambda]-$primitif :  "L'unicit\'e de la suite de Jordan-H{\"o}lder".

Les \'el\'ements de cette sous-classe que j'ai appel\'es des {\bf th\`emes} correspondent en fait \`a la construction "na{\"i}ve" suivante :\\
Consid\'erons un sous-ensemble fini \ $\Lambda \subset ]0,1] \cap \mathbb{Q}$ \ (r\'eduit \`a \ $\{\lambda\}$ \ dans le cas \ $[\lambda]-$primitif ) et un entier \ $N$, et consid\'erons l'espace des s\'eries formelles
$$ \Xi^{(N)}_{\Lambda} : = \sum_{\lambda \in \Lambda, j \in [0,N]} \ \C[[s]].s^{\lambda-1}\frac{(Log\, s)^j}{j!} .$$
D\'efinissons la \ $\C-$alg\`ebre non commutative \ $\A$ \ en posant
$$ \A : = \{ \sum_{\nu = 0}^{\infty} \ P_{\nu}(a).b^{\nu} , \quad P_{\nu} \in \C[x] \} $$
avec la relation de commutation \ $a.b - b.a = b^2$. \\
On a une action naturelle de \ $\A$ \ sur \ $\Xi^{(N)}_{\Lambda}$ \ via les actions donn\'ees respectivement par la multiplication par \ $s$ \ ($ a : = \times s$)  et  et l'int\'egration sans constante ( $b : = \int_0^s $). Un th\`eme sera alors un sous$-\A-$module monog\`ene d'un tel \ $\Xi^{(N)}_{\Lambda}$ \  c'est-\`a-dire le sous$-\A-$module \`a gauche engendr\'e par un \'el\'ement \ $\varphi \in \Xi^{(N)}_{\Lambda}$.\\
En pr\'esence d'un (a,b)-module g\'eom\'etrique \ $E$ \ et d'une application (a,b)-lin\'eaire \ $\Gamma : E \to \Xi^{(N)}_{\Lambda}$ \  l'image par \ $\Gamma$ du (a,b)-module monog\`ene \ $\A.x$ \ engendr\'e par \ $x$ \ dans \ $E$, sera un th\`eme.\\
Par exemple si \ $E$ \ est  le compl\'et\'e formel en \ $f$ \ du module de Brieskorn d'un germe de fonction \ $f$ \ holomorphe  \`a singularit\'e isol\'ee dans \ $\C^{n+1}$ \ l'application (a,b)-lin\'eaire \ $\Gamma : E \to \Xi^{(N)}_{\Lambda}$ \  associ\'ee \`a un cycle \'evanescent \ $\gamma$ \  qui  fait correspondre  \`a \ $[\omega]$ \ le d\'eveloppement asymptotique (formel) de la fonction multiforme de d\'etermination finie \ $s \to \int_{\gamma_s} \omega/df $, o\`u \ $(\gamma_s)_{s \in D^*}$ \ d\'esigne la famille horizontale multiforme associ\'ee \`a \ $\gamma$ \ dans les fibres de \ $f$, et o\`u l'on a choisi convenablement \ $\Lambda$ \ et \ $N$.\\

Le polyn\^ome de Bernstein d'un (a,b)-module monog\`ene r\'egulier \'etant d\'ecrit en terme du g\'en\'erateur de l'id\'eal annulateur d'un g\'en\'erateur du (a,b)-module monog\`ene consid\'er\'e, nous proposons dans cet article d'\'etudier les invariants (num\'eriques) plus fins que le polyn\^ome de Bernstein d'un th\`eme. En fait nous d\'ecrirons tous les invariants associ\'es \`a une classe d'isomorphisme de th\`eme primitif. Pour ce faire nous \'etudierons les familles holomorphes de th\`emes, ce qui correspond  \`a l'\'etude d'une p\'eriode \'evanescente d\'ependant holomorphiquement d'un param\`etre. Par exemple ce ph\'enom\`ene appara\^it dans le cas d'une famille \`a \ $\mu$ \ constant de fonctions holomorphes  \`a singularit\'es isol\'ees, quand on consid\`ere une forme holomorphe (relative) et une classe d'homologie fix\'ee dans une  fibre lisse.\\
Notre objectif principal sera de d\'ecrire  concr\`etement des familles holomorphes  {\bf verselles}  (et "minimales") pour les th\`emes \ $[\lambda]-$primitifs. Nous montrerons que dans le "cas stable", le seul o\`u l'on peut esp\'erer en g\'en\'eral l'existence d'une famille universelle, les familles d\'ecrites sont effectivement universelles.

\bigskip

Les principaux r\'esultats de ce travail sont les suivants.

\begin{enumerate}
\item Les th\'eor\`emes \ref{Stabilite par quotient} et \ref{dual d'un theme} de stabilit\'e des th\`emes  par quotient et dualit\'e "tordue" qui permettront de montrer qu'un(a,b)-module monog\`ene est un th\`eme si et seulement pour chacun de ses exposants \ $[\lambda]$ \ sa partie \ $[\lambda]-$primitive est un th\`eme \ref{caract. theme 2}.
\item La carat\'erisation des th\`emes stables et le th\'eor\`eme d'unicit\'e \ref{unicite} de l'\'ecriture du g\'en\'erateur de l'id\'eal annulant un g\'en\'erateur standard dans le cas d'un th\`eme primitif stable. Ceci donne l'universalit\'e de la famille standard quand elle ne contient que des th\`emes stables. Une condition suffisante simple (voir le corollaire \ref{Endo}) sur les invariants fondamentaux donn\'ee au th\'eor\`eme \ref{inclusion}  permet d'assurer que c'est souvent  le cas.
\item L'existence des bases standards qui donneront la construction de familles verselles de th\`emes \ $[\lambda]-$primitifs, une fois fix\'es les invariants fondamentaux.
\item Nous terminons par un exemple en rang 3  pour lequel nous montrons qu'il n'existe pas de famille universelle au voisinage de chaque th\`eme stable ayant ces invariants fondamentaux. Ces th\`emes stables sont param\'etr\'es dans cet exemple par une hypersurface (non vide) de la famille standard.\\
 Par contre, une fois enlever cette hypersurface, on peut construire une famille qui est universelle au voisinage de chacun de ses points et param\`etre tous les th\`emes instables ayant ces invariants fondamentaux.
\end{enumerate}

\section{D\'ecomposition primitive.}

\subsection{Rappels.}

Soit \ $A$ \ le quotient de l'alg\`ebre libre \ $\C<a,b>$ \ par l'id\'eal bilat\`ere engendr\'e par \ $a.b - b.a - b^2$. On notera que pour chaque \ $k \in \mathbb{N}$ \ $b^k.A = A.b^k$ \ est un id\'eal bilat\`ere de \ $A$. Soit \ $\A$ \ la compl\'et\'ee \ $b-$adique de \ $A$. On a alors
$$ \A : = \{ \sum_{\nu \geq 0} P_{\nu}(a).b^{\nu} , P_{\nu} \in \C[x] \} .$$
C'est une \ $\C-$alg\`ebre unitaire int\`egre qui contient la sous-alg\`ebre commutative \ $\C[[b]]$. On appelle (a,b)-module \ $E$ \ un \ $\A-$module \`a gauche qui est libre de type fini sur \ $\C[[b]]$. 
Se donner un (a,b)-module \'equivaut \`a la donn\'ee d'un \ $\C[[b]]-$module libre de rang fini \ $E$ \ et d'une application \ $\C-$lin\'eaire \ $ a : E \to E $ \ v\'erifiant la relation de commutation  \ $a.b - b.a = b^2$ ; elle est  continue pour la topologie \ $b-$adique de \ $E$.
Une telle application \ $\C-$lin\'eaire \ $ a $ \ est d\'etermin\'ee de fa{\c c}on unique par les valeurs de \ $a$ \ sur  une \ $\C[[b]]-$base de \ $E$, et elles peuvent \^etre choisies arbitrairement dans \ $E$ : en effet, si \ $E : = \oplus_{j=1}^k \C[[b]].e_j$ \ et si on s'est donn\'e arbitrairement des \'el\'ements \ $x_1, \dots, x_k$ \ dans \ $E$, l'application \ $a$ \ associ\'ee est bien d\'efinie sur \ $E_0 : = \oplus_{j=1}^k \C[b].e_j$ \ par les relations
$$ a.b^n.e_j = b^n.x_j + n.b^{n+1}.e_j \quad \forall n \in \mathbb{N}\quad  \forall j \in [1,k] $$
qui sont cons\'equences de \ $a.b - b.a = b^2 $ \ et de \ $a.e_j = x_j\quad \forall j \in [1,k] $.
 L'application \ $a$ \  se prolonge alors de fa{\c c}on unique \`a \ $E$ \ par continuit\'e.

\bigskip

On dit que le (a,b)-module \ $E$ \ est {\bf \`a p\^ole simple} si on a \ $a.E \subset b.E$. \\
On dit que \ $E$ \ est {\bf r\'egulier} s'il se plonge dans un (a,b)-module \`a p\^ole simple. Dans ce cas le plus petit (a,b)-module \`a p\^ole simple contenant \ $E$ \ est le satur\'e 
 $$E^{\sharp} \subset E[b^{-1}] :  = E \otimes_{\C[[b]]} \C[[b]][b^{-1}]$$
 de \ $E$ \ par \ $b^{-1}.a$. La r\'egularit\'e de \ $E$ \ est \'equivalente \`a la finitude sur \ $\C[[b]]$ \ de ce satur\'e 
 $$E^{\sharp} = \sum_{j\geq 0} \ (b^{-1}.a)^j.E \subset E[b^{-1}] .$$

Pour \ $E$ \ \`a p\^ole simple on d\'efinit le {\bf polyn\^ome de Bernstein} de \ $E$, not\'e \ $B_E$,  comme le polyn\^ome minimal de \ $-b^{-1}.a$ \ agissant sur l'espace vectoriel de dimension finie \ $E\big/b.E$.\\
Plus g\'en\'eralement, le polyn\^ome de Bernstein d'un (a,b)-module r\'egulier \ $E$ \ est d\'efini comme le polyn\^ome de Bernstein de son satur\'e par \ $b^{-1}.a$. Donc \ $B_E : = B_{E^{\sharp}}$. 

\bigskip

On dit qu'un (a,b)-module r\'egulier est {\bf g\'eom\'etrique} si toutes les racines de son polyn\^ome de Bernstein sont des rationnels strictement n\'egatifs.

\bigskip

Si \ $E$ \ et \ $F$ \ sont deux (a,b)-modules, on d\'efinit leur produit tensoriel \ $E \otimes_{a,b} F$ \ en consid\'erant le produit tensoriel des deux \ $\C[[b]]-$modules correspondants (qui est bien libre de type fini sur \ $\C[[b]]$), et en d\'efinissant \ $a : E \otimes_{a,b} F \to E \otimes_{a,b} F$ \ par la formule 
$$ a(e \otimes f) : = (a.e) \otimes f + e \otimes (a.f) .$$
On v\'erifie alors facilement que l'on a bien \ $a.b - b.a = b^2$.

De m\^eme, si \ $E$ \ et \ $F$ \ sont deux (a,b)-modules, on d\'efinit \ $Hom_{a,b}(E,F)$ \ en consid\'erant  le \ $\mathbb{C}[[b]]-$module \ $Hom_b(E,F)$ \ des applications \ $\C[[b]]-$lin\'eaires de \ $E$ \ dans \ $F$ \ et en d\'efinissant, pour \ $\varphi \in Hom_b(E,F)$ \ et \ $x \in E$ :\\
\begin{equation*}
(a.\varphi)(x) : = a.\varphi(x) - \varphi(a.x)  \tag{1}
\end{equation*}
on a alors la \ $\C[[b]]-$lin\'earit\'e de \ $(a.\varphi)$ \ et l'identit\'e
 \ $(a.b - b.a).\varphi = b^2.\varphi$.

\bigskip

On appellera dual de \ $E$ \ le (a,b)-module \ $Hom_{a,b}(E,E_0)$ \ o\`u \ $E_0 : = \A\big/\A.a$ \ est le (a,b)-module de rang 1 , de g\'en\'erateur \ $e_0$ \ v\'erifiant \ $a.e_0 = 0$. Le lecteur v\'erifiera facilement que pour \ $\lambda \in \C$ \ le dual de \ $E_{\lambda} : = \A\big/\A.(a - \lambda.b)$, est \ $(E_{\lambda})^* \simeq E_{-\lambda}$. On peut facilement en d\'eduire que pour \ $E$ \ r\'egulier on a canoniquement \ $(E^*)^* \simeq E$.

\parag{Exemple} Si \ $E$ \ est un (a,b)-module et si \ $E_{\delta}$ \ est le (a,b)-module de rang \ $1$ \ et de g\'en\'erateur \ $e_{\delta}$ \ tel que \ $a.e_{\delta} = \delta.b.e_{\delta} $ \ (donc \ $E_{\delta} \simeq \A\big/\A.(a - \delta.b)$), le (a,b)-module \ $E \otimes_{a,b} E_{\delta}$ \ peut \^etre identifi\'e au \ $\C[[b]]-$module \ $E$ \ dans lequel on a d\'efini l'action de \ $"a"$ \ par \ $\tilde{a} : = a + \delta.b$.\\
On remarquera que pour chaque \ $\delta \in \C$ \  il existe un unique automorphisme d'alg\`ebre unitaire \ $\theta_{\delta}$ \ de \ $\A$ \ envoyant \ $a$ \ sur \ $a + \delta.b$ \ et \ $b$ \ sur \ $b$. On peut donc voir \ $E \otimes_{a,b} E_{\delta}$ \ comme le \ $\A-$module obtenu en faisant agir \ $\A$ \ sur \ $E$ \ via \ $(\alpha,x) \mapsto \theta_{\delta}(\alpha).x$.

\subsection{Exposants.}

  \begin{defn}\label{Exp(E)}
  Soit \ $E$ \  un (a,b)-module r\'egulier. On notera \ $Exp(E) \subset \mathbb{C}\big/ \mathbb{Z}$ \  l'ensemble des classes modulo \ $\mathbb{Z}$ \  des nombres \ $-\alpha $ \ o\`u \ $\alpha$ \ d\'ecrit l'ensemble des  racines du polyn\^ome de Bernstein \ $B_E$ \ de \ $E$.
  \end{defn}
  
  On a donc toujours  \ $Exp(E) = Exp(E^{\sharp})$, puisque, par d\'efinition \ $B_E = B_{E^{\sharp}}$.
  
  \parag{Remarques} 
  \begin{enumerate}[1)]
  \item On notera que \ $[\lambda]$ \ est dans \ $Exp(E)$ \ si et seulement s'il existe \ $\lambda \in [\lambda]$ \ et une injection (a,b)-lin\'eaire de \ $E_{\lambda}$ \ dans \ $E$. En effet, il suffit de prouver cette assertion pour \ $E^{\sharp}$, et dans ce cas on peut prendre pour \ $\lambda$ \ le plus petit \'el\'ement de \ $[\lambda]$ \ pour lequel \ $a - \lambda.b$ \ n'est pas injectif, d'apr\`es la proposition 1.3  de [B.93]. 
  \item Soit \ $E^*$ \ le dual du (a,b)-module r\'egulier \ $E$. Alors 
  \ $ Exp(E^*) = - Exp(E) $.
  \item En utilisant la remarque pr\'ec\'edente et l'isomorphisme de \ $E$ \ avec son bi-dual, on voit que \ $[\lambda]$ \ est dans \ $Exp(E)$ \ si et seulement s'il existe \ $\lambda \in [\lambda]$ \ et une surjection (a,b)-lin\'eaire de \ $E$ \ dans \ $E_{\lambda}$. $\hfill \square$
  \end{enumerate}

   \begin{lemma}\label{Coprimitif 1}
  Soit \ $0 \to F \to E \overset{\pi}{\to} G \to 0 $ \ une suite exacte de (a,b)-modules r\'eguliers. Alors on a l'\'egalit\'e \ $Exp(E) = Exp(F) \cup Exp(G)$.
    \end{lemma}
  
  \parag{Preuve} Soit  \ $[\lambda] \in Exp(E)$. Alors il existe \ $\lambda \in [\lambda]$ \ et  une injection (a,b)-lin\'eaire \ $E_{\lambda} \hookrightarrow E$. Si on \ $\pi(E_{\lambda}) = \{0\}$ \ alors on a \ $E_{\lambda} \subset F$ \ et \ $[\lambda] \in Exp(F)$. Sinon, on a \ $\pi(E_{\lambda}) \simeq E_{\lambda+p}$ \ avec \ $p \in \mathbb{N}$, et on a donc \ $[\lambda] \in Exp(G)$.\\
  R\'eciproquement montrons que \ $Exp(G) \subset Exp(E)$ \ puisque l'inclusion de \ $Exp(F)$ \ dans \ $ Exp(E)$ \ est claire.\\
  Soit  \ $ f :  G \to E_{\mu}$ \ une application surjective. La compos\'ee \ $f \circ \pi $ \ est surjective ce qui montre que \ $[\mu] \in Exp(E)$, gr\^ace \`a la remarque 3) ci-dessus. $\hfill \blacksquare$
  
  \parag{Remarque} Une cons\'equence facile de ce lemme est que si on a deux sous-(a,b)-modules \ $F$ \ et \ $F'$ \ d'un (a,b)-module r\'egulier \ $E$ \ et si l'on a
    $$[\lambda] \not\in Exp(F) \cup Exp(F')$$
   alors \ $[\lambda]$ \ n'est pas dans \ $Exp(G)$ \ o\`u \ $G$ \ d\'esigne le plus petit sous-module normal de \ $E$ \ contenant \ $F + F'$.$\hfill \square$
   
   \subsection{Les  (a,b)-modules $[\p]-$primitifs.}

\bigskip

 \begin{defn} \label{primitif}
 Soit \ $\p \subset \mathbb{C}\big/\mathbb{Z}$. 
 On dira qu'un (a,b)-module r\'egulier \ $E$ \ est \ $[\p]-$primitif si toutes les racines du polyn\^ome de Bernstein de \ $E$ \ sont dans \ $- \p$, c'est-\`a-dire si \ $Exp(E) \subset \p$. 
 \end{defn}
 
 \parag{Notations} 
 Quand \ $\p = \{[\lambda]\}$ \ nous dirons que \ $E$ \ est \ $[\lambda]-$primitif. \\
 Si \ $M$ \ est le compl\'ementaire de \ $\p$ \ dans \ $\mathbb{C}\big/\mathbb{Z}$ \ nous dirons que \ $E$ \ est \ $[\not=\p]-$primitif pour dire qu'il est \ $M-$primitif.
 
\parag{Remarques}
\begin{enumerate}
\item  Avec notre d\'efinition le (a,b)-module nul est \ $\p-$primitif pour tout choix de \ $\p$. De plus c'est le seul (a,b)-module qui soit \`a la fois \ $\p-$primitif \ et \ $[\not=\p]-$primitif.
\item Une cons\'equence imm\'ediate de la  remarque 1)  qui suit  la d\'efinition \ref{Exp(E)}  est que tout sous-(a,b)-module (normal ou non)  d'un (a,b)-module  \ $\p-$primitif est  \ $\p-$primitif. 
\item Si on a une suite exacte \ $0 \to F \to E\to G\to 0$ \ avec \ $F$ \ et \ $G \quad \p-$primitifs, alors \ $E$ \ est \'egalement \ $\p-$primitif.\\
 Et r\'eciproquement si \ $E$ \ est \ $\p-$primitif dans une suite exacte de (a,b)-modules, alors \ $F$ \ et \ $G$ \ le sont \'egalement.
\item Si \ $f : E \to F$ \ est une application \ $\A-$lin\'eaire  entre  (a,b)-modules r\'eguliers et si \ $G \subset E$ \ est un sous-module \ $\p-$primitif, alors \ $f(G)$ \ est \ $\p-$primitif. En effet, sinon on peut trouver un sous-module isomorphe \`a \ $E_{\mu}$ \ dans \ $f(G)$ \ avec \ $\mu \not\in \p$ \ et donc un sous-module \ $H : = G \cap f^{-1}(E_{\mu})$ \ avec  une suite exacte
$$ 0 \to Ker(f)\cap G \to H \overset{f}{\to} E_{\mu} \to 0  $$
 ce qui contredit la remarque 3) pr\'ec\'edente.$\hfill \square$\\
\end{enumerate}
  
 \begin{prop} \label{partie primitive}
 Soit \ $E$ \ un (a,b)-module r\'egulier et soit \ $\p \subset  \mathbb{C}\big/ \mathbb{Z}$. Il existe un unique sous-(a,b)-module normal \ $E[\p] \subset E$ \ qui est \ $\p-$primitif et contient tout sous-module  \ $\p-$primitif de \ $E$.
 \end{prop}
 
 \parag{Preuve} Montrons l'assertion par r\'ecurrence sur le rang de \ $E$. Comme l'assertion est claire en rang 1, supposons l'assertion montr\'ee en rang \ $\leq k-1$ \ avec \ $k \geq 2$ \ et montrons-la en rang \ $k$.\\
 Si tout \ $[\lambda] \in -\p$ \ n'est pas la classe modulo \ $\mathbb{Z}$ \ d'une racine du polyn\^ome de Bernstein de \ $E$, il est clair que \ $\{0\}$ \ est le plus grand sous-(a,b)-module \ $\p-$primitif de \ $E$.\\
 Supposons donc qu'il existe une racine \ $-\lambda$ \ du polyn\^ome de Bernstein de \ $E$ \ telle que \ $\lambda \in \p$. On peut alors trouver un sous-module normal de \ $E$ \ isomorphe \`a \ $E_{\lambda'}$ \ avec \ $[\lambda'] = [\lambda]$ \ en normalisant l'image d'une injection \ $E_{\lambda} \hookrightarrow E$. On a alors la suite exacte
 $$ 0 \to E_{\lambda'} \to E \overset{\pi}{\to}  F \to 0 $$
 et \ $F$ \ est de rang \ $k-1$. L'hypoth\`ese de r\'ecurrence nous fournit un plus grand sous-module \ $\p-$primitif \ $F[\p]$ \ dans \ $F$ \ qui est normal. Montrons qu'alors \ $\pi^{-1}(F[\p])$ \ est le \ $E[\p]$ \ cherch\'e. D'abord il est normal dans \ $E$ \ puisque \ $F[\p]$ \ l'est dans \ $F$. De plus la suite exacte
 $$ 0 \to E_{\lambda'} \to \pi^{-1}(F[\p]) \overset{\pi}{\to}  F[\p] \to 0 $$
 montre qu'il est \ $\p-$primitif d'apr\`es la remarque 3) ci-dessus.\\
  Soit \ $G$ \ un sous-module   \ $\p-$primitif  de \ $E$. D'apr\`es la remarque 4)  faite plus haut son image par \ $\pi$ \ est  \ $\p-$primitive donc contenue dans \ $F[\p]$, ce qui montre que \ $G$ \ est bien contenu dans \ $\pi^{-1}(F[\p])$. $\hfill \blacksquare$
  
 \parag{Remarques}
 \begin{enumerate}
 \item On a \ $Exp(E[\p]) = Exp(E) \cap \p $.
 \item  Soit \ $E$ \ un (a,b)-module r\'egulier et \ $\p \subset \mathbb{C}\big/\mathbb{Z}$.  Pour \ $F \subset E$ \ un sous-(a,b)-module on a \ $F[\p] = F \cap E[\p]$. En effet l'inclusion de  \ $F \cap E[\p]$ \ dans \ $F[\p]$ \ r\'esulte de la maximalit\'e de \ $F[\p]$ \ puisque \ $F \cap E[\p]$ \ est \ $\p-$primitif. L'autre inclusion est \'evidente. $\hfill \square$
 \end{enumerate}  
  
  \begin{lemma} \label{coprimitif}
  Soit \ $E$ \ un (a,b)-module r\'egulier et \ $\p \subset \mathbb{C}\big/\mathbb{Z}$. D\'efinissons maintenant  \ $E\big/[\p] : = E\big/E[\not=\p]$. Alors  \ $E\big/[\p] $ \ est \ $\p-$primitif  et tout (a,b)-module quotient \ $\p-$primitif de \ $E$ \ est canoniquement un quotient de \ $E\big/[\p] $.
  \end{lemma}
  
   \parag{Preuve} Soit \ $\mu \in \p$ \ et consid\'erons une surjection \ $E \overset{\pi}{\to} E_{\mu}$. La restriction de \ $\pi$ \ \`a \ $E[\not= \p]$ \ est soit nulle, soit d'image \ $E_{\mu+p}$ \ pour un \ $p \in \mathbb{N}$. Mais ce second cas est exclu car il impliquerait que \ $[\mu] \in Exp(E[\not=\p])$, contredisant la remarque 1) ci-dessus.  On a donc \ $Exp(E[\not=\p]) = Exp(E)\setminus \p$. \\
   Supposons maintenant que \ $F \subset E$ \ est un sous-(a,b)-module normal tel que \ $E\big/F$ \ soit \ $\p-$primitif. Comme on a \ $F[\not=\p] = F \cap E[\not=\p]$ \ d'apr\`es la remarque 2) pr\'ec\'edente, on aura une injection de \ $E[\not=\p]\big/F[\not=\p]$ \ dans \ $E\big/F$ \ qui est suppos\'e \ $\p-$primitif. On en d\'eduit que \ $E[\not=\p]\big/F[\not=\p]$ \ est nul d'apr\`es les remarques 1 et 4 qui suivent la d\'efinition \ref{primitif}. Donc \ $E[\not=\p]\ \subset F$ \ et \ $E\big/F$ \ est un quotient de \ $E\big/[\p] $. $\hfill \blacksquare$
  
  \begin{defn}
  Nous appellerons partie \ $\p-$coprimitive de \ $E$ \ le quotient \ $E\big/[\p]$ \ introduit au lemme pr\'ec\'edent.
  \end{defn}

  \parag{Remarques}
  \begin{enumerate}
  \item Soit \ $ 0 \to F \overset{f}{\to} E \overset{g}{\to} G \to 0 $ \ une suite exacte de (a,b)-modules r\'eguliers. Pour tout \ $\p \subset \mathbb{C}\big/\mathbb{Z}$ \ on a la suite exacte :
  $$ 0 \to F[\p]  \overset{f}{\to} E[\p] \overset{g}{\to} G[\p]  .$$
  On n'a pas exactitude \`a droite en g\'en\'eral, ce que l'on peut d\'ej\`a v\'erifier sur la suite exacte\footnote{Rappelons que \ $E_{\lambda, \mu}$ \ est le (a,b)-module de rang 2  o\`u \ $a$ \ est d\'efini par  $$ a.e_1 = e_2 + (\lambda-1).b.e_1 \quad a.e_2 = \mu.b.e_2.$$}
  $$ 0 \to E_{\mu} \to E_{\lambda,\mu} \to E_{\lambda-1} \to 0  $$
  avec \ $\p = \{\lambda\}$ \ en supposant que \ $\mu \not\in [\lambda]$.
  
  \item  Soit \ $E$ \ un (a,b)-module r\'egulier et \ $\p \subset \mathbb{C}\big/\mathbb{Z}$. On a une fl\`eche naturelle
  $$ E[\p] \to  E\big/[\p] $$
  donn\'ee par composition de l'inclusion \ $E[\p] \hookrightarrow E $ \ et du quotient \ $ E \to E\big/[\p] $. Cette fl\`eche est injective, mais pas surjective en g\'en\'eral. En effet l'injectivit\'e r\'esulte de l'\'egalit\'e \ $E[\p] \cap E[\not=\p] = \{0\}$.  Elle n'est pas surjective d\'ej\`a pour \ $E_{\lambda,\mu}$ \ si \ $ [\lambda] \not= [\mu]$ \ et \ $\p = [\lambda]$ \ puisque \ $E[\p] = E_{\lambda}$ \ et \ $E\big/[\p] = E_{\lambda-1}$ \ dans ce cas.   En g\'en\'eral, on a donc \ $E[\p] \oplus E[\not=\p] \not=  E$. $\hfill \square$
  \end{enumerate}
  
  \bigskip
  
  \begin{cor} \label{decomp. primitive}
  Soit \ $E$ \ un (a,b)-module r\'egulier, et soient \ $Exp(E) : = \{\lambda_1, \dots, \lambda_d\}$, avec \ $[\lambda_i] \not= [\lambda_j]$ \ pour \ $i \not= j$, et rang\'es dans un ordre arbitraire. Alors on a une suite de composition unique (une fois l'ordre des \ $[\lambda_j]$ \ fix\'e)
  $$ 0 = F_0 \subset F_1 \subset \cdots \subset F_d = E$$
  de sous-modules normaux de \ $E$ \ tels que \ $F_{j}\big/ F_{j-1}$ \ soient \ $[\lambda_j]-$primitifs pour chaque \ $j \in [1,d]$.
  \end{cor}
  
 \parag{Preuve} La r\'ecurrence est imm\'ediate en posant  \ $  F_j : = E[\{\lambda_1, \dots, \lambda_j\}] $. $\hfill  \blacksquare$
 
 \bigskip
  
  On prendra garde que, en g\'en\'eral, pour \ $j \ge2$, le quotient \ $F_j\big/F_{j-1}$ \ n'est pas isomorphe \`a \ $E[\lambda_j]$, comme le montre l'exemple du (a,b)-module de rang 2 \ $E_{\lambda,\mu}$ \ quand \ $[\lambda] \not= [\mu]$.
  
 \bigskip

    \begin{cor}\label{Coprimitif 2}
  Soit \ $E$ \ un (a,b)-module r\'egulier et soit \ $\p \subset \mathbb{C}\big/\mathbb{Z}$. La dualit\'e des (a,b)-modules transforme la suite exacte
  $$ 0 \to E[\p] \to E \to E\big/E[\p] \to 0 $$
  en la suite exacte
  $$ 0 \to E^*[\not=-\p] \to E^* \to (E[\p])^* \to 0 $$
  ce qui montre que l'on a un isomorphisme canonique \ $(E[\p])^* \simeq E^*\big/[-\p]$.
  \end{cor}
  
  \parag{Preuve} Le dual d'un (a,b)-module \ $\p-$primitif est \ $[-\p]-$primitif. Donc la propri\'et\'e universelle de l'inclusion  \ $E[\p] \hookrightarrow E$ \ vis \`a vis des applications \ $\A-$lin\'eaires dans \ $E$ \ de modules \ $\p-$primitifs donne par dualit\'e que la surjection \ $E^* \to (E[\p])^*$ \ factorise toute application \ $\A-$lin\'eaire de \ $E^*$ \ dans un (a,b)-module \ $[-\p]-$primitif. Donc \ $(E[\p])^*$ \ est la partie \ $\p-$coprimitive de \ $E^*$. Ceci montre que le noyau de ce quotient est la partie primitive de \ $E^*$ \ pour \ $[\not=-\p]$. $\hfill \blacksquare$
  
  \bigskip
  
  Une cons\'equence simple de ce corollaire, puisque le dual d'un (a,b)-module monog\`ene r\'egulier est  monog\`ene r\'egulier (voir [B.09]), est que la partie \ $\p-$primitive d'un (a,b)-module monog\`ene r\'egulier est encore un (a,b)-module monog\`ene r\'egulier. \\
  En effet l'aspect monog\`ene est clair pour la partie coprimitive,  le corollaire ci-dessus donne alors  cette assertion par dualit\'e.
  
  \section{Th\`emes.}
  
   \subsection{D\'efinition et stabilit\'e par quotient et dualit\'e.}
   
   \subsubsection{D\'efinition et exemples.}

  \parag{Notations} Soit \ $\Lambda \subset ]0,1] \cap \mathbb{Q}$ \ un sous-ensemble fini et \ $N$ \ un entier. Nous consid\`ererons le \ $\C[[b]]-$module libre de type fini
\begin{equation*}
  \Xi_{\Lambda}^{(N)} : = \sum_{\lambda  \in \Lambda, j \in [0,N]} \quad \mathbb{C}[[b]].
  s^{\lambda -1}.\frac{(Logs)^j}{j!} \tag{@}
  \end{equation*}
  muni de la structure de \ $\A-$module (\`a gauche) donn\'ee par l'application \ $\C-$lin\'eaire \ $a$ \ qui est  la multiplication par \ $s$, la notation des g\'en\'erateurs correspondant au fait que l'on interpr\`ete \ $b$ \ comme l'int\'egration sans constante. Ceci correspond \`a l'\'egalit\'e
 $$\Xi_{\Lambda}^{(N)} =  \sum_{\lambda  \in \Lambda, j \in [0,N]} \quad \mathbb{C}[[s]].
  s^{\lambda -1}.\frac{(Logs)^j}{j!}  .$$

  Nous noterons ausssi \ $\Xi$ \ la somme de de tous les \ $\Xi_{\lambda}^{(N)}$ \ pour tous les \ $\lambda \in ]0,1] \cap \mathbb{Q}$ \ et tous les entiers \ $N$.\\

 \begin{defn} \label{theme}
 Nous dirons qu'un (a,b)-module monog\`ene est un {\bf th\`eme} s'il est isomorphe \`a un sous-(a,b)-module ( n\'ecessairement monog\`ene)  de  \ $\Xi_{\Lambda}^{(N)}$. 
   \end{defn}

 \parag{Remarques} \begin{enumerate}[1)]
 \item On notera d\'ej\`a qu'un th\`eme est toujours, par d\'efinition, un (a,b)-module monog\`ene {\bf g\'eom\'etrique}.
 \item  Le th\`eme \ $E$ \ est \ $[\lambda]-$primitif s'il  est isomorphe \`a un sous-module de 
 $$ \Xi_{\lambda}^{(N)} : = \sum_{j \in [0,N]} \   \mathbb{C}[[b]]. s^{\lambda -1}.\frac{(Logs)^j}{j!} $$
 pour \ $N \in \mathbb{N}$ \ assez grand,  o\`u \ $\lambda $ \ est dans \ $\mathbb{Q} \cap ]0,1]$.  $\hfill \square$
  \end{enumerate}
   
  \parag{Remarque} 
  En rang 1 tout (a,b)-module g\'eom\'etrique est un th\`eme, puisque pour \ $\lambda \in \mathbb{Q}^{+*}$ \ on a \ $E_{\lambda} \simeq \C[[b]].s^{\lambda-1} \subset \Xi $. $\hfill \square$

  \bigskip
  
  \begin{lemma}\label{rang 2}
  Un th\`eme \ $[\lambda]-$primitif de rang 2 est isomorphe soit  \`a \ $E_{\lambda,\lambda}$  \ soit \`a $E_{\lambda+n,\lambda}(\alpha)$ \ avec \ $\lambda \in 1 +\mathbb{Q}^{*+}$, $n \in \mathbb{N}^*$ \ et \ $\alpha \in \mathbb{C}^*$, c'est \`a dire isomorphe soit \`a \ $\A\big/\A.(a-\lambda.b).(a - (\lambda-1).b)$ \ soit  \`a \ $\A\big/\A.P_{n,\alpha}$ \ avec 
   $$P_{n,\alpha} : = (a - \lambda.b).(1 + \alpha.b^n)^{-1}.(a-(\lambda+n-1).b) .$$
  En particulier il contient un unique sous-module normal de rang 1 qui est isomorphe \`a \ $E_{\lambda}$.
  \end{lemma}
    
  \parag{Preuve} Dans la classification des (a,b)-modules r\'eguliers de rang 2 donn\'ee dans la proposition 2.4 de [B.93] p.34 on constate que les deux premiers types ne sont pas monog\`enes (car ils sont \`a p\^ole simple). Le quatri\`eme type est le second cas  donn\'e dans l'\'enonc\'e. Il se plonge dans \ $\Xi_{\lambda}^{(1)}$ \ pour \ $\lambda \in 1+\mathbb{Q}^{*+}$ \ sous la forme \ $\A.\psi$ \ o\`u 
  $$ \psi : = s^{\lambda+n-2}.Logs + \gamma.s^{\lambda-2}  $$
  avec 
   $$\gamma = -\frac{(\lambda-1).\lambda \dots (\lambda+n-2)}{n}.$$
  
  \smallskip
  
  Montrons qu'un (a,b)-module de rang deux primitif du troisi\`eme type n'est un th\`eme que dans le cas \ $E_{\lambda,\lambda}$. Montrons  donc que  \ $E : = E_{\lambda,\lambda+p}$ \ avec \ $\lambda \in \mathbb{Q}^{*+}$ \ et \ $p \in \mathbb{N}^*$ \ n'est pas un th\`eme. Raisonnons par l'absurde : si c'\'etait le cas, on aurait un plongement de \ $E_{\lambda,\lambda+p}$ \ dans \ $\Xi_{\lambda}^{(N)}$ \ et donc si l'on consid\`ere la \ $\mathbb{C}[[b]]-$base standard de \ $E_{\lambda,\lambda+p}$ \ donn\'ee par \ $ a.e_1 = e_2 + (\lambda-1).b.e_1 $ \ et \ $a.e_2 = (\lambda+p).b.e_2$, l'image de \ $e_2$ \ par ce plongement serait \'egale \`a \ $c.s^{\lambda+p-1}$ \ avec \ $c \in \mathbb{C}^*$. Soit \ $F$ \ l'image de \ $b.e_1$ \ par ce plongement. On aura alors
  $$ s.\frac{dF}{ds} = c.s^{\lambda+p-1} + (\lambda-1).F $$
  et la resolution de cette \'equation diff\'erentielle donne \ $F(s) = \frac{c}{p}.s^{\lambda+p-1} + \gamma.s^{\lambda-1} $. Pour \ $\lambda \in ]0,1]$ \ ceci impose \ $\gamma = 0$, puisque \ $F \in  b\Xi_{\lambda}^{(N)}$. On en d\'eduit que l'image dans \ $\Xi$ \ de \ $E_{\lambda,\lambda+p}$ \ est de rang 1  et \'egale \`a \ $\mathbb{C}[[b]].s^{\lambda-2}$ \ pour \ $\lambda >1$ \ et \`a \ $\mathbb{C}[[b]].s^{\lambda+p-1}$ \ pour \ $\lambda \in ]0,1]$. Ceci montre notre assertion.\\
  Par ailleurs on v\'erifie facilement que \ $\A.\varphi$ \ avec \ $\varphi : = s^{\lambda-2}.Logs $ \ est bien un plongement de \ $E_{\lambda,\lambda}$ \ dans \ $\Xi$ \ pour \ $\lambda \in  1+\mathbb{Q}^{*+}$.\\
 L'unicit\'e du sous-module normal de rang 1 pour les deux cas consid\'er\'es s'obtient facilement par un calcul direct ; le lecteur pourra aussi  se reporte \`a  [B.93]. $\hfill \blacksquare$
  
  \parag{Remarque} On notera que dans tous les cas la suite de Jordan-H{\"o}lder d'un th\`eme de rang 2 v\'erifie l'in\'egalit\'e  \ $\lambda_1 \leq \lambda_2-1$. On a m\^eme \ $\lambda_1 \leq \lambda_2$ \ sauf dans le cas de \ $E_{\lambda,\lambda}$. \ $\hfill \square$

   \bigskip
   
   \subsubsection{Quotient et dual d'un th\`eme.}
   
   La proposition suivante est la clef du th\'eor\`eme de stabilit\'e des th\`emes par quotient.
  
  \begin{prop}\label{JH unique}
  Soit \ $E$ \ un th\`eme \ $[\lambda]$ \ primitif non nul. Alors \ $E$ \ admet un unique sous-module normal de rang 1. Si \ $E_{\lambda} \subset E$ \ est ce sous-module normal, le quotient \ $E\big/ E_{\lambda}$ \ est un th\`eme \ $[\lambda]-$primitif.
  \end{prop}
  
  \parag{Preuve} L'existence des suites de Jordan-H{\"o}lder pour les (a,b)-modules r\'eguliers montre qu'il existe au moins un sous-(a,b)-module normal de rang 1 dans \ $E$, et comme \ $E$ \ est \ $[\lambda]-$primitif, il est isomorphe \`a \ $E_{\lambda_1}$ \ o\`u \ $\lambda_1 \in [\lambda]$.\\
  Supposons que l'on dispose de deux sous-modules normaux, not\'es respectivement \ $G_1 \simeq E_{\lambda_1}$ \ et \ $G_2 \simeq E_{\lambda_2}$. Posons \ $G : = G_1 + G_2$, et montrons que \ $G$ \ est n\'ecessairement de rang 2 si l'on suppose \ $G_1 \not= G_2$.\\
  En effet si \ $G$ \ est de rang 1 il est isomorphe \`a \ $E_{\lambda}$ \ et on a n\'ecessairement  \ $G_1 = b^p.G$ \ et \ $G_2 = b^q.G$. Mais la normalit\'e de \ $G_1$ \ et \ $G_2$ \ donne \ $p = q = 0$, c'est \`a dire \ $G = G_1 = G_2$.\\
  Donc \ $G$ \ est un th\`eme \ $[\lambda]-$primitif  de rang 2. Mais on a vu qu'un th\`eme \ $[\lambda]-$primitif  de rang 2 n'admet qu'un unique sous-module normal de rang 1. Donc on a \ $G_1 = G_2$, puisque la normalit\'e de \ $G_i$ \ dans \ $E$ \ implique sa normalit\'e dans \ $G$;  ceci prouve l'unicit\'e.
  
  \smallskip
  
  Pour montrer que le quotient \ $E\big/E_{\lambda}$ \ est un th\`eme, commen{\c c}ons par montrer que l'on a pour chaque \ $\lambda \in ]0,1] \cap  \mathbb{Q}$ \ une suite exacte de \ $\A-$modules \`a gauche
  $$ 0 \to \mathbb{C}[[b]].s^{\lambda-1} \to \Xi_{\lambda}^{(N)} \overset{f_{\lambda}}{\longrightarrow} \Xi_{\lambda}^{(N-1)} \to 0  $$
  o\`u \ $N \in \mathbb{N}^*$ \ et o\`u l'on rappelle que 
    $$ \Xi_{\lambda}^{(N)} : = \sum_{j \in [0,N]} \quad \mathbb{C}[[b]].s^{\lambda-1}.\frac{(Logs)^j}{j!} .$$
  D\'efinissons l'application \ $\mathbb{C}[[b]]-$lin\'eaire \ $f_{\lambda}$ \ en posant 
  \begin{align*}
  & f_{\lambda}(s^{\lambda-1}) = 0 \quad \quad  {\rm et } \\
  &   f_{\lambda}(s^{\lambda-1}.\frac{(Logs)^j}{j!}) = s^{\lambda-1}.\frac{(Logs)^{j-1}}{(j-1)!} \quad {\rm pour} \ j \geq 1 .
  \end{align*}
  On v\'erifie alors  que l'on a \ $f_{\lambda}(a.s^{\lambda-1}.\frac{(Logs)^j}{j!}) = a.f(s^{\lambda-1}.\frac{(Logs)^j}{j!})$ \ en utilisant la \ $\mathbb{C}[[b]]-$lin\'earit\'e de \ $f_{\lambda}$ \ et les relations
  $$ a.s^{\lambda-1}.\frac{(Logs)^j}{j!} = \lambda.b.s^{\lambda-1}.\frac{(Logs)^j}{j!} + b(s^{\lambda-1}.\frac{(Logs)^{j-1}}{(j-1)!}) $$
  pour \ $j \geq 1$ \ et \ $a.s^{\lambda-1} = \lambda.b.s^{\lambda-1}$ \ pour \ $j = 0 $. Le fait que \ $f_{\lambda}$ \ soit surjective et de noyau \ $\mathbb{C}[[b]].s^{\lambda-1} $ \ est alors imm\'ediat.
  
  \smallskip
  
  Consid\'erons alors un th\`eme \ $[\lambda]-$primitif \ $E \hookrightarrow \Xi_{\lambda}^{(N)}$, et soit \ $F$ \ son unique sous-module normal de rang 1. Il s'envoie bijectivement  sur \ $\mathbb{C}[[b]].s^{\lambda+p-1}$ \ pour un entier \ $p \geq 0$, par l'injection de \ $E$ dans \ $\Xi_{\lambda}^{(N)}$. Montrons qu'il est \'egal \`a \ $E \cap Ker(f_{\lambda})$. \\
  En effet il est contenu dans \ $Ker(f_{\lambda})$ \ d'apr\`es ce qui pr\'ec\`ede, et si \ $x \in Ker(f_{\lambda}) \cap E$, on a \ $ x = S(b).s^{\lambda-1}$. Notons \ $q$ \ la valuation de \ $S\in \mathbb{C}[[b]] \setminus \{0\}$ \ (le cas \ $x = 0$ \ est clair). Si on a \ $q < p$, alors \ $s^{\lambda+q-1} \in E$, puisque \ $E$ \ est un \ $\mathbb{C}[[b]]-$sous-module, et on obtient ainsi un \'el\'ement \ $y$ \ de \ $E$ \ tel que \ $b^{p-q}.y \in F$. Comme \ $F$ \ est normal, on a \ $y \in F$, ce qui est contredit l'hypoth\`ese \ $q < p$. Donc \ $S(b) \in b^p.\mathbb{C}[[b]]$ \ et \ $x \in F$.\\
  Donc le noyau de \ $f_{\lambda}$ \ restreinte \`a \ $E$ \ est \ $F$ \ et donc \ $f_{\lambda}$ \ induit une injection (a,b)-lin\'eaire de \ $E\big/F$ \ dans \ $\Xi_{\lambda}^{(N-1)}$. Donc \ $E\big/F$ \ est un th\`eme \ $[\lambda]-$primitif. $\hfill \blacksquare$

   \bigskip
  
  \begin{thm}\label{Stabilite par quotient}
  Soit  \ $E$ \ un th\`eme et \ $F$ \ un sous-module (a,b)-module monog\`ene de \ $E$;  alors \ $F$ \ est un th\`eme. Si \ $F$ \ est  un sous-(a,b)-module normal dans \ $E$, alors  le quotient  \ $E\big/F$ \ est un th\`eme.
  \end{thm}
  
  \parag{Remarque} Si \ $F$ \ est  un sous-(a,b)-module normal dans \ $E$, alors \ $F$ \ est n\'ecessairement monog\`ene  c'est  donc un sous-th\`eme normal de \ $E$.\\
  En effet  si \ $F$ \ est normal, \ $F\big/b.F \to E\big/b.E$ \ est injective. Comme l'action de  \ $a$ \ sur \ $E\big/b.E$ \ est un donn\'ee par un nilpotent principal, le sous-espace stable \ $F\big/b.F$ \ est \'egal \`a \ $Im (a^h)$ \ pour un entier \ $h$;  donc \ $F\big/a.F + b.F$ \ est de dimension 1, ce qui implique que \ $F$ \ est monog\`ene. $\square$
  
  \parag{Preuve} La premi\`ere assertion est imm\'ediate. \\
  Comme le quotient d'un (a,b)-module monog\`ene est monog\`ene et le quotient d'un (a,b)-module g\'eom\'etrique est g\'eom\'etrique, le quotient \ $E\big/F$ \ est monog\`ene et g\'eom\'etrique. \\
  Montrons que \ $E\big/F$ \ est un th\`eme par r\'ecurrence sur le rang de \ $F$. \\
  En rang 1, le r\'esultat est une cons\'equence imm\'ediate de la preuve de la proposition \ref{JH unique} : en consid\'erant l'application \ $ g_{\lambda} : \Xi_{\Lambda}^{(N)} \to \Xi_{\Lambda}^{(N)} $ \ qui est donn\'ee par la somme directe de \ $f_{\lambda}$ \ compos\'ee avec l'inclusion \ $\Xi_{\lambda}^{(N-1)} \hookrightarrow \Xi_{\lambda}^{(N)}$,  sur \ $\Xi_{\lambda}^{(N)}$ \ et l'identit\'e sur \ $\Xi_{\mu}^{(N)}$ \ pour \ $\mu \in \Lambda, \mu \not= \lambda$.\\
  Supposons maintenant le r\'esultat \'etabli pour \ $F$ \ de rang \ $\leq k-1$ \ et consid\'erons un sous-module normal \ $F$ \ de rang \ $k$ \ d'un th\`eme \ $E$. En choisissant un sous-module normal \ $E_{\lambda} \subset F$ \ qui est normal dans \ $F$ \  donc dans \ $E$, on constate que \ $E\big/F \simeq (E\big/E_{\lambda})\Big/(F\big/E_{\lambda})$ \ et donc que l'on a un quotient du th\`eme \ $E\big/E_{\lambda}$ \ par le sous module normal \ $F\big/E_{\lambda}$ \ qui est de rang \ $k-1$. L'hypoth\`ese de r\'ecurrence permet donc de conclure. $\hfill \blacksquare$
  
   \bigskip
   
   \begin{cor}\label{caract. theme 1}
   Soit \ $E$ \ un (a,b)-module monog\`ene. C'est un th\`eme si et seulement si pour chaque \ $[\lambda] \in Exp(E)$, la partie $[\lambda]-$coprimitive \ $E\big/[\lambda]$ \ est un th\`eme.
   \end{cor}
   
   \parag{Preuve} Le th\'eor\`eme  \ref{Stabilite par quotient}  de stabilit\'e des th\`emes par quotient implique que la condition est n\'ecessaire.  Montrons qu'elle est suffisante.\\
   Posons \ $Exp(E) = \{ [\lambda_1], \dots, [\lambda_d]\}$. Soit \ $\theta_i : E \to \Xi_{\lambda_i}^{(N)}$ \ pour \ $i \in [1,d]$ \ l'application (a,b)-lin\'eaire obtenue en composant le quotient \ $E \to E\big/[\lambda_i]$ \ avec une injection (a,b)-lin\'eaire du th\`eme primitif \ $E\big/[\lambda_i]$ \ dans \ $\Xi_{\lambda_i}^{(N)}$, o\`u \ $\{\lambda_i\} = [\lambda_i] \cap ]0,1] $.\\
    Posons alors \ $\theta : = \oplus_{i=1}^d  \theta_i : E \to \Xi_{\Lambda}^{(N)}$, o\`u \ $\Lambda : = \{\lambda_1, \dots, \lambda_d\}$, et montrons que \ $\theta $ \ est injective. Par construction, on a \ $Ker\, \theta_i = E[\not=\lambda_i]$ \ et comme \ $\cap_i E[P_i]  = E[\cap_i P_i] $ \ on obtient l'injectivit\'e puisque \ $\cap_{i=1}^d [\not=\lambda_i] = \emptyset $ \ dans \ $Exp(E)$.  $\hfill \blacksquare$
    
    \bigskip
   
  Nous d\'eduirons  plus loin, gr\^ace au th\'eor\`eme de dualit\'e, qu'un (a,b)-module monog\`ene r\'egulier  \ $E$ \ est un th\`eme si et seulement si pour chaque \ $[\lambda] \in Exp(E)$ \ la partie primitive \ $E[\lambda]$ \ de \ $E$ \ est un th\`eme. Le lecteur, \`a titre d'exercice, pourra montrer directement ce r\'esultat en s'inspirant de la m\'ethode de d\'emonstration utilis\'ee pour obtenir la caract\'erisation suivante des th\`emes primitifs. 
  
  \bigskip
  
   \begin{thm}\label{Reciproque de l'unicite de JH} 
  Soit \ $E$ \ un (a,b)-module monog\`ene g\'eom\'etrique poss\'edant un unique sous-module normal de rang 1. Alors \ $E$ \ est un th\`eme primitif.
  \end{thm}
  
  \parag{D\'emonstration} Par r\'ecurrence sur le rang de \ $E$. L'assertion \'etant claire en rang 1, supposons-la d\'emontr\'ee en rang \ $k \geq 1$ \ et consid\'erons un (a,b)-module monog\`ene g\'eom\'etrique \ $E$ \ de rang \ $k+1$ \ v\'erifiant notre hypoth\`ese. Soit \ $F$ \ un sous-module normal de rang \ $k$ \ de \ $E$. Alors c'est un th\`eme \ $[\lambda]-$primitif d'apr\`es l'hypoth\`ese de r\'ecurrence\footnote{Remarquer que si  \ $G$ \ est normal dans \ $F$ \ qui est normal dans \ $E$, $G$ \ est normal dans \ $E$.}. On a une suite exacte
  \begin{equation*}
   0 \to F \to E \to  E_{\lambda'} \to 0  \quad {\rm avec} \quad \lambda' \in [\lambda]  \tag{@}
   \end{equation*}
  car si \ $E$ \ n'\'etait pas primitif, il poss\`ederait deux sous-modules normaux de rang 1 correspondant \`a des exposants distincts modulo \ $\mathbb{Z}$, contredisant l'hypoth\`ese.

  Fixons une injection \ $\A-$lin\'eaire \ $j : F \to \Xi$ \ et consid\'erons la suite exacte d'espaces vectoriels
  $$ 0 \to Hom_{\A}(E_{\lambda'}, \Xi) \to  Hom_{\A}(E, \Xi) \to  Hom_{\A}(F, \Xi) \to 0  $$
  d\'eduite de \ $(@)$, l'exactitude r\'esultant de [B. 05] th. 2.2.1 p.24. Soit \ $\tilde{j} \in Hom_{\A}(E, \Xi)$ \ s'envoyant sur \ $j \in Hom_{\A}(F, \Xi) $. Si \ $\tilde{j}$ \ est injective, la d\'emonstration est termin\'ee par d\'efinition d'un th\`eme. Sinon, soit \ $G : = Ker \tilde{j} \not= \{0\}$. Comme \ $\tilde{j}$ \ induit \ $j$ \ sur \ $F$, on aura \ $G \cap F = \{0\}$, et donc \ $G$ \ est de rang 1 et normal comme noyau, puisque \ $\Xi$ \ n'a pas de \ $b-$torsion. Mais l'unique sous-module normal de rang 1 de \ $E$ \ est contenu dans \ $F$, puisque  \ $F$ \ est de rang \ $k \geq 1$. Contradiction.\\
   Donc \ $\tilde{j}$ \ est injectif . $\hfill \blacksquare$
   
   \bigskip
   
  \begin{cor}\label{sous-themes normaux}
   Soit \ $E$ \ un th\`eme primitif de rang \ $k$. Alors \ $E$ \ poss\`ede pour chaque \ $j \in [0,k]$ \  un unique sous-module \ $F_j$ \  qui est normal et de rang \ $j$.\\
    Choisissons, pour chaque \ $j \in [0,k-1]$, une injection \ $\A-$lin\'eaire \ $\theta_j : E\big/F_j \to \Xi$. Alors \ $\theta_0, \dots, \theta_{k-1}$ \ forment une base de l'espace vectoriel \ $Hom_{\A}(E, \Xi)$.
   \end{cor}
   
   \parag{Preuve} Montrons par r\'ecurrence sur \ $j \geq 1$ \ l'unicit\'e du sous-(a,b)-module normal de rang \ $j$ \ dans un th\`eme primitif \ $E$. Comme le cas \ $j = 1$ \ a \'et\'e montr\'e au th\'eor\`eme pr\'ec\'edent, supposons  \ $j \geq 2$ \ et l'assertion montr\'ee pour un sous-(a,b)-module normal de rang \ $j-1$ \ d'un th\`eme primitif. Soit \ $E$ \ un th\`eme primitif et notons \ $G$ \ son unique sous-(a,b)-module normal de rang 1. Alors \ $E\big/G$ \ est un th\`eme primitif et il admet donc un unique sous-(a,b)-module normal \ $F_0$ \ de rang \ $j-1$. Soit \ $\pi : E \to E\big/G$ \ l'application quotient, et notons \ $F : = \pi^{-1}(F_0)$. Alors \ $F$ \ est un sous-(a,b)-module normal de rang \ $j$ \ de \ $E$. \\
   Consid\'erons alors un sous-(a,b)-module normal \ $F_1$ \ dans \ $E$ \ de rang \ $j$. On a \ $G \subset F_1$, car si \ $G$ \ n'\'etait pas l'unique sous-(a,b)-module normal de rang 1  de \ $F_1$, cela contredirait  l'unicit\'e de \ $G$. Donc \ $\pi(F_1)$ \ est de rang \ $j-1$ \ dans \ $E\big/G$.\\
    Il est normal car si \ $y \in \pi(F_1) \cap b.(E\big/G) $, on peut \'ecrire \ $y = \pi(x)$ \ o\`u \ $x \in F_1$ \ et \ $y = \pi(b.z)$ \ o\`u \ $z \in E$. Donc \ $x = b.z + t$ \ avec \ $t \in G \subset F_1$. Alors \ $x-t \in F_1 \cap b.E = b.F_1$. Donc \ $y = \pi(x-t) $ \ est dans \ $b.\pi(F_1)$. On en d\'eduit que \ $\pi(F_1) = F_0$ \ ce qui implique \ $F_1 = F$.
   
   La seconde assertion du corollaire r\'esulte du fait que \ $\dim_{\C} (Hom_{\A}(E, \Xi)) = k$ \ d'apr\`es le th\'eor\`eme 2.2.1 de  [B. 05]  et du fait que les  \ $\theta_i$ \ sont lin\'eairement ind\'ependants. En effet si l'on a \ $\sum_{i=0}^{k-1} \ \alpha_i.\theta_i = 0$ \ et \ $i_0$ \ est le premier indice pour lequel \ $\alpha_{i_0} \not=0 $, alors, pour \ $x \in F_{i_0 +1} \setminus F_{i_0}$ \ on obtient \ $ \alpha_{i_0}.\theta_{i_0}(x) = 0$, ce qui est absurde. $\hfill \blacksquare$
   
   \parag{\bf Remarque importante} Un th\`eme \ $[\lambda]-$primitif poss\`ede une {\bf unique suite de Jordan-H{\"o}lder}, et  r\'eciproquement, l'unicit\'e de la suite de Jordan-H{\"o}lder pour un (a,b)-module monog\`ene g\'eom\'etrique \ $[\lambda]-$primitifs {\bf caract\'erise les th\`emes  \ $[\lambda]-$primitifs}.\\
En fait, quitte \`a fixer un ordre dans \ $Exp(E)$, on a \'egalement unicit\'e de la suite de Jordan-H{\"o}lder d'un th\`eme g\'en\'eral respectant l'ordre fix\'e. $\hfill \square$
   
  \bigskip
     
   \begin{lemma}\label{suite}
    Soit \ $E$ \ un th\`eme \ $[\lambda]-$primitif. Soit \ 
    $$0 = F_0 \subset F_1 \subset \cdots \subset F_{k-1} \subset F_k = E$$
    son unique  suite de Jordan-H{\"o}lder et posons \ $F_j/F_{j-1} \simeq E_{\lambda_j}$ \ pour \ $j \in[1,k]$. Alors pour chaque \ $j \in [1,k-1]$ \  le nombre \ $p_j = \lambda_{j+1} - \lambda_j +1$ \ est un entier naturel.
    \end{lemma}
    
    \parag{Preuve} Le fait que \ $p_j \in \mathbb{Z}$ \ est trivial. C'est un entier positif ou nul en raison de la proposition 3.5.2 de [B.09] qui donne l'in\'egalit\'e \ $\lambda_{j+1} \geq \lambda_j - 1$. $\hfill \blacksquare$
    
    \bigskip
    
    On remarquera que les in\'egalit\'es du lemme pr\'ec\'edent reviennent \`a dire que la suite \ $\lambda_j + j$ \ est croissante (comparer avec  la proposition 3.5.2  de [B.09]).
    
    \parag{Notation}  Soit \ $E$ \ un th\`eme \ $[\lambda]-$primitif. Nous noterons \ $\lambda_1, \dots, \lambda_k$ \ les nombres associ\'es aux quotients de son unique suite de Jordan-H{\"o}lder. Une fa{\c c}on \'equivalente de se donner la suite (ordonn\'ee) \ $\lambda_1, \dots, \lambda_k$ \ consiste \`a pr\'eciser \ $\lambda_1$ \ et \`a se donner les entiers (positifs ou nuls) \ $p_1, \cdots, p_{k-1}$ \ d\'efinis en posant \ $\lambda_{j+1} = \lambda_j + p_j -1$ \ pour \ $j \in [1,k-1]$.$\hfill \square$
    
    \bigskip
    
     \begin{defn}\label{Inv. fond.}
    Nous appellerons {\bf invariants fondamentaux} d'un th\`eme \ $E$ \ suppos\'e  \ $[\lambda]-$primitif de rang \ $k$ \ la donn\'ee de {\bf  la suite ordonn\'ee  \ $\lambda_1, \dots, \lambda_k$} \ ou bien, ce qui est \'equivalent,  de \ $\lambda_1$ \ et des entiers \ $p_1, \cdots, p_{k-1}$.
    \end{defn}

   \parag{Remarque}  On notera que l'on a \ $\lambda_j + j > \lambda_k + k > k $ \ pour chaque \ $j \in [1,k]$ \ puisque l'on a \ $\lambda_k > 0$. En particulier on a \ $\lambda_1> k-1$. $\hfill \square$
     
   \bigskip
    
    Pour un th\`eme \ $[\lambda]-$primitif \ $E$ \ la donn\'ee des invariants fondamentaux est plus fine que la donn\'ee du polyn\^ome de Bernstein \ $B_E$ \ qui revient \`a se donner l'\'el\'ement de Bernstein\footnote{voir [B.09] def. 3.3.1.} \ $P_E : = (a - \lambda_1.b)\dots (a - \lambda_k.b) \in \A$. En effet, le polyn\^ome de Bernstein ne pr\'ecise pas l'ordre de ses racines.

   \bigskip
   
      \begin{thm}\label{dual d'un theme}
   Soit \ $E$ \ un th\`eme \ $[\lambda]-$primitif de rang \ $k$ \ d'invariants fondamentaux  \ $\lambda_1, \dots, \lambda_k$. Alors pour tout nombre rationnel \ $\delta$ \ v\'erifiant  \ $\delta > \lambda_k + k-1$ \ le (a,b)-module  \ $E^*\otimes_{a,b} E_{\delta}$ \  est un th\`eme \ $[\delta-\lambda]-$primitif, o\`u \ $E^*$ \ d\'esigne le dual de \ $E$, d'invariants fondamentaux \ $\delta - \lambda_k, \dots, \delta - \lambda_1$.
   \end{thm}
   
   \parag{D\'emonstration} 
   D'abord le dual d'un (a,b)-module r\'egulier et monog\`ene est r\'egulier et monog\`ene d'apr\`es [B.09]. Le dual d'un th\`eme \ $[\lambda]-$primitif \ est \ $[-\lambda]-$primitif.  Par ailleurs les sous-modules normaux du dual correspondent bijectivement aux duaux des quotients. Comme on a exactement un seul quotient pour chaque \ $j \in [0,k]$ \ o\`u \ $k$ \ d\'esigne le rang de \ $E$, on en conclut que \ $E^*\otimes_{a,b} E_{\delta}$ \  sera un th\`eme gr\^ace au th\'eor\`eme \ref{Reciproque de l'unicite de JH} d\`es qu'il sera g\'eom\'etrique, c'est \`a dire d\`es que  le premier quotient de la  suite de Jordan-H{\"o}lder   de \ $E^*\otimes_{a,b} E_{\delta}$ \ sera \ $> k-1$.  Comme il  vaut \ $ \delta  -\lambda_k$ \ l'assertion est d\'emontr\'ee.$\hfill \blacksquare$

 \parag{Remarques}
  \begin{enumerate}[1)]
 \item Donc pour un th\`eme \ $[\lambda]-$primitif \ $E$ \ d'invariants fondamentaux \ $\lambda_1, p_1, \dots, p_{k-1}$ \ et \ $\delta$ \ rationnel tel que \ $\delta - \lambda_k > k-1$, le th\`eme \ $[\delta-\lambda]-$primitif \ $E^*\otimes_{a,b} E_{\delta}$ \  aura pour invariants fondamentaux \ $ \delta - \lambda_k, p_{k-1}, \dots, p_1$.
 \item  Pour un th\`eme  g\'en\'eral on en d\'eduit que pour \ $\delta \in \mathbb{N}$ \ assez grand, $E^*\otimes E_{\delta}$ \ est un th\`eme, en remarquant  que la partie \ $[\lambda]-$coprimitive de \ $E^* \otimes E_{\delta}$ \  est \ $(E[-\lambda] )^*\otimes E_{\delta}$ \ pour chaque \ $[\lambda] \in \C\big/\mathbb{Z}$. On conclut alors gr\^ace au corollaire \ref{caract. theme 1}. $\hfill \square$
 \end{enumerate}
  
   \bigskip
   
    On a alors le dual du corollaire \ref{caract. theme 1}.
   
    \begin{cor}\label{caract. theme 2}
   Soit \ $E$ \ un (a,b)-module monog\`ene g\'eom\'etrique. C'est un th\`eme si et seulement si pour chaque \ $[\lambda] \in Exp(E)$, la partie $[\lambda]-$primitive \ $E[\lambda]$ \ de \ $E$ \ est un th\`eme.
   \end{cor}

     \subsection{Structure des th\`emes \ $[\lambda]-$primitifs.}
     
     \subsubsection{Le th\'eor\`eme de structure.}

    En fait le th\'eor\`eme 3.4.1 de [B.09]  donne le th\'eor\`eme de structure suivant pour les th\`emes \ $[\lambda]-$primitifs :

  \begin{thm}\label{pre-invariants}
  Soit \ $E$ \ un th\`eme \ $[\lambda]-$primitif dont les invariants fondamentaux sont  \ $\lambda_1, p_1, \dots, p_{k-1}$. Alors il existe \ $S_1, \cdots, S_{k-1}$ \ des \'el\'ements de \ $\C[b]$ \ v\'erifiant \ $S_j(0) = 1$ \ et \ $\deg(S_j) \leq p_j+ \cdots + p_{k-1}$ , tels que  l'on ait \ $E \simeq \A\big/ \A.P$ \ avec
  $$ P = (a - \lambda_1.b).S_1^{-1}(a -\lambda_2.b) \dots S_{k-1}^{-1}.(a - \lambda_k.b). $$
 De plus,  pour chaque \ $j \in [1,k-1]$ \ le coefficient de \ $b^{p_j}$ \ dans \ $S_j$ \ est non nul.\\
    R\'eciproquement, pour tout choix de \ $[\lambda] \in  \mathbb{Q}\big/\mathbb{Z}$, tout choix de \ $\lambda_1 \in [\lambda] , \lambda_1 >  k-1$, d'entiers \ $p_1, \cdots, p_{k-1}$ \ positifs ou nuls, et d'\'el\'ements \ $S_1, \cdots, S_{k-1}$ \ dans \ $\mathbb{C}[[b]]$ \ v\'erifiant \ $S_j(0) = 1$, tels que le coefficient de \ $b^{p_j}$ \ dans \ $S_j$ \ soit non nul,  le quotient \ $\A\big/\A.P$ \ est un th\`eme \ $[\lambda]-$primitif.
    \end{thm}

  \parag{D\'emonstration du th\'eor\`eme} La partie directe est une cons\'equence imm\'ediate du th\'eor\`eme 3.4.1 et du lemme  3.5.1 de [B.09], compte tenu de l'unicit\'e de la suite de Jordan-H{\"o}lder d'un th\`eme \ $[\lambda]-$primitif.\\
  Montrons la r\'eciproque. Il est clair que le quotient \ $\A\big/\A.P$ \ est un (a,b)-module monog\`ene g\'eom\'etrique de rang \ $k$. Nous allons montrer que c'est un th\`eme par r\'ecurrence sur \ $k$. Comme le cas \ $k = 1$ \ est \'evident, supposons le r\'esultat montr\'e pour \ $k-1$. Soit \ $Q : = (a - \lambda_1.b).S_1^{-1}(a -\lambda_2.b) \dots S_{k-2}^{-1}.(a - \lambda_{k-1}.b) $. L'hypoth\`ese de r\'ecurrence donne alors que \ $F : = \A\big/\A.Q$ \ est un th\`eme. Soit donc \ $\varphi \in \Xi_{\lambda}^{(N)}$ \ tel que \ $\A.\varphi$ \ soit isomorphe \`a \ $F$. Pour construire \ $\psi \in \Xi_{\lambda}^{(N+1)}$ \  v\'erifiant \ $(a- \lambda_k.b).\psi = S_{k-1}.\varphi$, il suffit de r\'esoudre une \'equation diff\'erentielle \'el\'ementaire. Explicitement \ $ s.f'(s) - \lambda_k.f(s) = S_{k-1}(b).\varphi(s)$, o\`u l'on a pos\'e \ $b.\psi = f$. Le point important est que, comme \ $F$ \ est de rang \ $k-1$, on peut, quitte \`a multiplier \ $\varphi$ \ par un inversible de \ $\C[[b]]$, supposer que  \ $\varphi$ \ est un polyn\^ome  en \ $Log\, s$ \ \`a coefficient dans \ $\mathbb{C}[[b]]$,  de degr\'e \ $k-2$ \ avec un coefficient de \ $(Log\, s)^{k-2}$ \ \'egal \`a \ $s^{\lambda_{k-1}-1}$. On constate alors que le fait que le coefficient de \ $b^{p_{k-1}}$ \ dans \ $S_{k-1}$ \ soit non nul, assure que le degr\'e en \ $Log\, s$ \ de \ $\psi$ \ sera exactement \ $k-1$, puisque \ $\lambda_k = \lambda_{k-1} + p_{k-1} -1$. Alors le morphisme \ $\A/\A.P \to \Xi_{\lambda}^{(N+1)}$ \ d\'efini en envoyant \ $1$ \ sur \ $\psi$ \ aura une image  \ $\A.\psi$ \ qui sera de rang $k$ \ sur \ $\mathbb{C}[[b]]$. C'est donc un isomorphisme, et \ $\A/\A.P$ \ est donc un th\`eme. $\hfill \blacksquare$
      
  \bigskip
  
  \subsubsection{Bases standards.}
  
  Commen{\c c}ons par expliciter la suite de Jordan-H{\"o}lder  d'un th\`eme \ $[\lambda]-$primitif  de rang \ $k$ \ quand il est plong\'e dans \ $\Xi_{\lambda}^{(k-1)}$ \  l'espace des d\'eveloppements asymptotiques.  
  
  \begin{lemma}\label{devissage 0}
   On fixe \ $[\lambda] \in \mathbb{Q}\big/\mathbb{Z}$ \ et on note \ $\{\lambda\} = ]0,1] \cap [\lambda]$. Soit \ $\varphi \in \Xi_{\lambda}^{(k-1)} $ \ tel que \ $E : = \A.\varphi$ \ soit un th\`eme de rang \ $k$. Alors pour chaque \ $j \in [1,k]$ \ l'intersection \ $E \cap \Xi_{\lambda}^{(j-1)}$ \ est l'unique sous-th\`eme normal \ $F_j$ \  de rang \ $j$ \ de \ $E$. La restriction de \ $\pi_j : \Xi_{\lambda}^{(j-1)} \to \Xi_{\lambda}^{(j-1)}\big/\Xi_{\lambda}^{(j-2)} \simeq E_{\lambda} $ \ \`a   \ $F_j$ \ a pour image \ $E_{\lambda_j} $ \ et pour noyau \ $F_{j-1}$.
  \end{lemma}
  
  \parag{Preuve} Montrons que le sous-(a,b)-module \ $G_j : = E \cap \Xi_{\lambda}^{(j-1)}$ \ est normal : si \ $x \in E$ \ et v\'erifie \ $b.x \in \Xi_{\lambda}^{(j-1)}$ \ on a n\'ecessairement \ $x \in \Xi_{\lambda}^{(j-1)}$ \ puisque \ $b$ \ pr\'eserve le degr\'e en \ $Log\, s$. \\
  Comme le noyau de \ $(\pi_j)_{\vert G_j}$ \ est \ $G_{j-1}$, on a \ $rg(G_j) \leq rg(G_{j-1}) + 1$, pour chaque \ $j \in [1,k]$. Comme \ $G_0 = \{0\}$ \ et \ $rg(G_k) = k$ \ par hypoth\`ese, on a n\'ecessairement \ $rg(G_j) = j$ \ pour tout \ $j \in [1,k]$ \ et donc \ $G_j = F_j$. $\hfill \blacksquare$
  
  \bigskip
  
  \begin{cor}\label{devissage 1}
  Dans la situation du lemme pr\'ec\'edent, posons \\
   $\pi_k(\varphi) = S_k.e_{\lambda_k}$. Alors \ $S_k$ \ est un inversible de \ $\mathbb{C}[[b]]$ \ et   l'\'el\'ement 
    $$\varphi_{k-1} : = (a - \lambda_kb)).S_k^{-1}.\varphi $$
    est un g\'en\'erateur du sous-th\`eme normal  \ $F_{k-1}$ \ de \ $E = F_k$.
  \end{cor}
  
  \parag{Preuve} En fait \ $S_k$ \ est le coefficient de \ $s^{\lambda_k -1}.(Log\, s)^{k-1}/k!$ \ dans \ $\varphi$, ce qui montre que \ $\varphi_{k-1} $ \ est dans \ $F_{k-1} = Ker(\pi_k)\cap E$. C'est n\'ecessairement un g\'en\'erateur de \ $F_{k-1}$ \ car dans l'espace vectoriel  \ $F_{k-1}\big/b.F_{k-1} \subset E\big/b.E$ \ qui est de dimension \ $k-1$, les classes \ $\varphi_{k-1}, a.\varphi_{k-1}, \cdots, a^{k-2}.\varphi_{k-1}$ \ forment un syst\`eme libre, puisque les classes de  \ $\varphi, a.\varphi, \cdots, a^{k-1}.\varphi$ \ dans \ $E\big/b.E$ \ forment une base. $\hfill \blacksquare$
  
  \parag{\bf La base standard de \ $\A/\A.P$} Soit \ $E : = \A/\A.P$ \ un th\`eme de rang \ $k$, o\`u l'on suppose que \ $P =  (a - \lambda_1.b).S_1^{-1}(a -\lambda_2.b) \dots S_{k-1}^{-1}.(a - \lambda_k.b)$. Soit  \ $e$ \ un g\'en\'erateur de \ $E$ \ comme \ $\A-$module, d'annulateur \ $\A.P$ \ (par exemple \ $e = [1]$). D\'efinissons les \'el\'ements \ $e_k, e_{k-1}, \cdots, e_1$ \ de \ $E$ \ par les relations suivantes :
  \begin{enumerate}[i)]
  \item \ $e_k : = e $ ;
  \item  \ $e_j : = S_j^{-1}.(a - \lambda_{j+1}.b).e_{j+1} $ \ pour \ $j \in [1,k-1]$.
  \end{enumerate}

 Une cons\'equence imm\'ediate du corollaire pr\'ec\'edent est que \ $e_1, \cdots, e_k$ \ est une \ $\C[[b]]-$base de \ $E$. Elle sera appel\'ee {\bf la base standard} associ\'ee au g\'en\'erateur  \ $e$ \ et au choix de \ $P$ \ engendrant l'annulateur de \ $e$ \ dans \ $E$.\\
  On notera que les relations \ $(a - \lambda_{j+1}.b).e_{j+1} = S_j.e_j , j \in [0,k-1]$ \ avec la convention \ $e_0 = 0$ \ d\'eterminent le (a,b)-module \ $E$ \ de \ rang \ $k$. $\hfill \square$
  
  \bigskip   
  
\begin{lemma}\label{lemme manquant}
Soit \ $E$ \ un th\`eme \ $[\lambda]-$primitif et soit, pour \ $j \in [0,k-1]$
 $$P_j : = (a - \lambda_{j+1}.b).S_{j+1}^{-1}\dots S_{k-1}^{-1}.(a - \lambda_k.b) $$
 o\`u \ $P_0$ \ est le g\'en\'erateur de l'id\'eal annulateur d'un g\'en\'erateur standard  \ $e : = e_k$ \ de \ $E$. Notons \ $e_1, \dots, e_k$ \ la base standard associ\'ee \`a \ $e$.\\
 Si \ $x \in F_{k-j}$ \ v\'erifie \ $P_0.x = 0$ \ dans \ $E$, alors \ $P_j.x = 0 $ \ dans \ $E$, et il existe \ $\rho \in \C$ \ tel que  \ $ x - \rho.b^{\lambda_k-\lambda_{k-j}}.e_{k-j}$ \ soit dans \ $F_{k-j-1}$.
 \end{lemma}
 
 \parag{Preuve} En envoyant \ $e$ \ sur \ $x$ \ on d\'efinit un \'el\'ement de \ $Hom_{\A}(E,F_{k-j})$ \ dont le noyau contient n\'ecessairement \ $F_{j}$. C'est donc une application \ $\A-$lin\'eaire de \ $E\big/F_{j}$ \ dans \ $F_{k-j}$. Comme \ $P_j(e) = S_{j}.e_{j}$ \ est dans \ $F_j$, on aura \ $P_j.x = 0$.\\
  Mais alors  l'image de \ $x$ \ dans \ $E_{\lambda_{k-j}}$ \ via le quotient \ $F_{k-j}\big/F_{k-j-1} \simeq E_{\lambda_{k-j}}$ \ est dans le noyau de 
 $$ P_j : E_{\lambda_{k-j}} \to E_{\lambda_{k-j}} $$
 qui est \'egal \`a \ $\C.b^{\lambda_k-\lambda_{k-j}}.e_{\lambda_{k-j}}$ \ puisque ce noyau est isomorphe \`a \ $Hom_{\A}(F_j, E_{\lambda_{k-j}}) $ \ qui est de dimension  au plus \ $1$ \  puisque \ $F_j$ \ est un th\`eme, ce qui prouve notre seconde assertion. $\hfill \blacksquare$

    \section{Endomorphismes et th\`emes stables.}
   
   \subsection{Injections entre deux th\`emes primitifs de m\^eme rang.}
  
  Commen{\c c}ons par \'etudier les injections \ $\A-$lin\'eaires entre deux th\`emes \ $[\lambda]-$primitifs.

\begin{lemma}\label{necessaire}
Soitent \ $E' \subset E$ \ deux th\`emes \ $[\lambda]-$primitifs de m\^eme rang \ $k$. Soient \ $\mu_1, \dots, \mu_k$ \ et \ $\lambda_1, \dots, \lambda_k$ \ leurs invariants fondamentaux respectifs. Alors on a
\begin{enumerate}[i)]
\item $\forall j \in [1,k] \quad \mu_j \geq \lambda_j $ ;
\item \ $\dim_{\C}(E\big/E') = \sum_{j=1}^k \ \mu_j - \lambda_j .$
\end{enumerate}
\end{lemma}

\parag{Preuve} Montrons cela par r\'ecurrence sur le rang \ $k$. Le r\'esultat \'etant clair pour \ $k=1$, supposons-le montr\'e en rang \ $k-1$. Comme \ $F'_{k-1} \subset F_{k-1}$, on obtient imm\'ediatement les in\'egalit\'es \ $\mu_j \geq \lambda_j \quad \forall j \in [1,k-1]$. De plus, la restriction \`a \ $E'$ \ de l'application \ $\pi_k : E \to E\big/F_{k-1} \simeq E_{\lambda_k}$ \ n'est pas nulle, sinon \ $E'$ \ serait contenu dans \ $F_{k-1}$ \ et ne serait pas de rang \ $k$. On a donc un quotient de rang \ $1$ \ de \ $E'$ \ contenu dans \ $E_{\lambda_k}$. Ceci donne \ $\mu_k \geq \lambda_k$.\\
Enfin on a la suite exacte
$$ 0 \to F_{k-1}\big/F'_{k-1} \to E\big/E' \to E\big/ (F_{k-1}+E') \to 0 $$
puisque \ $F_{k-1} \cap E'$ \ est normal et de rang \ $k-1$ \ dans \ $E'$, comme noyau de \ $(\pi_k)_{\vert E'}$, et donc \'egal \`a \ $F'_{k-1}$. L'hypoth\`ese de r\'ecurrence donne \\ 
$\dim_{\C} F_{k-1}\big/F'_{k-1} = \sum_{j=1}^{k-1} \ \mu_j - \lambda_j $. De plus, puisque \ $\pi_k(E') = E_{\mu_k}$, le quotient \ $E\big/ (F_{k-1}+E') $ \  est de dimension \ $\mu_k - \lambda_k$ \ et la suite exacte permet de conclure. $\hfill \blacksquare$

  \bigskip
   
   \begin{thm}\label{inclusion}
   Soient \ $E'$ \ et \ $E$ \ deux th\`emes \ $[\lambda]-$primitifs de rang \ $k$.  L'espace vectoriel des morphismes \ $\A-$lin\'eaires de \ $E'$ \ dans \ $E$ \ modulo ceux qui sont de rang \ $\leq k-1$ \ est de dimension \ $\leq 1$.\\
   Supposons que les invariants fondamentaux  respectifs \ $\mu_1, \dots, \mu_k$ \ et \ $\lambda_1, \dots, \lambda_k$ \ de \ $E'$ \ et \ $E$ \ v\'erifient la condition
   $$ \mu_j - \lambda_j \geq k-1 \quad \forall j \in [1,k] .$$
   Alors il existe une injection \ $\A-$lin\'eaire \ $ i : E' \hookrightarrow E$.
      \end{thm}
   
   \parag{D\'emonstration} Montrons par r\'ecurrence sur le rang la premi\`ere assertion : comme elle est clair en rang 1, supposons-l\`a d\'emontr\'ee en rang \ $k-1$, et consid\'erons deux injections \ $\varphi_1$ \ et \ $\varphi_2$ \ de \ $E'$ \ dans \ $E$. Leurs restrictions \`a \ $F'_{k-1}$ \ sont des injections dans \ $F_{k-1}$, et l'hypoth\`ese de r\'ecurrence fournit un \ $\alpha \in \C^*$ \ tel que \ $\varphi_1 - \alpha.\varphi_2$ \ ne soit plus injective dans \ $F'_{k-1}$, donc \`a fortiori dans \ $E'$.
   
   \smallskip
   
  Montrons le second r\'esultat par r\'ecurrence sur \ $k \geq 1$. Comme le cas \ $k = 1$ \ est imm\'ediat, supposons \ $k \geq 2 $ \ et le r\'esultat prouv\'e en rang \ $k-1$ \ sous la forme suivante :
   Soient \ $F'_{k-1}$ \ et \ $F_{k-1}$ \ deux th\`emes de rang \ $k-1$ \ v\'erifiant \ $\mu_j - \lambda_j \geq k-2$ \ pour \ $j \in [1,k-1]$. Notons 
   \begin{align*}
   & Q' : = (a -\mu_1.b).T_1^{-1}\dots T_{k-2}^{-1}.(a - \mu_{k-1}) \\
   & Q : = (a - \lambda_1.b).S_1^{-1} \dots S_{k-2}^{-1}.(a -\lambda_{k-1})
   \end{align*}
   les g\'en\'erateurs respectifs des annulateurs dans \ $\A$ \ des g\'en\'erateurs respectifs \ $\varepsilon_{k-1}$ \ et \ $e_{k-1}$ \ de \ $F'_{k-1}$ \ et \ $F_{k-1}$. Alors il existe un \'el\'ement 
   $$ x : = \sigma.b^{\mu_{k-1} -\lambda_{k-1}}.e_{k-1} + \sum_{h=1}^{k-2} \ V_h.e_h $$
   de \ $F_{k-1}$ \ qui est annul\'e par \ $Q'$, avec \ $V_h \in \C[[b]], \forall h \in [1,k-2]$ \ et \ $\sigma \not= 0$. \\
   Ceci implique l'existence d'une injection \ $\A-$lin\'eaire de \ $F'_{k-1}$ \ dans \ $F_{k-1}$ \ donn\'ee en envoyant le g\'en\'erateur \ $\varepsilon_{k-1}$ \ de \ $F'_{k-1}$ \ sur \ $x$.
   
   \bigskip
   
   Appliquons cette hypoth\`ese de r\'ecurrence de la fa{\c c}on suivante : soient \ $F'_{k-1}$ \ et \ $F_{k-1}$ \ les  sous-th\`emes normaux de rang \ $k-1$ \ de \ $E'$ \ et \ $E$ \ respectivement. Comme on a \ $\mu_j - \lambda_j \geq k-1 $ \ pour \ $j \in [1,k-1]$, on peut appliquer l'hypoth\`ese de r\'ecurrence \`a \ $F'_{k-1}$ \ et \`a \ $b.F_{k-1}$. Ce qui signifie que l'\'el\'ement \ $x$ \ fournit par l'hypoth\`ese de r\'ecurrence est dans \ $b.F_{k-1}$, et l'on aura, puisque \ $\mu_{k-1} - \lambda_{k-1} \geq k- 1\geq 1 $ \ simplement  \ $V_h \in b.\C[[b]]$ \ pour chaque \ $h \in [1,k-1]$.\\

    Notons par \ $\varepsilon_k$ \ et \ $e_k$ \ les g\'en\'erateurs de \ $E'$ \ et \ $E$, et par \ $P' : = Q'.T_{k-1}^{-1}.(a - \mu_k.b)$ \ et \ $P : = Q.S_{k-1}^{-1}.(a - \lambda_k.b)$ \ les g\'en\'erateurs des annulateurs respectifs de \ $\varepsilon_k$ \ et \ $e_k$. Cherchons alors un \'el\'ement \ $y \in E$ \ de la forme
   $$ y = \tau.b^{\mu_k-\lambda_k}.e_k + \sum_{h=1}^{k-1} \ W_h.e_h $$
   v\'erifiant les conditions suivantes :
   \begin{enumerate}[i)]
   \item  \ $\tau \not= 0 $.
   \item \ $W_h \in \C[[b]] \quad  \forall h \in [1,k-1]$.
   \item \ $ (a - \mu_k.b).y = T_{k-1}.x $.
   \end{enumerate}
   Il donnera alors une injection de \ $E'$ \ dans \ $E$, puisque \ $y \not\in F_{k-1}$ \ et que \ $P.y = 0$.\\
   Remarquons que comme la suite \ $\lambda_j + j$ \ est  croissante, on a
   $$ \mu_k \geq \lambda_k + k-1 \geq \lambda_j + j -1  \geq \lambda_j  \quad \forall j \in [1,k] .$$
   La relation iii) donne les \'equations suivantes :
   \begin{align*}
   & b.W'_{k-1} - (\mu_k - \lambda_{k-1}).W_{k-1} = \sigma.T_{k-1}.b^{\mu_{k-1}-\lambda_{k-1}-1} - \tau.S_{k-1}.b^{\mu_k-\lambda_k-1} \\
   & b^2.W'_h - (\mu_k-\lambda_h).b.W_h = T_{k-1}.V_h - S_h.W_{h+1} \tag{h} 
   \end{align*} 
   La premi\`ere \'equation aura une solution dans \ $\C[[b]]$, unique \`a \ $\C.b^{\mu_k-\lambda_{k-1}}$ \ pr\`es, pourvu que le coefficient de \ $b^{\mu_k-\lambda_{k-1}}$ \ soit nul dans le membre de droite. Si \ $\alpha' \not= 0$ \ est le coefficient de \ $b^{p'_{k-1}}$ \ dans \ $T_{k-1}$ \ et \ $\alpha \not= 0$ \ celui de \ $b^{p_{k-1}}$ \ dans  \ $S_{k-1}$, il nous suffit de choisir \ $\sigma = \tau.\alpha/\alpha'$ \ pour assurer l'existence de \ $W_{k-1} \in b^{k-2}.\C[[b]]$, en fait unique \`a \ $\C.b^{\mu_k-\lambda_{k-1}}$ \ pr\`es.
   
   \bigskip
   
   Supposons prouv\'e l'existence de \ $W_{h+1} \in b^h.\C[[b]]$, unique modulo \ $\C.b^{\mu_k-\lambda_{h+1}}$, et \ $h \geq 1$. Comme \ $V_h $ \ et \ $W_{h+1}$ \ sont dans \ $ b.\C[[b]]$, pour que l'\'equation \ $(h)$ \ ait une solution, unique modulo \ $\C.b^{\mu_k-\lambda_h}$,  il suffit de s'assurer que le coefficient de \ $b^{\mu_k-\lambda_h +1}$ \ dans  \ $T_{k-1}.V_h - S_h.W_{h+1}$ \ est nul. Mais comme le coefficient de \ $b^{p_h}$ \ de \ $S_h$ \ est non nul et que l'on peut fixer arbitrairement le coefficient de \ $b^{\mu_k-\lambda_{h+1}}$ \ dans \ $W_{h+1}$, ceci ne pose pas de probl\`eme \`a l'aide d'un choix convenable de \ $W_{h+1}$ puisque l'on a \ $\lambda_{h+1} = \lambda_h + p_h -1 $ \ qui donne   \ $ \mu_k-\lambda_h +1 =  \mu_k-\lambda_{h+1} + p_h $. $\hfill \blacksquare$
   
   \parag{Remarques}
   \begin{enumerate}[i)]
   \item Une cons\'equence de la d\'emonstration du th\'eor\`eme est que s'il existe une injection de \ $E'$ \ dans \ $E$, elle envoie le g\'en\'erateur \ $e'$ \ de \ $E'$ \ sur
   $$ \tau.b^{\mu_k - \lambda_k}.e \quad modulo \  F_{k-1} $$
   o\`u \ $e'$ \ et \ $e$ \ sont les g\'en\'erateurs "standards" de \ $E'$ \ et \ $E$, et \ $\tau \in \C^*$ \ est arbitraire.
   \item L'exemple 5.1.1  de l'appendice  5.1 nous fournit deux th\`emes primitifs de rang 3 \ $E' : = E\big/F_1$ \ et \ $E : = F_3$,  v\'erifiant \ $\mu_j - \lambda_j \geq k-2 = 1 $ \ et tels que \ $E'$ \ ne s'injecte pas dans \ $E$. En effet on a dans cet exemple \ $\mu_1 = \lambda_1+1, \mu_2 = \lambda_1 + 3, \mu_3 = \lambda_1 + 4$ \ et \ $\lambda_1 = \lambda_1, \lambda_2 = \lambda_1+ 1, \lambda_3 = \lambda_1 +3$. $\hfill \square$
   \end{enumerate}
   
   \bigskip
   
   \begin{cor}\label{drapeau}   
   Soit \ $E$ \ un th\`eme \ $[\lambda]-$primitif, et soit \ $R_j \subset Hom_{\A}(E,E)$ \ le sous-espace vectoriel des endomorphismes de rang \ $\leq j$ \ de \ $E$. Alors pour chaque \ $j \in [0,k-1]$ \ l'espace vectoriel  complexe \ $R_{j+1}\big/R_j$ \ est de dimension \ $\leq 1$.\\
En particulier, on a toujours \ $ \dim_{\C}(Hom_{\A}(E,E)) \leq k$ \ avec \'egalit\'e si et seulement si \ $(R_j)_{j \in[1,k]}$ \ est un drapeau complet de \ $Hom_{\A}(E,E)$, c'est-\`a-dire que chaque quotient \ $R_{j+1}\big/R_j$ \ est de dimension 1 pour \ $j \in [0,k-1]$.
\end{cor}

\parag{Preuve} Soit \ $\varphi : E \to E$ \ un morphisme de rang \ $j$. Alors son noyau est \ $F_{k-j}$ \ puisque ce noyau est normal et de rang \ $k-j$. De plus, le normalis\'e de son image est \ $F_j$. Donc \ $\varphi$ \ se factorise de la fa{\c c}on suivante :
$$ E \to E\big/F_{k-j} \overset{f}{\to} F_j \hookrightarrow E $$
o\`u la premi\`ere fl\`eche est le quotient et la derni\`ere l'injection naturelle. La fl\`eche \ $f$ \ est injective, et la correspondance \ $\varphi \to f$ \ induit une bijection \ $\C-$lin\'eaire entre \ $R_j\big/R_{j-1}$ \ et l'espace vectoriel des morphismes de \ $E\big/F_{k-j}$ \ dans \ $F_j$, modulo ceux qui sont de rang \ $\leq j-1$ \ (ou non injectifs, ce qui revient au m\^eme). La premi\`ere assertion du th\'eor\`eme permet alors de conclure. $\hfill \blacksquare $
   
  \bigskip
  
   \begin{cor}\label{Endo}
   Soit \ $E$ \ un th\`eme \ $[\lambda]-$primitif de rang \ $k$. Une condition \\
   suffisante pour qu'il existe une injection \ $\A-$lin\'eaire de \ $E\big/F_j$ \ dans \ $F_{k-j}$ \ est que pour chaque \ $h \in [1,k-j]$ \ on ait
   $$ p_h + \cdots + p_{h+j-1} \geq k-1 .$$
   En particulier, pour \ $p_h \geq k-1 \quad \forall h \in [1,k-1]$ \ l'espace vectoriel \ $Hom_{\A}(E,E)$ \ sera de dimension \ $k$.
   \end{cor}
   
   \parag{Preuve} Comme \ $E\big/F_j$ \ est un th\`eme  \ $[\lambda]-$primitif de rang \ $k-j$ \ et d'invariants fondamentaux \ $\lambda_{j+1}, \dots, \lambda_k$ \ et que \ $F_{k-j}$ \ est \'egalement th\`eme  \ $[\lambda]-$primitif de rang \ $k-j$ \ et d'invariants fondamentaux \ $\lambda_1, \dots, \lambda_{k-j}$, nous pouvons conclure gr\^ace au th\'eor\`eme d\`es que l'on a
   $$ \lambda_{j+h} - \lambda_h \geq k-j-1 \quad \forall h \in [1,k-j].$$
   Mais \ $\lambda_{j+h} - \lambda_h  = p_h+ \dots +p_{j+h-1} -j $ \ et donc le corollaire se d\'eduit du th\'eor\`eme \ref{inclusion} et de son premier corollaire \ref{drapeau}. $\hfill \blacksquare$
   
   \bigskip
   
   On voit facilement que la condition {\em  n\'ecessaire} pour avoir une injection de \ $E\big/F_j $ \ dans \ $F_{k-j}$ \ donn\'ee par le lemme \ref{necessaire} correspond aux  in\'egalit\'es
   $$ p_h + \dots + p_{h+j-1} \geq \ j \quad {\rm pour} \quad h \in [1,k-j].$$
   Elle est trivialement v\'erifi\'ee si on a \ $p_j \geq 1, \forall j \in [1,k]$, c'est-\`a-dire si la suite \ $\lambda_1, \dots \lambda_k$ \ est croissante (large).

  \subsection{Th\`emes primitifs stables.}

Commen{\c c}ons par rappeler deux remarques simples.

\parag{Remarques}
\begin{enumerate}[1)]
\item Si \ $E$ \ est un th\`eme \ $[\lambda]-$primitif, l'espace vectoriel \ $Hom_{\A}(E, E_{\lambda})$ \ est de dimension au plus \'egale \`a \ $1$. En effet, comme \ $E$ \ admet un unique quotient de rang \ $1$, \`a savoir \ $E\big/F_{k-1} \simeq E_{\lambda_k}$, un morphisme non nul est n\'ecessairement une injection de \ $E_{\lambda_k}$ \ dans \ $E_{\lambda}$. L'espace vectoriel \ $Hom_{\A}(E, E_{\lambda})$ \  sera donc nul pour \ $\lambda_k < \lambda $ \ et de dimension \ $1$ \ pour \ $\lambda_k = \lambda + q$ \ avec \ $q \in \mathbb{N}$. 
\item Soient \ $E_1$ \ et \ $E_2$ \ deux th\`emes \ $[\lambda]-$primitifs de rangs \ $k_1$ \ et \ $k_2$, et soit \ $i : E_1 \hookrightarrow E_2$ \ une injection \ $\A-$lin\'eaire. Alors on a \ $k_2 \geq k_1$, et le normalis\'e de \ $i(E_1)$ \ est le sous-th\`eme normal de rang \ $k_1$ \ de \ $E_2$. En particulier on aura \ $i(E_1)$ \ qui sera contenu dans le sous-th\`eme normal de rang \ $k_1$ de \ $E_2$. Cette image est m\^eme de codimension finie dans ce sous-th\`eme.$\hfill \square$
\end{enumerate}

\begin{prop}\label{stable 0}
Soit \ $E$ \ un th\`eme \ $[\lambda]-$primitif de rang \ $k \geq 1$ \ et soit \ $\varphi_0 $ \ un endomorphisme de \ $E$ \ de rang \ $k-1$. Les propri\'et\'es suivantes sont v\'erifi\'ees:
\begin{enumerate}[i)]
\item Pour chaque \ $j \in [0,k]$ \ le rang de \ $\varphi_0^j$ \ est \ $k-j$.
\item Une base de l'espace vectoriel \ $End_{\A}(E)$ \ est donn\'ee par \ $\id, \varphi_0, \dots, \varphi_0^j, \dots, \varphi_0^{k-1}$. En particulier on a \ $\dim_{\C}(End_{\A}(E)) = k$, et cette alg\`ebre est commutative et  isomorphe \`a \ $\C[x]\big/(x^k)$.
\item Pour chaque \ $j \in [1,k-1]$ \ la restriction \`a \ $F_j$ \ de \ $ \varphi_0$ \ est de rang \ $j-1$ ; donc la restriction \ $End_{\A}(E) \to End_{\A}(F_j)$ \ est surjective.
\item Pour chaque \ $j \in [1,k-1]$ \ on a un endomorphisme de \ $E\big/F_j$ \  induit par \ $\varphi_0$, et il est de rang \ $k-j-1$. On a donc \'egalement une application lin\'eaire surjective de \ $End_{\A}(E) $ \ dans \ $End_{\A}(E\big/F_j)$.
\end{enumerate}
\end{prop}

\parag{Preuve}  
Montrons l'assertion i) par r\'ecurrence sur \ $j$. Comme elle est claire pour \ $j = 0, j = 1$, supposons montr\'e que \ $\varphi_0^j$ \ est de rang \ $k-j$ \ pour \ $j\in [1,k-1]$, et montrons que \ $\varphi_0^{j+1}$ \ est de rang \ $k-j-1$.\\
 Le noyau de \ $\varphi_0^j$ \ est un sous-th\`eme normal de rang \ $j$. Il est donc \'egal \`a \ $F_j$, l'unique sous-th\`eme normal de rang\ $ j$ \  de \ $E$. Donc \ $\varphi_0^j$ \ induit une injection de \ $E\big/F_j$ \ dans \ $F_{k-j}$ \ dont l'image \ $\Phi_j$ \ est de codimension finie dans \ $F_{k-j}$ \ (voir la remarque 2) ci-dessus). En composant \`a nouveau avec \ $\varphi_0$ \ dont le noyau \ $F_1$ \ rencontre \ $\Phi_j$ \ en un sous-(a,b)-module de rang \ $1$, puisque \ $F_1 \subset F_{k-j}$,  on en d\'eduit que le noyau de \ $\varphi_0$ \ restreinte \`a \ $\Phi_j$ \ est de rang \ $1$, et donc que son image, qui est \ $\Phi_{j+1}$, c'est-\`a-dire l'image de \ $\varphi_0^{j+1}$, est de rang \ $k-j-1 = rg(\Phi_j) - 1$. Donc i) est d\'emontr\'ee.\\
 Si on a une relation lin\'eaire dans \ $End_{\A}(E)$
 $$ \sum_{j=0}^{k-1} \ \alpha_j.\varphi_0^j = 0 $$
   les nombres complexes \ $\alpha_j$ \ n'etant pas tous nuls, soit \ $j_0$ \ le premier entier pour lequel on a \ $\alpha_{j_0} \not= 0 $. Alors on aura \ $\Phi_{j_0} \subset \sum_{h=j_0+1}^{k-1} \ \Phi_h $. Mais pour chaque \ $j$ \ on a \ $\Phi_{j} \subset F_{k-j}$, et donc \ $\Phi_{j_0} \subset F_{k-j_0-1}$, ce qui contredit le fait que \ $\varphi_0^{j_0}$ \ soit de rang \ $k-j_0$. Donc on a \ $k$ \ vecteurs ind\'ependants dans \ $End_{\A}(E)$, et comme on sait (voir corollaire 4.0.13) que cet espace vectoriel est de dimension au plus \'egale \`a \ $k$, il est de dimension \ $k$ \ et on a une base de \ $End_{\A}(E)$.\\
Montrons iii). Comme la restriction de \ $\varphi_0$ \ \`a \ $F_j$ \ est de noyau \ $F_1 \subset F_j$, le rang est bien \ $j-1$. La surjectivit\'e annonc\'ee r\'esulte alors de ii) appliqu\'e au th\`eme \ $F_j$.\\
Comme la restriction de \ $\varphi_0$ \ \`a \ $F_j$ \ est de rang \ $j-1$, on a \ $\varphi_0(F_j) \subset F_{j-1} \subset F_j$, et \ $\varphi_0$ \ induit un endomorphisme de \ $E\big/F_j$. Comme l'image \ $\Phi_1$ \ de \ $\varphi_0$ \ est de codimension finie dans \ $F_{k-1}$, le quotient \ $\Phi_1 \big/F_j \cap \Phi_1$ \ est de rang \ $k-j-1$, puisque \ $F_j \subset F_{k-1}$ \ montre que l'on quotiente un \ $\C[[b]]-$module de rang \ $k-1$ \ par un sous$-\C[[b]]-$module de rang \ $j$. La surjectivit\'e de l'application lin\'eaire  \ $End_{\A}(E) \to End_{\A}(E\big/F_j)$ \ se d\'eduit alors de ii) appliqu\'e au th\`eme \ $E\big/F_j$.  $\hfill \blacksquare$
   
\parag{\bf Exemple important} Soit \ $e \in \Xi_{\lambda}^{(k-1)}$ \ et supposons que le th\`eme \ $\A.e $ \ soit de rang \ $k$ \ et {\bf stable par la monodromie} \ $\mathcal{T}$ \ de \ $\Xi_{\lambda}^{(k-1)}$. Rappelons que  pour \ $j \in [0,k-1]$
$$\mathcal{T}\big(s^{\lambda-1}.\frac{(Log\,s)^j}{j!}\big)= \exp(2i\pi.\lambda).s^{\lambda-1}\frac{(Log\,s + 2i\pi)^j}{j!}$$ 
et que \ $\mathcal{T}$ \ commute \`a l'action de \ $\A$ \ sur \ $\Xi_{\lambda}^{(k-1)}$. Alors \ $\mathcal{T} - \exp(2i\pi.\lambda).\id$ \ induit un endomorphisme de rang \ $k-1$ \ sur \ $\A.e$.\\
En effet, on peut supposer que 
 $$e = s^{\lambda_k-1}.\frac{(Log\,s)^{k-1}}{(k-1)!}  \quad modulo \quad \Xi_{\lambda}^{(k-1)}$$
  et donc que  \ 
 $\varepsilon : = (\mathcal{T} - \exp(2i\pi.\lambda).\id)(e)$ \ sera dans \ $\Xi_{\lambda}^{(k-2)}\setminus \Xi_{\lambda}^{(k-3)} $, c'est \`a dire dans \ $F_{k-1}\setminus F_{k-2}$, puisque l'on a suppos\'e la stabilit\'e de \ $E : = \A.e$ \ par \ $\mathcal{T}$ \ et que le terme en 
  $$s^{\lambda_k-1}.\frac{(Log\,s)^{k-2}}{(k-2)!}$$
  ne peut dispara\^itre, \'etant donn\'e que \ $\mathcal{T} - exp(2i\pi.\lambda).\id$ \ fait strictement d\'ecro\^itre le degr\'e en \ $Log\, s$ \ dans \ $\Xi_{\lambda}$.\\
 Mais si \ $P \in \A$ \ engendre l'annulateur de \ $e$ \ dans \ $\Xi_{\lambda}^{(k-1)}$, on aura \ $P.\varepsilon = 0$, puisque \ $\mathcal{T}$ \ commute \`a l'action de \ $\A$ ; ceci montre que l'on a bien un endomorphisme de \ $E$ \  en posant \ $\varphi_0(e) = \varepsilon$, et que cet endomorphisme est de rang \ $k-1$. En effet, pour \ $\varepsilon \in \Xi_{\lambda}$, le th\`eme \ $[\lambda]-$primitif \ $\A.\varepsilon$ \  est de rang \ $l$ \ si et seulement si  \ $l-1$ \ est le degr\'e en \ $Log\, s$ \ de \ $\varepsilon$. $\hfill \square$

\bigskip

\begin{defn}\label{stable 1} On dira qu'un th\`eme \ $[\lambda]-$primitif \ $E$ \ de rang \ $k$ \ est {\bf stable} s'il admet un endomorphisme de rang \'egal \`a \ $k-1$.
\end{defn}

La proposition \ref{stable 0} implique imm\'ediatement le corollaire suivant :

\begin{cor}[de la proposition \ref{stable 0}]
Soit \ $E$ \ un th\`eme \ $[\lambda]-$primitif \ $E$ \ de rang \ $k$. Si \ $E$ \ est stable, tout sous-th\`eme normal et tout th\`eme quotient de \ $E$ \ est stable.
\end{cor}

\begin{lemma}\label{stable 2}
 Les propri\'et\'es suivantes sont \'equivalentes pour un th\`eme \ $[\lambda]-$primitif \ $E$ \ de rang \ $k$ :
\begin{enumerate}[i)]
\item \ $E$ \ est stable.
\item La dimension de \ $End_{\A}(E)$ \ est \'egale \`a \ $k$.
\item L'image d'une injection \ $\A-$lin\'eaire de \ $E$ \ dans \ $\Xi_{\lambda}^{(k-1)}$ \ est ind\'ependante de l'injection choisie.
\item Il existe une injection \ $\A-$lin\'eaire de \ $E$ \ dans \ $\Xi_{\lambda}^{(k-1)}$ \ dont l'image est stable par la monodromie  \ $\mathcal{T}$.  $\hfill \square$
\end{enumerate}
\end{lemma}

\parag{Preuve} L'implication \ $i) \Rightarrow ii)$ \ est montr\'ee dans la proposition \ref{stable 0}. L'implication \ $ii) \Rightarrow iii)$ \ r\'esulte du fait que si \ $i$ \ est une injection de \ $E$ \ dans \ $\Xi_{\lambda}^{(k-1)}$, la composition par \ $i$ \ donne une injection lin\'eaire de \ $End_{\A}(E)$ \ dans \ $Hom_{\A}(E, \Xi_{\lambda}^{(k-1)})$. Ces deux espaces vectoriels \'etant de m\^eme dimension \ $k$, le premier par hypoth\`ese, le second d'apr\`es le th\'eor\`eme 2.2.1 de [B.05], en remarquant qu'une application \ $\A-$lin\'eaire d'un th\`eme \ $[\lambda]-$primitif de rang \ $k$ \ dans \ $\Xi$ \ a toujours son image dans \ $\Xi_{\lambda}^{(k-1)}$, on en d\'eduit que toute injection de \ $E$ \ dans \ $\Xi_{\lambda}^{(k-1)}$ \ est de la forme \ $\varphi\circ i$ \ o\`u \ $\varphi \in Aut_{\A}(E)$, ce qui prouve iii).\\
L'implication \ $iii) \Rightarrow iv)$ \ est facile puisque \ $\mathcal{T}\circ i$ \ est encore une injection \ $\A-$lin\'eaire de \ $E$ \ dans \ $\Xi_{\lambda}^{(k-1)}$ \ quand \ $i$ \ l'est.\\
L'implication  \ $iv) \Rightarrow i)$ \ r\'esulte de l'exemple important trait\'e plus haut. $\hfill \blacksquare$

\begin{cor}\label{dual stable}
Soit \ $E$ \ un th\`eme \ $[\lambda]-$primitif stable d'invariants fondamentaux \ $\lambda_1, \dots, \lambda_k$. Alors pour \ $\delta \in \mathbb{Q}$ \ v\'erifiant \ $\delta - \lambda_k > k-1$ \ le th\`eme \ $E^* \otimes E_{\delta}$ \ est stable.
\end{cor}

\parag{Preuve} Il nous suffit de montrer que \ $End_{\A}(E^* \otimes E_{\delta})$ \ est de dimension \\
 $k : = rg(E)$. Mais cet espace vectoriel est isomorphe \`a \ $End_{\A}(E^*)$ \ puisque pour un (a,b)-module \ $F$, le produit tensoriel \ $F \otimes_{a,b} E_{\delta}$ \ consiste \`a regarder \ $F$ \ en changeant \ $a$ \ en \ $a + \delta.b$, ce qui ne change pas les endomorphismes \ $\A-$lin\'eaires.\\
Par ailleurs la transposition donne une application \ $\C-$lin\'eaire \ $End_{\A}(E) \to End_{\A}(E^*)$ \ qui est clairement bijective, puisque \ $(E^*)^*$ \ est canoniquement isomorphe \`a \ $E$.\\
 Le r\'esultat en d\'ecoule. $\hfill \blacksquare$

\bigskip

\begin{lemma}\label{stable 3}
Soit \ $E$ \ un th\`eme \ $[\lambda]-$primitif  de rang \ $k$, et supposons que l'on ait \ $p_{k-1} = 0$ \ avec \ $k \geq 2$,  ou bien  \ $p_{k-1} = 1$ \ et \ $p_{k-2} \geq 2$ \ avec \ $k \geq 3$. Alors \ $E$ \ n'est pas stable.
\end{lemma}

\parag{Preuve} Soit \ $e$ \ un g\'en\'erateur standard de \ $E$, et soit \ $\A.P$ \ son annulateur. Il nous suffit de montrer qu'il n'existe pas d'\'el\'ement \ $x \in F_{k-1} \setminus F_{k-2}$ \ qui soit annul\'e par \ $P$. Un tel \'el\'ement doit v\'erifier
\begin{equation*}
 (a - \lambda_k.b).x = S_{k-1}.y \quad {\rm avec} \quad y \in F_{k-2}  \tag{*}
 \end{equation*}
et \ $Q.y = 0 $ \ o\`u l'on a pos\'e \ $P : = Q.S_{k-1}^{-1}.(a - \lambda_k.b) $. On sait (voir le lemme \ref{lemme manquant})  que l'on   peut \'ecrire
\begin{align*}
&  x = b^{\lambda_k - \lambda_{k-1}}.e_{k-1} + \sum_{j=1}^{k-2} \quad U_j.e_j  \quad {\rm avec} \quad U_j \in \C[[b]] \quad \forall j \in [1,k-2] \\
&  y - \rho.b^{\lambda_{k-1}-\lambda_{k-2}}.e_{k-2} \in  F_{k-3} \quad {\rm avec} \quad \rho \in \C^*
\end{align*}
Si on a \ $\lambda_k - \lambda_{k-1} = -1$, un tel \ $x$ \ ne peut exister. Supposons donc \ $k \geq 3$ \ et  \ $\lambda_k = \lambda_{k-1}$, c'est \`a dire \ $p_{k-1} = 1$. 
En rempla{\c c}ant dans l'\'equation \ $(^*)$, on obtient que \ $U_{k-2}$ \ doit v\'erifier l'\'equation suivante :
\begin{equation*}
 S_{k-2} + b^2.U'_{k-2} - (\lambda_k - \lambda_{k-2}).b.U_{k-2} = \rho.S_{k-1}.b^{\lambda_{k-1}-\lambda_{k-2}} . \tag{**}
 \end{equation*}
Comme on a suppos\'e \ $\lambda_{k-1} > \lambda_{k-2}$, c'est-\`a-dire \ $p_{k-2} \geq 2$, l'\'equation \ $(^{**})$ \ ne peut avoir de solution, puisque \ $S_{k-2}(0) = 1$. $\hfill \blacksquare$

\bigskip

Le lemme pr\'ec\'edent admet la cons\'equence imm\'ediate suivante :

\begin{cor}\label{stable 4}
Soit \ $E$ \ un th\`eme \ $[\lambda]-$primitif  de rang \ $k$ \ stable. Alors ou bien la suite \ $\lambda_1, \dots, \lambda_k$ \ est strictement croissante, ou bien elle est constante.\\
Dans le cas  o\`u l'on a \ $\lambda_1 = \dots = \lambda_k$, nous dirons que le th\`eme stable \ $E$ \ est {\bf sp\'ecial}.
\end{cor}

\parag{Preuve} Comme le cas o\`u le rang est \ $\leq 2$ \ est clair d'apr\`es le lemme \ref{rang 2}, nous pouvons supposer \ $k \geq 3$. 
Commen{\c c}ons par montrer qu'il existe \ $j_0 \in [1,k]$ \ tel que l'on ait 
$$ \lambda_1 = \dots = \lambda_{j_0} < \lambda_{j_0+1} < \dots < \lambda_k .$$
Soit \ $E$ \ un th\`eme stable de rang \ $k$ \ et d'invariants fondamentaux \ $\lambda_1, p_1, \dots, p_{k-1}$. Chaque \ $p_j$ \ est au moins \'egal \`a 1. En effet si on a \ $p_j = 0$ \ cela revient \`a dire que \ $F_{j+1}\big/F_{j-1}$ \  est isomorphe \`a \ $E_{\lambda_j,\lambda_j}$ \ qui n'est pas stable. Ceci contredit le corollaire \ref{stable 0}. \\
Si tous les \ $p_j$ \ sont au moins \'egaux \`a  2, la suite des \ $\lambda_j$ \  est strictement croissante et on pose \ $j_0 = 1$. Sinon soit  \ $j_0$ \ le plus grand  entier dans \ $[1, k-1]$ \ tel que \ $p_{j} = 1$ \ pour \ $j \leq j_0-1$. On a donc \ $\lambda_1 = \dots = \lambda_{j_0} < \lambda_{j_0+1}$. Donc le th\`eme stable \ $E\big/F_{j_0-1}$ \ admet comme invariants fondamentaux \ $\mu_1 : = \lambda_{j_0} < \mu_2 = \lambda_{j_0+1} \leq \dots \leq \mu_{k- j_0-1} = \lambda_k $. Nous voulons montrer qu'alors la suite \ $\mu_1, \dots, \mu_{k- j_0-1} $ \ est strictement croissante. Supposons qu'elle croisse strictement jusqu'\`a \ $\mu_{h-1} < \mu_h = \mu_{h+1}$, o\`u l'on pose \ $\mu_0 = \lambda_{j_0-1}$ \ dans le cas \ $h = 1$. Alors la seconde assertion du lemme \ref{stable 3} appliqu\'ee au th\`eme stable \ $F_{j_0+h+1}\big/F_{j_0+h-2}$ \ de rang 3 donne la contradiction cherch\'ee.\\
On conclut alors  gr\^ace \`a la remarque 1) qui suit la d\'emonstration du  th\'eor\`eme de dualit\'e \ref{dual d'un theme}. En effet pour \ $\delta \in \mathbb{Q}$ \ assez grand \ $E^*\otimes_{a,b} E_{\delta}$ \ est un th\`eme, et il est stable d'apr\`es le corollaire \ref{dual stable}. On a donc ou bien \ $j_0 = 1$, ou bien \ $j_0 = k$. $\hfill \blacksquare$

\parag{Remarques}\begin{enumerate}[i)]
\item Le cas sp\'ecial impose les  \'egalit\'es \ $p_j = 1 \quad \forall j \in [1,k-1]$.
\item  Le lemme \ref{non stables} de l'appendice 5.1 montre qu'il existe des th\`emes stables sp\'eciaux de rang \ $3$.$\hfill \square$ 
\end{enumerate}

\subsection{Forme canonique pour un th\`eme primitif.}

\subsubsection{Suppl\'ementaires.}
     
     \begin{prop}\label{Supplementaires}
     Soit \ $E : = \A\big/\A.P$ \ un th\`eme \ $[\lambda]-$primitif de rang \ $k$ \ o\`u l'on a pos\'e :
     \begin{enumerate}[1)]
     \item \ $ P : = (a - \lambda_1.b).S_1^{-1} \dots S_{k-1}^{-1}.(a - \lambda_k.b) $.
     \item \ $\lambda_1, p_1, \cdots, p_{k-1}$ \ sont les invariants fondamentaux de \ $E$.
     \item \ $S_1, \cdots, S_{k-1}$ \ sont des \'el\'ements inversibles de \ $\mathbb{C}[[b]]$ \ de termes constants \'egaux \`a \ $1$ tels que le coefficient de \ $b^{p_j}$ \ dans \ $S_j$ \ soit non nul.
     
     \end{enumerate}
     Pour \ $j \in [1,k-1]$ \ d\'efinissons \ $P_j : = (a - \lambda_{j+1}.b).S_{j+1}^{-1} \cdots S_{k-1}^{-1}.(a - \lambda_k.b) $.\\
     Si \ $ p_j + \cdots + p_{k-1} \geq k-j $ \ notons \ $q_j : = p_j + \cdots + p_{j+h}$, o\`u \ $h$ \ est le plus petit entier tel que \ $ p_j + \cdots + p_{j+h} \geq k-j$, posons  \ $V_j : = \oplus_{i = 0}^{k-j-1} \ \mathbb{C}.b^i.e_{\lambda_j} \oplus \mathbb{C}.b^{q_j}.e_{\lambda_j}$.\\
     Si \ $p_j + \cdots + p_{k-1} < k-j $ \ posons \ $V_j : = \oplus_{i = 0}^{k-j-1} \ \mathbb{C}.b^i.e_{\lambda_j}$.\\
     Alors on a
     $$ E_{\lambda_j} = P_j.E_{\lambda_j} \oplus V_j .$$
     \end{prop}
     
     \parag{Preuve} Commen{\c c}ons par remarquer que \ $E\big/F_j \simeq \A\big/\A.P_j $ \ est un th\`eme \ $[\lambda]-$primitif de rang \ $k-j$. On a donc, d'apr\`es lemme  \ref{Ext}  rappel\'e au d\'ebut de l'appendice,
     $$ dim_{\mathbb{C}}(Ext^1_{\A}(E\big/F_j, E_{\lambda_j})) - dim_{\mathbb{C}}(Hom_{\A}(E\big/F_j, E_{\lambda_j})) = k - j .$$
     Mais on a \ $Hom_{\A}(E\big/F_j, E_{\lambda_j}) \simeq \mathbb{C}$ \  ou \ $ \{0\} $ \ suivant que \ $\lambda_j \leq \lambda_k$ \ ou bien que \ $\lambda_j > \lambda_k$. En effet le seul quotient de rang 1 de \ $E\big/F_j$ \ est \ $E\big/F_{k-1} \simeq E_{\lambda_k}$ \ d'apr\`es le corollaire  \ref{sous-themes normaux}.\\
     On en d\'eduit que l'on a \ $dim_{\mathbb{C}}(Ext^1_{\A}(E\big/F_j, E_{\lambda_j}) = k-j +1 $ \ ou bien \ $k - j$ \ suivant que \ $\lambda_j \leq \lambda_k$ \ ou bien \ $\lambda_j > \lambda_k$.\\
     La r\'esolution \ $\A-$libre de \ $E\big/F_j \simeq \A\big/\A.P_j$ \ montre que l'espace vectoriel \ $Ext^1_{\A}(E\big/F_j, E_{\lambda_j})$ \ est isomorphe au conoyau de \ $P_j$ \ agissant sur \ $E_{\lambda_j}$. La codimension de \ $P_j.E_{\lambda_j}$ \ dans \ $E_{\lambda_j}$ \ est donc \ $k-j$ \ si \ $\lambda_j > \lambda_k$, et \ $k-j+1$ \ si \ $\lambda_j \leq \lambda_k$. Dans le premier cas, l'inclusion de \ $P_j.E_{\lambda_j}$ \ dans \ $b^{k-j}.E_{\lambda_j}$ \ suffit pour prouver notre assertion.\\
     Dans le cas \ $\lambda_j \leq \lambda_k$, ce qui \'equivaut \`a \ $p_j+ \cdots p_{k-1} \geq k-j$, il s'agit de montrer que toute combinaison lin\'eaire
     $$ \sum_{i=0}^{k-j-1} \ c_i.b^i.e_{\lambda_j} + \gamma.b^{q_j}.e_{\lambda_j} $$
     qui est dans \ $P_j.E_{\lambda_j}$ \ est nulle. L'inclusion \ $P_j.E_{\lambda_j} \subset b^{k-j}.E_{\lambda_j}$ \ montre d\'ej\`a que l'on doit avoir \ $c_i = 0 \quad \forall i \in [0,k-j-1]$. Il reste donc \`a montrer que \ $b^{q_j}.e_{\lambda_j} \not\in P_j.E_{\lambda_j}$.\\
     Pour cela remarquons d\'ej\`a que si \ $x \in E_{\lambda_j}$ \ est de valuation \ $b-$adique \'egale \`a \ $q$, alors \ $P_j.x$ \ sera de valuation \ $b-$adique exactement \ $q+k-j$ \ si \ $q$ \ n'est pas de la forme \ $p_j+ \cdots + p_{j+h} - (k-j)$ \ pour un entier \ $h \in [0,k-j-1]$. En effet, on peut ignorer les inversibles \ $S_{j+1}, \cdots, S_{k-1}$ \ qui ne changent pas la valuation \ $b-$adique, et constater qu'apr\`es l'action de \ $(a - \lambda_{j+h+1}.b)\dots (a- \lambda_k)$ \ ou bien on arrive  \`a une valuation exactement \'egale \`a \ $q + k - (j+h)$ \ ou bien la valuation finale ne sera pas \ $q + k-j$. L'action de \ $(a - \lambda_{j+h}.b)$ \ sur \ $b^{q+k-(j+h)}.e_{\lambda_j}$ \ donnera 
     $$ (q + \lambda_j + k - (j+h) - \lambda_{j+h}).b^{q+k-(j+h)+1}.e_{\lambda_j}$$
     et l'on a 
     \begin{align*}
     & (q + \lambda_j + k - (j+h) - \lambda_{j+h}) = q  -\big[ p_j + \cdots + p_{j+h-1} - h \big] + k - (j+h) \\
     & \qquad  =  q -\big[ p_j + \cdots + p_{j+h-1}\big] + k - j
     \end{align*}
     qui ne s'annule pas tant que \ $q \not=  p_j + \cdots + p_{j+h-1} - ( k-j) $.\\
     Montrons maintenant que \ $b^{q_j}.e_{\lambda_j} \not\in P_j.E_{\lambda_j}$. Pour cela raisonnons par l'absurde, et consid\'erons \ $x \in E_{\lambda_j}$ \ de valuation \ $b-$adique \'egale \`a \ $q \geq 0 $ \ et tel que \ $P_j.x = b^{q_j}.e_{\lambda_j}$. Si on a \ $q \not= p_j+ \cdots + p_{j+h} - (k-j)$ \ pour chaque \ $h \in [0,k-j-1]$, alors on aura \ $q_j = q + k-j $ \ ce qui contredit la d\'efinition de \ $q_j$. \\
     On a donc  pour un \ $h \in [0,k-j-1]$ \ tel que 
     \ $ q = p_j+ \cdots + p_{j+h} - (k-j) \geq 0  $ \ ce qui implique \ $q \geq q_j - (k-j)$.\\
     Mais si \ $q = q_j -(k-j)$ \ ceci contredit le fait que la valuation de \ $P_j.x$ \ soit \ $q_j$ \ d'apr\`es le calcul pr\'ec\'edent, et si on a \ $q > q_j$ \ la valuation de \ $P_j.x$ \ est strictement plus grande que \ $q_j$.  On a donc bien la contradiction d\'esir\'ee. $\hfill \blacksquare$ 
     
     \parag{Remarques}
     \begin{enumerate}
     \item Si on a \ $p_j \geq k-j$, alors on a \ $q_j = p_j$ \ et  \ $V_j = \oplus_{i = 0}^{k-j-1} \ \mathbb{C}.b^i.e_{\lambda_j} \oplus \mathbb{C}.b^{p_j}.e_{\lambda_j}$. Si \ $p_j \leq k-j-1 $ \ alors on a encore \ $\mathbb{C}.b^{p_j}.e_{\lambda_j} \subset V_j$. On a donc toujours \ $b^{p_j}.e_{\lambda_j} \in V_j$, pour chaque \ $j \in [1, k-1]$.
     \item Les sous-espaces vectoriels \ $V_j \subset E_{\lambda_j}$ \ sont ind\'ependants des inversibles \ $S_1, \cdots, S_{k-1}$. Ils sont d\'efinis uniquement \`a partir des invariants fondamentaux \ $\lambda_1, p_1, \cdots, p_{k-1}$ \ du th\`eme \ $[\lambda]-$primitif \ $E$. 
     \item Remarquons \'egalement que \ $P_j$ \ ne d\'epend que de \ $\lambda_{j+1}, \dots, \lambda_k$ \ et \ $S_{j+1}, \dots, S_{k-1}$, donc de l'id\'eal annulateur de la classe induite dans \ $E\big/F_j$ \  par le g\'en\'erateur fix\'e dans \ $E$.  $\hfill \square$
     \end{enumerate}

\bigskip

\subsubsection{Unicit\'e dans le cas stable.}

La proposition suivante est la clef du th\'eor\`eme d'unicit\'e.

\begin{prop}\label{Uni. stable}
Soit \ $E$ \ un th\`eme \ $[\lambda]-$primitif stable, et soit \ $e$ \ et \ $e'$ \ deux g\'en\'erateurs de \ $E$.
Soit \ $P _1 : =  (a -\lambda_2.b).S_2^{-1} \cdots S_{k-1}^{-1}.(a - \lambda_k.b) $ \ et supposons que 
\begin{enumerate}[i)]
\item \ $P : = (a - \lambda_1.b).S_1^{-1}.P_1 $ \ engendre l'id\'eal annulateur de \ $e$ \ dans \ $E$;
\item \ $P_1.e' = T_1.e_1$ \ o\`u \ $e_1 : = S_1^{-1}.P_1.e $;
\item \ $e - e' \in F_{k-1}$.
\end{enumerate}
Alors on a \ $T_1 - S_1 \in P_1.F_1$. En particulier, si \ $T_1.e_1$ \ et \ $S_1.e_1$ \ sont dans un m\^eme suppl\'ementaire de \ $P_1.F_1$ \ on aura \ $S_1 = T_1$.
\end{prop}

\parag{Preuve} Notons \ $[e]$ \  et \ $[e']$ \ les classes de \ $e$ \ et \ $e'$ \ dans \ $E\big/F_1$. Ce sont deux g\'en\'erateurs de ce th\`eme de rang \ $k-1$ \ ayant le m\^eme id\'eal annulateur  $\A.P_1$ \ dans \ $E\big/F_1$. Comme \ $[e - e'] \in F_{k-1}\big/F_1$, l'endomorphisme \ $\psi \in End_{\A}(E\big/F_1)$ \ d\'efini par \ $\psi([e]) = [e - e']$ \ est de rang \ $\leq k-2$, et il existe, puisque \ $E$ \ est stable, un \'el\'ement \ $\varphi \in End_{\A}(E)$ \ de \ rang \ $rang(\psi)+1 \leq k-1$ \ qui induit \ $\psi$. Posons \ $\varepsilon : = \varphi(e)$. On aura \ $P.\varepsilon = 0$ \ et comme \ $\varphi$ \ est de rang \ $\leq k-1$ \ on aura m\^eme \ $P_1.\varepsilon = 0$ \ d'apr\`es le lemme \ref{lemme manquant}.\\
Mais dire que \ $\varphi$ \ induit \ $\psi$ \ signifie que l'on a 
\begin{equation*}
\varepsilon = e - e' + U.e_1 \quad {\rm avec} \quad U \in \C[[b]] \tag{@}
\end{equation*}
puisque \ $F_1 = \C[[b]].e_1$. On en d\'eduit que
$$ P_1.\varepsilon = 0 = S_1.e_1 - T_1.e_1 + P_1.U.e_1 $$
ce qui prouve notre assertion. $\hfill \blacksquare$

\parag{Notation}Soient \ $\lambda_1, p_1, \dots, p_{k-1}$ \ les invariants fondamentaux d'un th\`eme \ $[\lambda]-$primitif. \\
Notons  \ $W_j$ \ l'ouvert affine de l'espace vectoriel \ $V_j \subset E_{\lambda_j}$ \ d\'efini par les deux conditions:  \ $ S_j(0) = 1$ \ et le coefficients de \ $b^{p_j}$ \ dans \ $S_j$ \ est non nul.

\bigskip

\begin{thm}\label{unicite}
Soit \ $E$ \ un th\`eme \ $[\lambda]-$primitif de rang \ $k \geq 2$. Si \ $E$ \  est stable on a unicit\'e de \ $P$ \ tel que \ $E \simeq \A\big/\A.P$ \ avec 
$$P : =  (a - \lambda_1.b).S_1^{-1}\cdots S_{k-1}^{-1}.(a - \lambda_k.b) \quad {\rm et} \quad  S_j \in W_j, \ S_j(0) = 1 \quad \forall j \in [1,k-1].$$
\end{thm}

\parag{D\'emonstration} C'est une r\'ecurrence imm\'ediate sur le rang du th\`eme stable consid\'er\'e en utilisant le fait que \ $E$ \ stable implique la stabilit\'e de \ $E\big/F_1$ \ et la proposition pr\'ec\'edente \ref{Uni. stable}. $\hfill \blacksquare$

\parag{Remarque} Le g\'en\'erateur \ $e$ \ d'annulateur \ $\A.P$ \ est unique modulo le groupe des automorphismes de \ $E$ \ qui est isomorphe au groupe des inversibles de l'alg\`ebre \ $\C[x]\big/(x^k)$. $\hfill \square$\\

\subsubsection{La propri\'et\'e d'unicit\'e.}

Le probl\`eme consistant \`a caract\'eriser les th\`emes primitifs poss\'edant cette propri\'et\'e d'unicit\'e est assez d\'elicat. Donnons d\'ej\`a un crit\`ere de non unicit\'e.

\bigskip

\begin{lemma}\label{instable 1}
Soit \ $E$ \ un th\`eme \ $[\lambda]-$primitif de rang \ $k \geq 3$ \ non stable mais tel que \ $E\big/F_1$ \ soit stable. Notons \ 
\begin{itemize}
\item $e$ \ un g\'en\'erateur de \ $E$, 
\item  $P_1 : = (a - \lambda_2.b).S_2^{-1}\cdots S_{k-1}^{-1}.(a - \lambda_k.b) $ \ le g\'en\'erateur de l'annulateur de \ $[e]$ \ dans \ $E\big/F_1$ \   et
\item  \ $P : = (a -\lambda_1.b).S_1^{-1}.P_1$ \ un g\'en\'erateur de l' annulateur de \ $e$ \ dans \ $E$. 
\end{itemize}
Soit  \ $e_{\lambda_1}$ \ le g\'en\'erateur standard de \ $F_1$. Alors il existe \ $e'$ \ un g\'en\'erateur de \ $E$ \ dont l'annulateur est \ $Q : = (a - \lambda_1.b).T_1^{-1}.P_1$, o\`u \ $T_1 \in \C[[b]]$ \ v\'erifie \ $T_1(0) = 1$ \ et  \ $(S_1 - T_1).e_{\lambda_1} \not\in P_1.F_1$.
\end{lemma}

\parag{Preuve} Comme nous avons suppos\'e \ $E\big/F_1$ \ stable, il existe \ $\psi \in End_{\A}(E\big/F_1)$ \ de rang \ $k-2$. Posons \ $\psi([e]) : = [\eta]$ ; alors \ $\eta \in F_{k-1}\setminus F_{k-2}$, et la relation \ $P_1.[e] = 0 $ \ dans \ $E\big/F_1$ \ donne que \ $P_1.\eta \in F_1$. Posons \ $P_1.\eta = Z_1.e_{\lambda_1}$, o\`u \ $Z_1 \in \C[[b]]$. Si on peut trouver \ $U \in \C[[b]]$ \ tel que \ $P_1.U.e_{\lambda_1} = Z_1.e_{\lambda_1}$, alors \ $\eta - U.e_{\lambda_1}$ \ qui est dans \ $F_{k-1}\setminus F_{k-2}$ \ puisque \ $k \geq 3$, v\'erifiera \ $P_1.(\eta - U.e_{\lambda_1}) = 0$ \ et \`a fortiori \ $P.(\eta - U.e_{\lambda_1}) = 0$, nous fournissant un endomorphisme de rang \ $k-1$ \ de \ $E$ \ ce qui contredit notre hypoth\`ese.\\
D'autre part \ $P_1\eta$ \ est dans \ $b.E$ \ car \ $\eta \in F_{k-1} \subset a.E + b.E$, et on a donc \ $Z_1(0) = 0$, puisque \ $F_1$ \ est normal. Posons  \ $ e' : = .e - \eta $ ; c'est un g\'en\'erateur de \ $E$, et il v\'erifie \ $P_1.e' = T_1.e_{\lambda_1}$ \ o\`u \ $T_1 : = S_1 - Z_1$. On a \ $T_1(0) = S_1(0) = 1$ \ et \ $(S_1 -T_1).e_{\lambda_1} = Z_1.e_{\lambda_1}$ \ n'est pas dans \ $P_1.F_1$, d'apr\`es ce qui pr\'ec\`ede. $\hfill \blacksquare$

\parag{Remarque} On notera que dans situation du lemme ci-dessus si \ $t \in \C$ \ le g\'en\'erateur \ $e_t : = t.e + (1-t).e' $ \ v\'erifiera \ $P_1.e_t = (t.S_1 + (1-t).T_1).e_{\lambda_1}$. Donc  si les coefficients de \ $b^{p_1}$ \ dans \ $S_1$ \ et \ $T_1$ \ \'etaient diff\'erents, on pourrait trouver \ $t \in \C$ \ tel que ce coefficient devienne nul. Mais ceci est impossible pour un th\`eme. Donc m\^eme si on trouve tout un sous-espace  affine de dimension \ $> 0$ \ de \ $S_1$ \ possibles, le coefficient (non nul) de \ $b^{p_1}$ \ est ind\'ependant   des choix. $\hfill \square$

\bigskip

 \begin{defn}\label{pte Uni.1}
   Soit \ $E$ \ un th\`eme \ $[\lambda]-$primitif de rang \ $k \geq 2$. On dira que \ $E$ \ a la propriet\'e U si
   on a unicit\'e de \ $P$ \ tel que \ $E \simeq \A\big/\A.P$ \ avec 
$$P : =  (a - \lambda_1.b).S_1^{-1}\cdots S_{k-1}^{-1}.(a - \lambda_k.b) \quad {\rm et} \quad  S_j \in W_j, \ S_j(0) = 1 \quad \forall j \in [1,k-1].$$
\end{defn}

\bigskip

On remarquera qu'en rang 1 et 2 tout th\`eme v\'erifie la propri\'et\'e  U.

\bigskip

Comme tout \ $E$ \ stable a cette propri\'et\'e nous allons explorer quels sont les th\`emes \ $[\lambda]-$primitifs {\em instables} (c'est-\`a-dire non stables) qui ont cependant cette propri\'et\'e. Nous verrons qu'il y en a peu.
   
   \begin{prop}\label{pte Uni.2}
   Si \ $E$ \ est  un th\`eme \ $[\lambda]-$primitif instable v\'erifiant la propri\'et\'e  U, alors pour tout \ $j \in [1,k-2]$ \ le th\`eme \  $E\big/F_j$ \ est instable et v\'erifie \'egalement la propri\'et\'e U. \\
   En particulier on a  \ $p_{k-1} = 0$.\\
   R\'eciproquement, si le quotient \ $E\big/F_1$ \ v\'erifie la propri\'et\'e U et v\'erifie\footnote{ce qui montre qu'il est instable.}  de plus   l'\'egalit\'e \ $End_{\A}(E\big/F_1) = \C.\id$,  alors \ $E$ \ v\'erifie la propri\'et\'e U ( et il est instable).
   \end{prop}
   
    \parag{Preuve} Le fait que \ $E\big/F_j$ \ v\'erifie la propri\'et\'e  U  si \ $E$ \ la v\'erifie est imm\'ediat. Si \ $E\big/F_1$ \ \'etait stable, le lemme \ref{instable 1} montrerait que \ $E$ \ ne v\'erifie pas la propri\'et\'e  U pourvu que le rang soit \ $\geq 3$. Par r\'ecurrence sur \ $j \in [1,k-2]$, on en d\'eduit que tous les \ $E\big/F_j$ \ sont instables et v\'erifient la propri\'et\'e U, pour \ $j \in [1,k-2]$. Le cas \ $j = k-2$ \ donne alors \ $p_{k-1} = 0$.\\
   Pour montrer  la r\'eciproque consid\'erons deux g\'en\'erateurs \ $e$ \ et \ $e'$ \ de \ $E$ \ v\'erifiant les propri\'et\'es suivantes :
   \begin{enumerate}[i)]
  \item  L'annulateur \ $\A.P$ \ de \ $e$ \ dans \ $E$ \ est de la forme donn\'ee dans le th\'eor\`eme \ref{unicite}.
  \item Les images de \ $e$ \ et \ $e'$ \ dans \ $E\big/F_1$ \ ont m\^eme annulateur \ $\A. P_1$.
  \item La diff\'erence \ $e - e'$ \ est dans \ $F_{k-1}$.
  \item On a  \ $P_1.e = S_1.e_{\lambda_1}$ \ et \ $P_1.e' = T_1.e_{\lambda_1}$ \ avec \ $S_1(0) = T_1(0) = 1$.
  \end{enumerate}
  Il s'agit alors de montrer que l'on a \ $S_1 - T_1 \in P_1.F_1$. L'endomorphisme de \ $E\big/F_1$ \ donn\'e en envoyant \ $[e]$ \ sur \ $[e-e']$ \ n'est pas surjectif, puisque \ $[e - e'] \in F_{k-1}\big/F_1$. Il est donc nul d'apr\`es notre hypoth\`ese,  ce qui signifie que \ $e - e' \in F_1$;  ceci  donne la conclusion cherch\'ee. $\hfill \blacksquare$  
  
  \bigskip 

 Un corollaire facile d\'ecrit compl\`etement la situation en rang 3.
   
   \begin{cor}\label{rg 3 et pte U}
   Les th\`emes \ $[\lambda]-$primitifs de rang 3 v\'erifiant la propri\'et\'e  U  sont les th\`emes stables et ceux qui v\'erifient \ $p_2 = 0$ \  qui sont n\'ecessairement instables.
   \end{cor}
   
   \parag{Preuve} Comme tout th\`eme de rang 2 v\'erifie la propri\'et\'e  U, si \ $E$ \ est un th\`eme \ $[\lambda]-$primitif de rang 3, qui est instable et v\'erifie la propri\'et\'e U, alors \ $E\big/F_1$ \ est instable donc isomorphe \`a  \ $E_{\lambda.\lambda}$ \ ce qui impose \ $p_2 = 0$. Mais r\'eciproquement, si on a \ $p_2 = 0$, alors \ $E\big/F_1$ \ est isomorphe \`a  \ $E_{\lambda.\lambda}$ \ qui v\'erifie la condition \ $End_{\A}(E\big/F_1)  \simeq \C.\id $ \ de la proposition pr\'ec\'edente. Donc tout th\`eme primitif de rang 3 v\'erifiant \ $p_2 = 0$ \ v\'erifie la propri\'et\'e  U. $\hfill \blacksquare$

\bigskip

Terminons par le cas extr\^eme o\`u \ $p_1 = \dots = p_{k-1} = 0$.

\smallskip

\begin{lemma}\label{unicite speciale}
Soit \ $E$ \ un th\`eme  \ $[\lambda]-$primitif de rang \ $k$. Supposons que l'on ait \ $p_1 = \dots = p_{k-1} = 0 $. Alors on a unicit\'e de \ $P$ \ tel que \ $E \simeq \A\big/\A.P$ \ avec 
$$P : =  (a - \lambda_1.b).S_1^{-1}\cdots S_{k-1}^{-1}.(a - \lambda_k.b) \quad {\rm et} \quad  S_j \in W_j, \ S_j(0) = 1 \quad \forall j \in [1,k-1].$$
\end{lemma}

\parag{Preuve} Prouver l'assertion suivante par r\'ecurrence sur le rang \ $k$  \  est une cons\'equence facile de la remarque iv)  qui suit le corollaire \ref{stable 4} et de la proposition \ref{pte Uni.2}:
\begin{itemize}
\item Soit \ $E$ \ un th\`eme  \ $[\lambda]-$primitif de rang \ $k$ \ v\'erifiant \ $p_1= \dots = p_{k-1} = 0$. Alors il v\'erifie la propri\'et\'e U et l'\'egalit\'e \ $End_{\A}(E) = \C.\id$.
\end{itemize}
$\hfill \blacksquare$
 
\subsection{ Th\`emes stables g\'en\'eraux.}

Nous allons \'etendre une partie des consid\'erations pr\'ec\'edentes aux th\`emes g\'en\'eraux.

\begin{lemma}\label{hom(E,E)}
Soit \ $E$ \ un th\`eme de rang \ $k$. Alors l'espace vectoriel \ $Hom_{\A}(E,E)$ \ est de dimension au plus \'egale \`a \ $k$.
\end{lemma}

\parag{Preuve} Comme le r\'esultat est connu dans le cas \ $[\lambda]-$primitif, nous pouvons faire une r\'ecurrence sur le cardinal \ $q$ \  de l'ensemble  \ $Exp(E)$, le cas \ $q =1$ \ \'etant acquis. Supposons  le r\'esultat connu pour \ $q\geq 1$ \ et montrons-le pour \ $q+1$. Soit donc \ $E$ \ un th\`eme avec \ $Card\{Exp(E)\} = q+1 $ \ et fixons \ $[\lambda] \in Exp(E)$. On a une suite exacte
\begin{equation*}
0 \to E[\not= \lambda] \to E \to E\big/[\lambda] \to 0 \tag{1} 
\end{equation*}
o\`u \ $E\big/[\lambda]$ \ est \ $[\lambda]-$primitif, et o\`u \ $Card\{ Exp(E[\not=\lambda]\} = q $. On a alors le diagramme d\'eduit de \ $(1)$\\
$$\xymatrix{0 \ar[r] &Hom_{\A}(E,E[\not=\lambda]) \ar[r] \ar[d]^{\simeq} &  Hom_{\A}(E,E) \ar[r] &  Hom_{\A}(E, E\big/[\lambda]) \ar[r] \ar[d]^{\simeq} & 0 \\
\quad & Hom_{\A}(E[\not=\lambda],E[\not=\lambda]) & \quad &Hom_{\A}(E\big/[\lambda], E\big/[\lambda]) & \quad }$$
et de l'additivit\'e de la dimension  permet de conclure gr\^ace \`a l'hypoth\`ese de r\'ecurrence. $\hfill \blacksquare$

\bigskip

\begin{prop}\label{stable gen. 0}
Soit \ $E$ \ un th\`eme de rang \ $k$ ; les propri\'et\'es suivantes sont \'equivalentes :
\begin{enumerate}[1)]
\item Il existe une injection \ $\A-$lin\'eaire de \ $E$ \ dans \ $\Xi$ \ dont l'image est invariante par la monodromie \ $\mathcal{T}$.
\item Il existe un unique sous-th\`eme de \ $\Xi$ \ isomorphe \`a \ $E$.
\item L'espace vectoriel \ $Hom_{\A}(E,E)$ \ est de dimension \ $k$.
\end{enumerate}
\end{prop}

\parag{Preuve} L'implication \ $2) \Rightarrow 1)$ \ est claire car \ $\mathcal{T} \circ j$ \ est une  injection \ $\A-$lin\'eaire de \ $E$ \ dans \ $\Xi$ \ d\`es que \ $j$ \ est  injection \ $\A-$lin\'eaire de \ $E$ \ dans \ $\Xi$. \\
Montrons \ $3) \Rightarrow 2) $. Soit \ $j : E \to \Xi $ \ une  injection \ $\A-$lin\'eaire de \ $E$ \ dans \ $\Xi$. La composition avec \ $j$ \ donne une injection \ $\C-$lin\'eaire
$$ Hom_{\A}(E, E) \to Hom_{\A}(E, \Xi) .$$
Comme ces deux espaces vectoriels ont m\^eme dimension \ $k$, le premier par hypoth\`ese, le second en vertu du th\'eor\`eme 2.2.1  de [B.05], c'est une bijection. En particulier toute  injection \ $\A-$lin\'eaire de \ $E$ \ dans \ $\Xi$ \ a son image contenue dans \ $j(E)$, ce qui donne  la propri\'et\'e 2).\\
Montrons enfin que \ $1) \Rightarrow 3)$. Pour cela montrons d'abord  par r\'ecurrence sur l'entier  \\ 
$q : = Card\{Exp(E)\}$ \ que si \ $E$ \ v\'erifie 3) alors   le polyn\^ome minimal de l'action de \ $\mathcal{T}$ \ sur \ $E$ \ est de degr\'e exactement le rang de \ $E$. L'assertion \'etant connue pour \ $q = 1$, d'apr\`es la proposition \ref{stable 0} et le fait que \ $\mathcal{T} - exp(2i\pi.\lambda).\id$ \ induise un endomorphisme de rang \ $k-1$,  supposons-la v\'erifi\'ee pour \ $q \geq1$ \ et montrons-l\`a pour \ $Card\{Exp(E)\} = q+1$, en reprenant les notations utilis\'ees dans la preuve du lemme ci-dessus.\\
La suite exacte \ $(1)$ \ montre d\'ej\`a, gr\^ace \`a l'hypoth\`ese de r\'ecurrence, que 3) \ est v\'erifi\'ee pour \ $E[\not=\lambda]$ \ et aussi pour \ $E\big/[\lambda]$ \ puisque qu'en composant l'injection consid\'er\'ee de \ $E$ \ dans \ $\Xi$ \ avec le quotient \ $ \Xi \to (\Xi\big/\Xi_{\not=\lambda}) \simeq \Xi_{\lambda}$ \ on obtient une injection de \ $E\big/[\lambda]$ \ dans \ $\Xi_{\lambda}$ \ qui est stable par \ $\mathcal{T}$.\\
On conclut alors en remarquant que le polyn\^ome minimal de \ $\mathcal{T}$ \ agissant sur \ $E$\footnote{remarquer que comme \ $Hom_{\A}(E,E)$ \ est de dimension au plus \'egale \`a \ $k$ \ ce polyn\^ome minimal est de degr\'e au plus \ $k$.} divise les polyn\^omes minimaux de \ $\mathcal{T}$ \ agissant sur \ $E[\not=\lambda]$ \ et \ $E\big/[\lambda]$ \ respectivement. Comme ils sont premiers entre eux, il divise le produit qui, gr\^ace \`a l'hypoth\`ese de r\'ecurrence est de degr\'e \ $rg(E[\not=\lambda]) + rg(E\big/[\lambda]) = rg(E)$. $\hfill \blacksquare$

\bigskip

\begin{defn}\label{stable gen. 1}
On dira qu'un th\`eme est {\bf stable} s'il v\'erifie les propri\'et\'es 1) , 2), 3)  de la proposition pr\'ec\'edente.
\end{defn}

\parag{Remarques}
 \begin{enumerate} [1)]
 \item Il est clair que cette d\'efinition est compatible avec celle donn\'ee dans le cas \ $[\lambda]-$primitif.
 \item Le dual d\'ecal\'e \ $E^*\otimes E_{\delta}$ \ pour \ $\delta \in \mathbb{Q}$ \ assez grand d'un th\`eme stable est encore un th\`eme stable : en effet d\`es que le d\'ecalage sera suffisant pour avoir un (a,b)-module monog\`ene g\'eom\'etrique, la  condition 3) de la proposition pr\'ec\'edente sera v\'erifi\'ee, puisque \ $\varphi \mapsto \varphi^*\otimes \id $ \ donne un isomorphisme lin\'eaire  de \ $End_{\A}(E)$ \ sur \ $End_{\A}(E^*\otimes E_{\delta})$.
  \end{enumerate}
 
 \begin{lemma}
 Tout sous-th\`eme normal et tout th\`eme quotient d'un th\`eme stable est stable.
 \end{lemma}
 
 \parag{Preuve} D'apr\`es la remarque 2) ci-dessus il suffit de traiter le cas des quotients. Par r\'ecurrence sur le rang du sous-th\`eme normal par lequel on quotiente, on se ram\`ene au cas o\`u l'on quotiente par un sous-th\`eme normal de rang \ $1$. Mais si \ $E_{\lambda+p}$ \ est un sous-th\`eme normal de rang \ $1$ \ de \ $E \subset \Xi$ \ stable par \ $\mathcal{T}$, l'image de \ $E$ \ par \ $f_{\lambda}\oplus \id_{\mu \not= \lambda} : \Xi \to \Xi$ \ de \ $E$ \ est un sous-th\`eme de \ $\Xi$ \ isomorphe \`a \ $E\big/E_{\lambda + p}$ \ qui est stable par \ $\mathcal{T}$ \ car \ $f_{\lambda}$ \ commute \`a \ $\mathcal{T}$. 
 Donc \ $E\big/E_{\lambda+ p}$ \ est stable. $\hfill \blacksquare$
 
 \parag{Remarques}
 \begin{enumerate}[1)]
 \item Si chaque partie coprimitive d'un th\`eme est stable, alors le th\`eme est stable : en effet, il suffit de composer les injections des parties $[\lambda]-$coprimitives avec les quotients et de prendre la somme directe des morphismes dans  les \ $\Xi_{\lambda}$ \ ainsi obtenus pour avoir une injection dont l'image dans  \ $\Xi$ \  est stable par la monodromie; en effet la monodromie de \ $\Xi$ \ se d\'ecompose en somme directe des monodromies de chacun des \ $\Xi_{\lambda}$.
 \item Un argument simple de dualit\'e montre, \`a partir de la remarque pr\'ec\'edente que si chaque partie \ $[\lambda]-$primitive d'un th\`eme est stable, alors le th\`eme est stable.
 \end{enumerate}
 
\section{Familles holomorphes de th\`emes \ $[\lambda]-$primitifs.}

\subsection{D\'efinitions et premiers exemples.}

\subsubsection{D\'efinitions.}

Soit \ $X$ \ un espace complexe. Nous noterons \ $\mathcal{O}_X[[b]]$ \ le faisceau sur \ $X$ \ d\'efini par le pr\'efaisceau
$$ U \to \mathcal{O}_X(U)[[b]] .$$
C'est un faisceau de \ $\mathcal{O}_X-$alg\`ebres. Pour \ $\mathcal{J} \subset (\mathcal{O}_X)^p$ \ un sous-faisceau de \ $\mathcal{O}_X-$modules (resp. $\mathcal{O}_X-$coh\'erent) de \ $(\mathcal{O}_X)^p$, on notera \ $\mathcal{J}[[b]]$ \ le sous-faisceau de \ $\mathcal{O}_X[[b]]-$modules  (resp. $\mathcal{O}_X[[b]]-$coh\'erent) de \ $(\mathcal{O}_X[[b]])^p$ \ qui est engendr\'e par \ $\mathcal{J}$.

\bigskip

On notera que pour \ $X$ \ de Stein on a le th\'eor\`eme B de Cartan pour le faisceau \ $\mathcal{O}_X[[b]]$\footnote{et m\^eme pour tout \ $\mathcal{O}_X[[b]]$-module coh\'erent. }, \`a savoir que \ $H^i(X, \mathcal{O}_X[[b]]) = 0, \forall i \geq 1$.

\bigskip 

\begin{defn}\label{Faisc. (a,b)}
Soit \ $X$ \ un espace complexe. Un faisceau de \ $\mathcal{O}_X-$(a,b)-modules \ $\mathbb{E}$ \ sur \ $X$ \ est la donn\'ee d'un faisceau localement libre de type fini  de \ $\mathcal{O}_X[[b]]-$modules muni d'un morphisme de faisceaux 
$$ a : \mathbb{E} \to \mathbb{E} $$
qui est \ $\mathcal{O}_X-$lin\'eaire, continu pour la topologie \ $b-$adique de \ $\mathbb{E}$, et satisfait la relation de commutation \ $ a.b - b.a = b^2 $.\\
Un morphisme entre deux faisceaux de \ $\mathcal{O}_X-$(a,b)-modules sur \ $X$ \ sera un morphisme de faisceaux de \ $\mathcal{O}_X[[b]]-$modules qui commute aux actions respectives de \ $a$.
\end{defn}

\parag{Exemple}
 Pour \ $\lambda \in \mathbb{Q}^{+*}$ \ et \ $N$ \ entier consid\'erons
$$ \Xi_{X,\lambda}^{(N)} : = \oplus_{j=0}^{N} \mathcal{O}_X[[b]].e_{\lambda,j}\quad {\rm avec} \quad e_{\lambda,j} : = s^{\lambda-1}\frac{Log\, s)^j}{j!}  $$
muni de l'op\'eration \ $a$ \ d\'efinie par r\'ecurrence sur \ $j \geq 0$ \  de la fa{\c c}on suivante :
\begin{align*}
& a.e_{\lambda, 0} = \lambda b.e_{\lambda, 0} \\
& a.e_{\lambda, j} = \lambda b.e_{\lambda, j} + b.e_{\lambda, j-1}, \quad \forall j \geq 1.
\end{align*}

\parag{Notation} On notera \ $\Xi_{X,\lambda}^{(N)} : = \mathcal{O}_X \ptc \ \Xi_{\lambda}^{(N)} $.
On  notera aussi \ $\Xi_{X,\lambda} : =  \bigcup_{N \in \mathbb{N}} \Xi_{\lambda}^{(N)}$ \ et \ $\Xi_{X} : = \oplus_{\lambda \in ]0,1] \cap \mathbb{Q}} \ \Xi_{X,\lambda} $.\\
On remarquera que \ $\Xi_{X,\lambda}^{(N)}$ \ est \`a p\^ole simple, c'est-\`a-dire que \ $a.\Xi_{X,\lambda}^{(N)} \subset b.\Xi_{X,\lambda}^{(N)}$. $\hfill \square$

\bigskip

Soit \ $x \in X$ \ un point (ferm\'e). On a un morphisme d'\'evaluation en \ $x$ \ des fonctions holomorphes
$$ \mathcal{O}(X) \to \mathbb{C}\simeq  \mathcal{O}(X)\big/\mathcal{M}_x  $$
o\`u \ $\mathcal{M}_x \subset \mathcal{O}_X$ \ est le sous-faisceau des fonctions holomorphes nulles en \ $x$.\\
Si \ $\mathbb{E}$ \ est un faisceau de  \ $\mathcal{O}_X-$(a,b)-modules sur \ $X$, on aura, de fa{\c c}on analogue une application d'\'evaluation en \ $x$
$$ \mathbb{E} \to  E(x) : = \mathbb{E}\big/\mathcal{M}_x[[b]].\mathbb{E} $$
o\`u \ $E(x)$ \ est un (a,b)-module qui sera appel\'e la fibre en \ $x$ \ du faisceau \ $\mathbb{E}$.\\
On consid\'erera un faisceau de (a,b)-modules sur \ $X$ \ comme une famille de (a,b)-modules param\'etr\'ee par \ $X$.

\parag{Convention} Nous appellerons application holomorphe d'un espace complexe \ $X$ \ \`a valeurs dans \ $\Xi_{\lambda}^{(N)}$ (resp. \ $\Xi_{\lambda} : = \bigcup_{N \in \mathbb{N}} \Xi_{\lambda}^{(N)}$, \ $\Xi$) \ une section globale du faisceau \ $\Xi_{X,\lambda}^{(N)} $ \ (resp. \ $\Xi_{X,\lambda}$, $\Xi_X$). $\hfill \square$

\bigskip

\begin{defn}\label{k-thematique}
Une application holomorphe \ $\varphi : X \to \Xi$ \ d'un espace complexe \ $X$ \`a valeurs dans \ $\Xi$ \ sera dite {\bf k-th\'ematique} si la condition suivante est satisfaite:
\begin{itemize}
\item Le sous-$\mathcal{O}_X[[b]]-$module \ $\mathbb{E}_{\varphi}$ \ de \ $\Xi_X$ \ engendr\'e par les  \ $ a^{\nu}.\varphi , \ \nu \in \mathbb{N}$ \ est libre de rang \ $k$ \ et de base \ $\varphi, a.\varphi, \dots, a^{k-1}.\varphi$. 
\end{itemize}
\end{defn}

\bigskip

 Pour chaque \ $x \in X$ \ notons \ $ E(\varphi(x)) : = \mathbb{E}_{\varphi}\big/\mathcal{M}_x.\mathbb{E}_{\varphi} \simeq \A.\varphi(x) \subset \Xi$. C'est un th\`eme de rang k. On a, de plus, la restriction suivante:
 
 \bigskip

\begin{lemma} \label{Bernstein constant}
Soit \ $X$ \ un espace complexe r\'eduit et soit  \ $\varphi : X \to \Xi$ \  une application holomorphe k-th\'ematique;  le polyn\^ome de Bernstein \ $B_{\varphi(x)}$ \ de \ $E(\varphi(x))$ \ est localement constant sur \ $X$.
\end{lemma}

\parag{Preuve} \'Ecrivons sur \ $X$ :
$$ a^k.\varphi = \sum_{j=0}^{k-1} \quad S_{k-j}.a^j.\varphi $$
o\`u \ $S_1, \cdots, S_k$ \ sont des sections sur \ $X$ \ du faisceau \ $\mathcal{O}_X[[b]]$. Comme pour chaque \ $x \in X$ \ $\A.\varphi(x) \subset \Xi$ \ est un th\`eme de rang k, son \'el\'ement de Bernstein est donn\'e par \ $a^k - \sum_{j=0}^{k-1} \sigma_{k-j}(x).b^{k-j}.a^j $, o\`u \ $\sigma_{k-j}(x)$ \ est le coefficient de \ $b^{k-j}$ \ dans \ $S_{k-j}(x)$. On notera que \ $ \sigma_{k-j}(x).b^{k-j}$ \ est la forme initiale de \ $S_{k-j}(x)$ \ quand celle-ci n'est pas de degr\'e strictement plus grand que \ $ k-j$.\\
Mais \ $x \mapsto \sigma_{k-j}(x)$ \ est une fonction holomorphe sur \ $X$ qui ne prend que des valeurs dans \ $\mathbb{Q}$. Elle est donc localement constante. $\hfill \blacksquare$

\bigskip

\begin{cor}\label{inv. fond. loc. cst.}
Soit \ $X$ \ un espace complexe r\'eduit et soit  \ $\varphi : X \to \Xi_{\lambda}$ \ une application holomorphe k-th\'ematique, les invariants fondamentaux \ $\lambda_1, p_1, \dots, p_{k-1}$ \ des th\`emes \ $\A.\varphi(x) \subset \Xi_{\lambda}$ \ sont localement constants sur \ $X$.
\end{cor}

\parag{Preuve} Comme pour un th\`eme \ $[\lambda]-$primitif de rang \ $k$ \  les racines du polyn\^ome de Bernstein sont les nombres \ $k - (\lambda_j+j)$ \ et que la suite des \ $\lambda_j + j$ \ est croissante, le fait que le polyn\^ome de Bernstein soit localement constant sur \ $X$ \ implique la locale constance des invariants fondamentaux.$\hfill \blacksquare$

\parag{Remarque} M\^eme quand \ $X$ \ est r\'eduit, il ne suffit pas, en g\'en\'eral,  de v\'erifier que pour chaque \ $x \in X$ \ le (a,b)-module \ $E(\varphi(x))$ \ est un th\`eme de rang \ $k$ \ pour satisfaire la condition de la d\'efinition \ref{k-thematique}  comme le montre l'exemple suivant :\\
Soit \ $\lambda \in \mathbb{Q} , \lambda > 1$, et posons pour \ $z \in \mathbb{C}$ :
 $$\varphi(z) : = s^{\lambda-1}.Log\,s + (z + b).s^{\lambda-2} = s^{\lambda-1}.Log\,s + z.s^{\lambda-2} + \frac{1}{\lambda -1}.s^{\lambda-1} .$$
 Alors l'\'el\'ement de Bernstein de \ $E(z) : = \A.\varphi(z) $ \ pour \ $z \not= 0$ \ est \ $(a-\lambda.b)(a - \lambda.b)$ \ alors que l'\'el\'ement de Bernstein de \ $E(0)$ \ vaut \ $(a -(\lambda+1).b)(a - \lambda.b)$. On conclut gr\^ace au lemme pr\'ec\'edent.$\hfill \square$

\bigskip

\parag{Exemples}
Soit \ $\varphi : X \to \Xi_{\lambda}^{(k-1)} = \oplus_{j=0}^{k-1} \ \mathbb{C}[[b]].e_{\lambda,j} $ \ une application holomorphe. Supposons que le coefficient de \ $e_{\lambda,k-1}$ \ soit de la forme \ $b^n.S(b,x)$ \ o\`u \ $S$ \ est un inversible de l'alg\`ebre\footnote{ce qui revient \`a dire que l'\'el\'ement de \ $\mathcal{O}(X)$ \ qui est le terme constant en \ $b$ \ de \ $S$ \ est un inversible de \ $\mathcal{O}(X)$.} \ $\mathcal{O}(X)[[b]]$, et que la valuation en \ $b$ \ de \ $\varphi - b^n.S.e_{\lambda,k-1}$ soit strictement plus grande que \ $n$. Alors  le sous-faisceau
$$ \mathbb{E}_{\varphi} : = \sum_{i=0}^{k-1} \ \mathcal{O}_X[[b]].a^i.\varphi \ \subset \ \Xi_{X,\lambda}^{(k-1)} $$
est libre de rang \ $k$ \ sur \ $\mathcal{O}_X[[b]]$, stable par \ $a$. \\
En effet, on se ram\`ene imm\'ediatement au cas o\`u \ $S \equiv 1$, et on constate alors que
$$\psi : =  (a - (\lambda + n).b).\varphi $$
v\'erifie la m\^eme hypoth\`ese que \ $\varphi$ \ en rempla{\c c}ant \ $k$ \ par \ $k-1$ \ et \ $n$ \ par \ $n+1$. On conclut par une r\'ecurrence facile.\\
On notera que dans cet exemple (qui est bien particulier) on a \ $\lambda_k = \lambda + n$ \ puis \ $ \lambda_{k-1} = \lambda + n + 1, \dots $ \ ce qui signifie que \ $p_1 = p_2 = \dots = p_{k-1} = 0$. 
\\
Le lecteur trouvera dans l'appendice \ref{App.2} dans le corollaire \ref{exist. k-thm.} une m\'ethode g\'en\'erale et syst\'ematique pour construire des applications \ $k-$th\'ematiques. $\hfill \square$

\bigskip

\begin{defn}\label{Familles holomorphes}
Soit \ $X$ \ un espace complexe r\'eduit et soit \ $\mathbb{E}$ \ un faisceau de (a,b)-modules sur \ $X$. Nous dirons que \ $\mathbb{E}$ \ est une {\bf famille holomorphe de th\`emes  de rang $k$ \ param\'etr\'ee par} \ $X$ \ si la condition suivante est remplie :
\begin{itemize}
\item Il existe un recouvrement ouvert \ $(\mathcal{U}_{\alpha})_{\alpha \in A}$ \ de \ $X$ \ et pour chaque \ $\alpha \in A$ \ une application holomorphe th\'ematique 
$$\varphi_{\alpha} : \mathcal{U}_{\alpha} \to \Xi $$
et un isomorphisme de faisceaux de \ $\mathcal{O}_{\mathcal{U}_{\alpha}}[[b]]-$modules
 $$\mathbb{E}_{\vert \mathcal{U}_{\alpha}} \simeq \mathbb{E}_{\varphi_{\alpha}}$$
 compatible aux \ $\A-$structures.
 \end{itemize}
\end{defn}

\parag{Remarques}
\begin{enumerate}[i)]
\item Dans une famille holomorphe de th\`emes de rang \ $k$, le polyn\^ome de Bernstein est localement constant d'apr\`es \ref{Bernstein constant}. 
\item  Si, de plus,  $E(x)$ \ est un th\`eme \ $[\lambda]-$primitif  pour chaque \ $x \in X$,  les invariants fondamentaux sont localement constants sur \ $X$ \ d'apr\`es \ref{inv. fond. loc. cst.}.
\item Quand on consid\`ere une famille holomorphe de th\`emes \ $[\lambda]-$primitifs de \\
rang k on peut supposer que chaque application \ $\varphi_{\alpha}$ \ est \`a valeurs dans \ $\Xi_{\lambda}^{(k-1)}$ \ o\`u \ $[\lambda] \, \cap \, ]0,1] = \{\lambda \}$. $\hfill \square$
\end{enumerate}

\subsubsection{Premiers exemples : Famille holomorphes de th\`emes \\ $[\lambda]-$primitifs de rang 1 et 2.}

Le cas du rang 1 se d\'eduit de la remarque simple suivante :\\
Soit \ $X$ \ un espace complexe r\'eduit et connexe et soit \ $\varphi : X \to \Xi_{\lambda}^{(N)} $ \ une application \ $1-$th\'ematique. Alors il existe une section \ $S_0 \in \Gamma(X, \mathcal{O}_X[[b]])$ \ qui est inversible et un entier \ $n$ \ tel que \ $S_0.\varphi = s^{\lambda+n-1}$. On en d\'eduit que le faisceau \ $\mathbb{E}_{\varphi}$ \ est  isomorphe au faisceau \ $\mathcal{O} \ptc E_{\lambda+n}$ \ et donc \`a l'application holomorphe constante \ $X \to E_{\lambda+n}$ \ dont la valeur en chaque point est un g\'en\'erateur standard \ $e_{\lambda+n}$ \ de \ $E_{\lambda+n}$ \ qui v\'erifie \ $a.e_{\lambda+n} = (\lambda+n).b.e_{\lambda+n}$.

\bigskip

Le cas du rang 2 est d\'ecrit par  la proposition et le lemme  qui suivent. 

\begin{prop}\label{universelle rg 2}
Fixons \ $[\lambda] \in \mathbb{Q}\big/\mathbb{Z}, \lambda_1 > 1, \lambda_1 \in [\lambda]$ \ ainsi que \ $p \in \mathbb{N}^*$. Soit  \ $X$ \ un espace complexe r\'eduit connexe et soit \ $\varphi : X \to \Xi_{\lambda}$ \ une application holomorphe 2-th\'ematique, telle que les invariants fondamentaux des th\`emes associ\'es soient \ $\lambda_1 $ \ et \ $p_1 : = p$. Alors il existe deux  applications holomorphes \ $\alpha, \beta : X \to \C^*$ \ telles  que l'on ait, pour chaque \ $x \in X$ \ l'\'egalit\'e
$$ \A.\varphi(x) = \A.\psi(x) $$
o\`u 
\begin{equation*}
\psi(x) : = s^{\lambda_1+p-2}.Log\,s + \beta(x).(1 + \alpha(x).b^p).s^{\lambda_1-1}. \tag{@}
\end{equation*}
\end{prop}

\parag{Preuve} Il n'est pas restrictif de supposer que l'on a
\begin{equation*}
\varphi(x)  =  s^{\lambda_1+p-2}.Log\,s +  \Sigma(x).s^{\lambda_1-2} \tag{1}
\end{equation*}
o\`u \ $\Sigma : X \to \C[[b]]$ \ est holomorphe et \ $\Sigma(x)$ \ est un inversible de \ $\C[[b]]$ \ pour chaque \ $x \in X$. Ceci r\'esulte du fait que le coefficient de \ $ s^{\lambda_1+p-2}.Log\,s$ \ doit \^etre un inversible de \ $\C[[b]]$ \ d\'ependant holomorphiquement de \ $X$, et que les autres termes dans \ $\varphi(x)$ \ doivent \^etre dans \ $\C[[b]].s^{\lambda_1 - 2}$ \ pour avoir un th\`eme de rang 2 avec \ $\lambda_2 = \lambda_1+p-1$. La d\'efinition de  \ $\lambda_1 \in [\lambda]$ \ impose alors l'inversibilit\'e de \ $\Sigma(x)$ \ pour tout \ $x \in X$.\\
On d\'eduit de \ $(1)$ \ la relation
\begin{align*}
& (a - (\lambda_1+p-1).b).\varphi(x) = \frac{s^{\lambda_1+p-1}}{\lambda_1+p-1} + \Sigma(x).s^{\lambda_1-1} +  b^2.\Sigma(x)'.s^{\lambda_1-2} + \\
& \qquad \qquad \qquad \qquad \qquad \qquad \qquad  - (\lambda_1+p-1).b.\Sigma(x).s^{\lambda_1-2}  \\
&\qquad \qquad  \qquad \qquad \qquad  = \big(b^2.\Sigma(x)' - p.b.\Sigma(x) + \gamma.b^{p+1}\big).s^{\lambda_1-2}
\end{align*}
o\`u \ $\gamma : = (\lambda_1-1).\lambda_1\dots (\lambda_1+p-2)$ \ et o\`u l'on a not\'e \ $\Sigma(x)'$ \ la d\'eriv\'ee en \ $b$ \ de \ $\Sigma(x) \in \C[[b]]$. On a donc 
 $$(a - (\lambda_1+p-1).b).\varphi(x) = S(x).s^{\lambda_1-1} \quad {\rm avec} \quad (\lambda_1-1).S(x) : = b.\Sigma(x)' - p.\Sigma(x) + \gamma.b^p .$$
Posons \ $S(x)  = S_0(x) + S_p(x).b^p + b.\tilde{S}(x) $, o\`u \ $S_0(x) \in \C$ \ et o\`u \ $\tilde{S}(x) \in \C[[b]]$ \ n'a plus de terme en \ $b^{p-1}$. On peut alors trouver une application holomorphe \ $T : X \to \C[[b]]$ \ v\'erifiant :
$$ b^2.T(x)' - (p-1).b.T(x) = b.\tilde{S}(x) .$$
Soit \ $\psi : X \to \Xi_{\lambda}^{(1)} $ \ l'application holomorphe d\'efinie en posant
$$ \psi(x) : = \varphi(x) - T(x).s^{\lambda_1-1} .$$
Comme \ $\A.\varphi(x)$ \ contient \ $\C[[b]].s^{\lambda_1-1}$ \ pour chaque \ $x \in X$, ce qui se traduit par l'inversibilit\'e de la fonction holomorphe \ $S_0 : X \to \C$, qui se d\'eduit de l'inversibilit\'e de \ $\Sigma$ \ dans \ $\mathcal{O}(X)[[b]]$, on aura l'\'egalit\'e des th\`emes \ $\A.\varphi(x)$ \ et \ $\A.\psi(x)$ \ pour chaque \ $x \in X$.\\
Mais par construction on a
$$ (a - (\lambda_1+p-1).b).\psi(x) = (S_0(x) + S_p(x).b^p).s^{\lambda_1-1} .$$
Comme \ $S_0(x) \not= 0 $ \ et \ $S_p(x) \not= 0 $ \ pour chaque \ $x \in X$ \  on peut finalement d\'efinir, gr\^ace aux relations \ $ (\lambda_1 - 1).S_p(x) = \gamma$ \ et \ $(\lambda_1 - 1).S_0(x) = -p.\Sigma_0(x)$ : 
 $$\alpha(x) : = -\frac{\gamma}{p.\Sigma_0(x)}, \quad \beta(x) : = S_0(x)$$
ce qui donne l'identit\'e \ $(@)$. $\hfill \blacksquare$

\begin{lemma}\label{rg 2 p =0 }
Fixons \ $[\lambda] \in \mathbb{Q}\big/\mathbb{Z}, \lambda_1 > 1, \lambda_1 \in [\lambda]$. Soit  \ $X$ \ un espace complexe r\'eduit connexe et soit \ $\varphi : X \to \Xi_{\lambda}$ \ une application holomorphe 2-th\'ematique, telle que les invariants fondamentaux des th\`emes associ\'es soient \ $\lambda_1 $ \ et \ $p_1 : = 0$. Alors il existe une   application holomorphe \ $ \beta : X \to \C^*$ \ telles  que l'on ait, pour chaque \ $x \in X$ \ l'\'egalit\'e
$$ \A.\varphi(x) = \A.\psi(x) $$
o\`u 
\begin{equation*}
\psi(x) : = s^{\lambda-2}.Log\,s + \beta(x).s^{\lambda-1}. \tag{@}
\end{equation*}
\end{lemma}

\parag{Preuve} C'est une variante simple de la preuve de la proposition pr\'ec\'edente qui est laiss\'ee au lecteur. $\hfill \blacksquare$

\parag{Remarque} On verra que la proposition pr\'ec\'edente signifie  que la famille holomorphe \ $(E_{\lambda,p}(\alpha))_{\alpha \in \C^*}$ \ de th\`emes de rang 2 est universelle pour \ $\lambda > 1$ \ et \ $p \geq 1$ \  au sens de la d\'efinition \ref{verselle}. \\
De m\^eme,  lemme pr\'ec\'edent signifiera  que la famille constante (param\'etr\'ee par un point !) est universelle au sens de la d\'efinition \ref{verselle}. $\hfill \square$

\bigskip

\begin{defn}\label{Inv. fond.} Pour un th\`eme \ $[\lambda]-$primitif \ $E$ \ de rang 2 et d'invariant fondamentaux \ $\lambda_1, p_1$ \ avec \ $p_1 \geq 1$ \ nous appellerons {\bf invariant holomorphe} le nombre complexe (non nul) \ $\alpha$ \ tel que \ $E$ \ soit isomorphe \`a \ $E_{\lambda_1+p, \lambda_1}(\alpha)$.
\end{defn}

 Donc \ $\alpha$ \ est le nombre complexe donn\'e par la proposition pr\'ec\'edente appliqu\'ee \`a la famille constante \'egale \`a \ $E$.\\
 Il sera  commode de convenir que pour \ $p = 0$ \ l'invariant holomorphe de \ $E_{\lambda,\lambda}$ \  est \'egal \`a  \ $1$.
 
 \parag{Remarque} Comme le dual de \ $E_{\lambda,\lambda}$ \ est \ $E_{-\lambda+1,-\lambda+1}$, on constate que si le rationnel \ $\delta$ \ v\'erifie \ $\delta > \lambda$ \ alors \ $(E_{\lambda,\lambda})^*\otimes_{a,b} E_{\delta}$ \ est un th\`eme primitif de rang 2 d'invariants fondamentaux \ $\mu, 0$ \ (donc  isomorphe \`a \ $E_{\mu,\mu}$) \ avec \ $\mu = -\lambda+1 +\delta $.\\
 
 De m\^eme le dual de \ $E_{\lambda+p,\lambda}(\alpha)$ \ est \ $E_{-\lambda-p+1,-\lambda}((-1)^p.\alpha)$, et pour \ $\delta$ \ rationnel v\'erifiant \ $\delta > \lambda+p$, le (a,b)-module \ $(E_{\lambda+p,\lambda}(\alpha))^*\otimes_{a,b} E_{\delta}$ \ sera un th\`eme primitif de rang 2 d'invariants fondamentaux \ $-\lambda-p+\delta+1, p$ \ et d'invariant holomorphe \ $(-1)^p.\alpha$.\\
 
 Donc quitte \`a tensoriser\footnote{pour rendre les (a,b)-modules g\'eom\'etriques.} par \ $E_{\delta}$ \ avec \ $\delta$ \ rationnel assez grand (en fait plus grand que \ $\lambda_1$), la famille duale d'une famille holomorphe de th\`emes \ $[\lambda]-$primitifs de rang 2 est holomorphe.

\subsubsection{Crit\`ere d'holomorphie.}

La proposition ci-dessous montre que dans une famille holomorphe de th\`eme \ $[\lambda]-$primitif, la suite de Jordan-H{\"o}lder est "holomorphe".\\

\begin{prop}\label{J-H parametres}
Soit \ $X$ \ un espace complexe r\'eduit connexe et soit \ $\mathbb{E}$ \ une famille holomorphe de th\`emes \ $[\lambda]-$primitifs de rang \ $k$ \ param\`etr\'ee par  \ $X$.\\
 Notons \ $\lambda, p_1, \dots, p_{k-1}$ \ les invariants fondamentaux communs \`a chaque th\`eme de cette famille. Pour chaque \ $j \in [0,k]$ \ il existe une famille holomorphe unique  \ $\mathbb{F}_j$ \ de th\`emes \ $[\lambda]-$primitifs de rang \ $j$ \ param\`etr\'ee par \ $X$ \ et v\'erifiant les propri\'et\'es suivantes :
\begin{enumerate}[i)]
\item \ $\mathbb{F}_j \subset \mathbb{F}_{j+1} $ \ et \ $\mathbb{F}_0 = 0$ \ et \ $\mathbb{F}_k = \mathbb{E}$ ;
\item  pour chaque \ $x \in X$ le th\`eme \ $F_j(x)$ \ est {\bf le} sous-th\`eme normal de rang \ $j$ \ de \ $E(x)$.
\item La famille  \ $\mathbb{E}/\mathbb{F}_j$ \ des th\`emes  quotients \ $E(x)\big/F_j(x)$ \ est holomorphe.
\end{enumerate}
\end{prop}

\parag{Preuve} Le probl\`eme est local sur \ $X$, et l'on peut supposer que l'on a une application holomorphe \ $k-$th\'ematique
$$ \varphi : X \to \Xi_{\lambda}^{(k-1)} $$
telle que \ $\mathbb{E} = \mathbb{E}_{\varphi}$. Posons
$$ \varphi(x) : = \sum_{j=0}^{k-1} \quad \Sigma_j(x).s^{\lambda -1}.\frac{(Log\, s)^j}{j!} $$
o\`u \ $\Sigma_j \in \Gamma(X, \mathcal{O}_X[[b]])$. On a \ $\Sigma_{k-1} = b^{p_1+\dots+p_{k-1} -(k-1)}.\tilde{\Sigma}_{k-1}$ \ avec \ $\tilde{\Sigma}_{k-1}$ \ inversible dans \ $\mathcal{O}(X)[[b]]$. Donc, quitte \`a remplacer \ $\varphi$ \ par \ $\tilde{\Sigma}_{k-1}^{-1}.\varphi$ \ ce qui ne change pas \ $\mathbb{E}_{\varphi}$, on peut supposer que \ $\tilde{\Sigma}_{k-1} \equiv 1$, c'est \`a dire que l'on a
$$ \varphi -  s^{\lambda_k - 1}.\frac{(Log\, s)^{k-1}}{(k-1)!} \in \Gamma(X,\Xi_{X,\lambda}^{(k-2)}).$$
Posons alors \ $\psi : = (a - \lambda_k.b).\varphi $. Alors \ $\psi$ \ est \ $(k-1)-$th\'ematique, puisque \ $\psi, a.\psi, \dots, a^{k-2}.\psi$ \ est \ $\mathcal{O}(X)[[b]]-$libre et engendre \ $\mathbb{E}_{\psi}$ :\\
si on a \ $\sum_{j=0}^{k-2} \ U_j.a^j.\psi = 0 $ \ avec \ $U_j \in \mathcal{O}(X)[[b]]$, on aura
$$ \sum_{j=0}^{k-2} \ U_j.a^{j+1}.\varphi  - \sum_{j=0}^{k-2} \  U_j.a^j.\lambda_k.b.\varphi = 0 $$
ce qui impose successivement, puisque \ $a^j.b.\varphi \in \sum_{h=0}^j \mathcal{O}(X)[[b]].a^h.\varphi$,   les relations
\begin{align*}
 & U_{k-2} = 0, U_{k-3} = 0 \dots, U_1 = 0
  \end{align*}
 car \ $\varphi$ \ est th\'ematique. \\
 On obtient ainsi la famille holomorphe \ $\mathbb{F}_{k-1} : = \mathbb{E}_{\psi}$, et on conclut par une r\'ecurrence imm\'ediate. \\
 La propri\'et\'e \ $ iii)$ \  se d\'eduit facilement par r\'ecurrence du cas \ $j =1$. Dans ce cas il suffit de montrer que la compos\'ee \ $\theta : = f_{\lambda}\circ \varphi : X \to \Xi_{\lambda}^{(k-2)}$ \ est \ $(k-1)-$th\'ematique, puisque \ $F_1(x) = Ker\, f_{\lambda} \cap \A.\varphi(x)$ \ d'apr\`es le lemme \ref{devissage 0}. Pour cela montrons que \ $\sum_{j=0}^{k-2} \ S_j(x).a^j.\varphi(x)\in \Xi_{\lambda}^{(0)} $ \ implique \ $S_j(x) = 0, \forall j\in[0,k-2]$. En effet, sinon on aurait un entier \ $q \geq 0$ \ et un inversible \ $T $ \ de \ $ \C[[b]]$ \ qui v\'erifieraient
 $$  \sum_{j=0}^{k-2} \ S_j(x).a^j.\varphi(x) = T.s^{\lambda+q-1}.$$
 Alors \ $(a - (\lambda+q).b).T^{-1}.\big( \sum_{j=0}^{k-2} \ S_j(x).a^j\big)$ \ qui est un polyn\^ome en \ $a$ \ de degr\'e inf\'erieur ou \'egal \`a \ $k-1$ \ annulerait \ $\varphi(x)$ \ contredisant le fait que \ $E(x)$ \ est de rang \ $k$.  $\hfill \blacksquare$

 \bigskip

\begin{thm}[Crit\`ere d'holomorphie]\label{reciproque}
Soit \ $E(\sigma)_{\sigma \in X}$ \ une famille de th\`emes \ $[\lambda]-$primitifs d'invariants fondamentaux \ $\lambda_1, p_1, \dots, p_{k-1}$, o\`u l'on suppose \ $k \geq 2$. Soit \ $s_{k-1} : X \to \C$ \ l'application d\'efinie en associant \`a \ $\sigma \in X$ \ l'invariant holomorphe du  th\`eme \ $[\lambda]-$primitif de rang \ $2$ \ $E(\sigma)\big/F_{k-2}(\sigma)$.  Alors la famille \ $E(\sigma)_{\sigma \in X}$ \  est holomorphe si et seulement si 
\begin{enumerate}[i)]
\item \ $s_{k-1}$ \ est holomorphe sur \ $X$ ;
\item  la famille \ $F_{k-1}(\sigma))_{\sigma \in X}$ \ est holomorphe.
\end{enumerate}
\end{thm}

\parag{Remarque} Il est \'equivalent de demander \`a la fonction \ $s_{k-1}$ \ d\^etre holomorphe ou de demander que la famille de th\`eme \ $[\lambda]-$primitifs de rang 2 \ $(E(\sigma)\big/F_{k-2}(\sigma))_{\sigma \in X}$ \ soit holomorphe. \\
On notera que l'on aura pour chaque \ $\sigma \in X$ \ un isomorphisme 
 $$E(\sigma)\big/F_{k-2}(\sigma) \simeq E_{\lambda_k, \lambda_{k-1}}(s_{k-1}(\sigma)) \simeq \A\big/\A.(a - \lambda_{k-1}.b).(1 + s_{k-1}(\sigma).b^{p_{k-1}})^{-1}.(a - \lambda_k.b) .$$
 Donc le th\'eor\`eme pr\'ec\'edent est un crit\`ere n\'ecessaire et suffisant d'holomorphie, qui permet, par r\'ecurrence sur le rang, de se ramener au cas du rang \ $2$. $\hfill \square$

\bigskip

La d\'emonstration de ce th\'eor\`eme utilisera le lemme suivant :

\begin{lemma}\label{equ. diff.}
Soit \ $j, q \in \mathbb{N}$ \ et \ $\lambda \in ]0,1] \cap \mathbb{Q}$. Notons \ $H(j,q))$ \ l'hyperplan de \ $\Xi^{(j)}_{\lambda}$ \ correspondant \`a l'annulation du coefficient de \ $b^q.e_0$.\\
Alors l'application \ $ (a - (\lambda+q).b) : H(j,q) \oplus  \C.b^q.e_{j+1} \to b.\Xi^{(j)}_{\lambda} $ \ est un isomorphisme \ $\C-$lin\'eaire d'espaces de Frechet. En cons\'equence l'inverse est une application \ $\C-$lin\'eaire continue.
\end{lemma}

\parag{Preuve} On v\'erifie imm\'ediatement l'\'egalit\'e suivante pour tout couple d'entiers\ $(h, m) \in \mathbb{N}^2$ \ 
$$ (a - (\lambda+q).b).b^m.e_h = (m-q).b^{m+1}.e_h + b^{m+1}.e_{h-1} $$
toujours avec la convention \ $e_{-1} = 0$. On en d\'eduit que l'image de \ $\Xi^{(j)}_{\lambda}$ \ par \\ 
$(a - (\lambda+q).b)$ \ est l'hyperplan de \ $b.\Xi^{(j)}_{\lambda}$ \ donn\'e par l'annulation du coefficient de \ $b^{q+1}.e_j$ \ et que son noyau est \ $\C.b^q.e_0$. On conclut ais\'ement. $\hfill \blacksquare$

\parag{Remarque} Si l'on part d'un \'el\'ement de \ $b.\Xi^{(j)}_{\lambda}$ \ pour lequel le coefficient de \ $b^{q+1}.e_j$ \ vaut \ $\rho$, alors le coefficient de \ $b^q.e_{j+1}$ \ dans son image par l'application inverse sera \'egal \`a \ $\rho$. En particulier, il sera non nul quand \ $\rho \not= 0$. $\hfill \square$

\parag{Preuve du th\'eor\`eme \ref{reciproque}} Le probl\`eme est local et on peut donc supposer que l'on a une application holomorphe \ $k-$th\'ematique \ $\psi : X \to \Xi^{(k-2)}_{\lambda}$ \ telle que \ $\mathbb{E}_{\psi}$ \ donne l'holomorphie de la famille \ $F_{k-1}(\sigma))_{\sigma \in X}$. Il n'est pas restrictif, quitte \`a multiplier \ $\psi$ \ par un inversible de \ $\mathcal{O}(X)[[b]]$,  de supposer que l'on a
$$ \psi(\sigma) = b^{\lambda_{k-1}-\lambda}.e_{k-2} \quad modulo \quad \Xi^{(k-3)}_{\lambda}.$$
Posons \ $q : = \lambda_k -\lambda, S_{k-1} : = 1 + s_{k-1}.b^{p_{k-1}}$, et d\'efinissons
$$ \varphi : X \to \Xi^{(k-1)}_{\lambda} $$
en composant l'application holomorphe \ $S_{k-1}.\psi$ \ avec l'inverse de l'application  inverse construite dans le lemme pour \ $j : = k-2$ \  et l'inclusion   \'evidente de l'espace de Frechet
 $ H(k-2, \lambda_k-\lambda) \oplus \C.b^{\lambda_k-\lambda}.e_{k-1}$ \ dans \ $\Xi^{(k-1)}_{\lambda}$, en remarquant que l'on a \ $\psi(\sigma) \in b.\Xi^{(k-2)}_{\lambda}$ \ pour chaque \ $\sigma \in X$, puisque l'on a \ $F_{k-1}(\sigma) \subset a.E(\sigma) + b.E(\sigma)$, que \ $a.\Xi^{(j)}_{\lambda} \subset b.\Xi^{(j)}_{\lambda}$ \ pour tout \ $j \geq 0$, et que toute application\ $\A-$lin\'eaire de \ $F_{k-1}(\sigma)$ \  dans \ $\Xi$ \ est restriction d'une application \ $\A-$lin\'eaire de \ $E(\sigma)$ \ dans \ $\Xi$\footnote{Ceci r\'esulte de l'exactitude de
$$ 0 \to Hom_{\A}(E\big/F_{k-1}, \Xi) \to Hom_{\A}(E, \Xi) \to Hom_{\A}(F_{k-1}, \Xi) \to 0 $$
qui est un ph\'enom\`ene g\'en\'eral pour les suites exactes de (a,b)-modules g\'eom\'etriques (voir [B.05]).}.\\
On remarquera enfin que le coefficient de \ $b^{\lambda_k -\lambda}.e_{k-2}$ \ dans \ $S_{k-1}.\psi$ \ co{\"i}ncide avec celui de \ $b^{p_{k-1}}$ \ dans \ $S_{k-1}$, c'est-\`a-dire est \'egal \`a \ $s_{k-1}$. Il est donc non nul pour chaque \ $\sigma \in X$, ce qui montre que le coefficient de \ $b^{\lambda_k-\lambda}.e_{k-1}$ \ dans \ $\varphi(\sigma)$ \ est non nul pour chaque \ $\sigma \in X$, gr\^ace \`a la relation \ $\lambda_k = \lambda_{k-1} + p_{k-1} - 1$ \ qui implique \ $\lambda_k - \lambda + 1 = \lambda_{k-1} - \lambda + p_{k-1}$. Ceci est \'evidemment n\'ecessaire pour que \ $E(\sigma)$ \ soit un th\`eme de rang \ $k$. $\hfill \blacksquare$

\bigskip

\subsubsection{Le th\'eor\`eme de dualit\'e.}

\bigskip

\begin{prop}[D\'ecalage]\label{decal.}
Soit \ $(E(\sigma))_{\sigma \in X}$ \ une famille holomorphe de th\`emes \ $[\lambda]-$primitifs. Soit \ $\delta  \in \mathbb{Q}$ \ tel que \ $E(\sigma)\otimes_{a,b} E_{\delta}$ \ soit un th\`eme pour chaque \ $\sigma \in X$. Alors la famille \ $E(\sigma)\otimes_{a,b} E_{\delta}$ \  est une famille holomorphe. En particulier pour \ $r \in \mathbb{Z}$ \ tel que pour tout \ $\sigma \in X$ \ $b^r.E(\sigma)$ \ soit un th\`eme, alors  \ $(b^r.E(\sigma))_{\sigma \in X}$ \ est une famille holomorphe.
\end{prop}

\parag{Preuve} On rappelle que le (a,b)-module \ $E \otimes_{a,b} E_{\delta}$ \ est le (a,b)-module obtenu en rempla{\c c}ant l'action de \ $a$ \ par \ $a + \delta.b$. Donc la condition pour que les  \ $E(\sigma)$ \ soient g\'eom\'etriques\footnote{qui est la seule condition qui peut ne pas \^etre r\'ealis\'ee pour avoir un th\`eme.} est que \ $\lambda_1 + \delta > k-1 $. C'est-\`a-dire que \ $ \delta > -  \lambda_1 + k -1  $. Comme \ $\lambda_1 - k + 1 > 0 $ \ ceci a lieu en particulier pour tout \ $\delta \in \mathbb{Q}^+$.\\
La d\'emonstration de la proposition est triviale. $\hfill \blacksquare$

\bigskip

\begin{thm}[Th\'eor\`eme de dualit\'e]\label{th. dual}
Soit \ $(E(\sigma))_{\sigma \in X}$ \ une famille holomorphe de th\`emes \ $[\lambda]-$primitifs. Soit \ $\delta \in \mathbb{Q}$ \ un rationnel assez grand pour que \ chaque \ $E(\sigma)^*\otimes_{a,b} E_{\delta}$ \ soit un th\`eme. Alors la famille \ $((E(\sigma)^*\otimes E_{\delta})_{\sigma \in X}$ \ est holomorphe.
\end{thm}

 \parag{D\'emonstration du th\'eor\`eme \ref{th. dual}} Nous allons faire une r\'ecurrence sur le rang des th\`emes de la famille holomorphe consid\'er\'ee. Comme en rang 1 et 2 le th\'eor\`eme est d\'ej\`a d\'emontr\'e (voir la remarque qui suit la d\'efinition \ref{Inv. fond.} ), supposons le th\'eor\`eme d\'emontr\'e en rang \ $k-1 \geq 2$ \ et montrons-le en rang \ $k$.\\
Soit \ $(F_1(\sigma))_{\sigma \in X}$ \ la famille des sous-th\`emes normaux de rang \ $1$ \ des th\`emes \ $(E(\sigma))_{\sigma \in X}$. Il r\'esulte du iv) de la proposition  \ref{J-H parametres} que la famille \ $(E(\sigma)\big/F_1(\sigma))_{\sigma \in X}$ \  est une famille holomorphe de th\`emes.\\
 De m\^eme, il r\'esulte  de la proposition \ref{J-H parametres}  que la famille \ $(F_2(\sigma))_{\sigma \in X}$ \ est holomorphe. L'hypoth\`ese de r\'ecurrence donne alors que, pour \ $\delta \in \mathbb{Q}$ \ assez grand, les familles de th\`emes  \ $\big((E(\sigma)\big/F_1(\sigma))^*\otimes_{a,b} E_{\delta}\big)_{\sigma \in X}$ \ et \ $\big((F_2(\sigma))^*\otimes E_{\delta}\big)_{\sigma \in X}$ \ sont holomorphes. Mais alors le crit\`ere d'holomorphie \ref{reciproque} s'applique \`a la famille de th\`emes \ $\big((E(\sigma)^*\otimes_{a,b} E_{\delta}\big)_{\sigma \in X}$ \ puisque la famille des sous-th\`emes de rang 
 \ $k-1$ \ associ\'ee est pr\'ecis\'ement la famille \ $\big((E(\sigma)\big/F_1(\sigma))^*\otimes_{a,b} E_{\delta}\big)_{\sigma \in X}$ \ et que la famille des quotients de rang 2  associ\'ee est pr\'ecis\'ement la famille \ $\big((F_2(\sigma))^*\otimes_{a,b} E_{\delta}\big)_{\sigma \in X}$. $\hfill \blacksquare$

\subsection{Familles standards de th\`emes \ $[\lambda]-$primitifs.}

Nous fixerons dans ce paragraphe \ $\lambda_1 \in k-1 + \mathbb{Q}^{*+}$ \ et les entiers \ $p_1, \dots , p_{k-1}$. Pour \ $j \in [1,k-1]$ \ nous d\'efinirons l'ouvert affine  \ $W_j \subset V_j$ \ de l'espace vectoriel \ $V_j$ \ d\'efini dans la proposition \ref{Supplementaires}  de la fa{\c c}on suivante :\\

si \ $p_j + \cdots + p_{k-1} < k-j $ 
\begin{equation*}
W_j : = \{ S_j \in \C[b] \ / \  S_j(0) = 1, \ \deg(S_j) \leq k-j-1 \quad {\rm et} \quad  {\rm coeff} \ b^{p_j} \not= 0 \}  \tag{@}
\end{equation*}
si \ $p_j + \cdots + p_{k-1} \geq k-j$ \ d\'efinissons l'entier \ $q_j \geq k-j$ \ comme le plus petit entier de la forme \ $p_j + \cdots + p_{j+h}$ \ qui v\'erifie \ $p_j + \cdots + p_{j+h} \geq k-j$, et posons 
\begin{equation*}
W_j : = \{ S_j \in \C[b] \ / \  S_j(0) = 1, S_j \in \sum_{h=0}^{k-j-1} \C.b^h + \C.b^{q_j} \quad {\rm et} \quad {\rm coeff} \ b^{p_j} \not= 0 \}  \tag{@@}
\end{equation*}
Posons alors  
$$ \mathcal{S}(\lambda_1, p_1, \dots, p_{k-1}) : = \{ (S_1, \dots, S_{k-1} \big/ \  S_j \in W_j \quad \forall j \in [1,k-1] \} .$$
Pour \ $\sigma \in \mathcal{S}(\lambda_1, p_1, \dots, p_{k-1})$ \ notons \ $E(\sigma)$ \ le th\`eme \ $[\lambda]-$primitif d'invariants fondamentaux \ $(\lambda_1, p_1, \dots, p_{k-1})$ \ d\'efini par 
\begin{align*}
&  E(\sigma) : = \A \big/\A.P(\sigma) \quad {\rm avec} \tag{$\sigma$} \\
&  P(\sigma) : = (a - \lambda_1.b).S_1^{-1}\cdots S_{k-1}^{-1}.(a - \lambda_k.b) \tag{$\sigma'$}
\end{align*}
o\`u nous avons pos\'e \ $\sigma : = (S_1, \dots, S_{k-1})$.

\bigskip

\begin{defn}\label{famille standard}
Nous appellerons {\bf famille standard d'invariants fondamentaux} \ $(\lambda_1, p_1, \dots, p_{k-1})$ \ la famille \ $E(\sigma)_{\sigma \in  \mathcal{S}(\lambda_1, p_1, \dots, p_{k-1})}$.
\end{defn}

\bigskip

\parag{Exemple}
Pour \ $k = 1$ \ chaque \ $\lambda_1 \in \mathbb{Q}^{+*}$ \ la famille standard associ\'ee est r\'eduite au th\`eme \ $E_{\lambda_1}$.\\ 
Pour \ $k = 2$, et les invariants fondamentaux \ $(\lambda_1, p_1)$ \ on a
\begin{enumerate}[i)]
\item Pour \ $p_1 = 0$, \ $\mathcal{S}(\lambda_1, p_1) = \{1\}$ \ et le th\`eme correspondant est \\
 $E : = \A\big/\A.(a-\lambda_1.b).(a - (\lambda_1-1).b$.
\item Pour \ $p_1 \geq 1$ \ on a \ $\mathcal{S}(\lambda_1, p_1) = \{ 1 + \alpha.b^{p_1}, \alpha \in \C^* \}$ \ et le th\`eme associ\'e \`a \ $\alpha$ \ est \ $E = \A\big/\A.(a - \lambda_1.b).(1 + \alpha.b^{p_1})^{-1}.(a - (\lambda_1+p_1 -1).b) $. $\hfill \square$
\end{enumerate}

\bigskip

\begin{thm}\label{anal. standard}
Quelques soient les invariants fondamentaux fix\'es, la famille de th\`emes param\'etr\'ee par \ $\mathcal{S}(\lambda_1, p_1, \dots, p_{k-1})$ \ est holomorphe.
\end{thm}

\parag{D\'emonstration} Ce th\'eor\`eme s'obtient imm\'ediatement \`a partir de la proposition \ref{reciproque} par une r\'ecurrence sur le rang gr\^ace au th\'eor\`eme de dualit\'e et \`a la proposition de d\'ecalage. En effet, la famille \ $\mathbb{E}\big/\mathbb{F}_1$ \ correspond \`a la famille de th\`emes param\'etr\'ee par \ $\mathcal{S}(\lambda_2, p_2, \dots, p_{k-1})$ \ qui est holomorphe par hypoth\`ese de r\'ecurrence, et la famille \ $\mathbb{F}_2$ \ correspond soit \`a la famille constante \ $\mathcal{S}(\lambda_1,p_1 = 0)$ \ soit au cas trait\'e au lemme \ref{universelle rg 2}. On conclut en appliquant la proposition \ref{reciproque}  \`a la famille duale suffisamment d\'ecal\'ee et en utilisant \`a nouveau le th\'eor\`eme de dualit\'e \ref{th. dual} et la proposition de d\'ecalage \ref{decal.}. $\hfill \blacksquare$

\bigskip

\subsection{Les d\'eformations standards sont verselles.}

Commen{\c c}ons par deux  d\'efinitions.

\begin{defn}\label{pull-back}
Soit \ $X$ \ un espace complexe r\'eduit et soit \ $\mathbb{E}$ \ un faisceau sur \ $X$ \  de \ $\mathcal{O}_X-(a,b)-$modules. Soit \ $f : Y \to X$ \ une application holomorphe d'un espace complexe r\'eduits \ $Y$ \ dans \ $X$. On appellera {\bf image r\'eciproque de \ $\mathbb{E}$ \ par \ $f$}, not\'e \ $f^*\mathbb{E}$, le faisceau de \ $\mathcal{O}_Y-(a,b)-$modules d\'efini comme suit :\\

Si \ $\mathbb{E}$ \ est \ $\mathcal{O}_X[[b]]-$libre de rang \ $p$ \ sur l'ouvert \ $U$ \ et de base \ $e_1, \dots, e_p$, alors \ $f^*\mathbb{E}$ \ est \ $ \mathcal{O}_Y[[b]]-$libre de rang \ $p$ sur l'ouvert \ $f^{-1}(U)$ \ et de base \ $f^*e_1, \dots ,f^*e_p$. L'application \ $a$ \ sur un tel ouvert est d\'efinie par la formule
$$ a.f^*e_j = \sum_{i=1}^p \  f^*S_{i,j}.f^*e_i $$
si l'on a sur \ l'ouvert \ $U$ \ la formule \ $a.e_j = \sum_{i=1}^p S_{i,j}.e_i $. Ici les \ $S_{i,j}$ \ sont dans \ $\mathcal{O}_X(U)[[b]]$ \ et \ $f^*S$ \ pour \ $S \in \mathcal{O}_X(U)[[b]]$ \ d\'esigne l'\'el\'ement de \ $\mathcal{O}_Y(f^{-1}(U)[[b]]$ \ d\'eduit de \ $S : = \sum_{\nu = 0}^{\infty}\  s_{\nu}.b^{\nu}$ \ via la formule
\ $ f^*S = \sum_{\nu = 0}^{\infty} \ f^*s_{\nu}.b^{\nu} $.
\end{defn}

\bigskip

\begin{defn}\label{verselle}
Soit \ $X$ \ un espace complexe r\'eduit et soit \ $\mathbb{E}$ \ une famille holomorphe de th\`emes \ $[\lambda]-$primitifs param\'etr\'ee par \ $X$. Soit \ $x_0 \in X$. On dira que la famille \ $E$ \ est {\bf verselle} au voisinage de \ $x_0$ \ si la condition suivante est r\'ealis\'ee :\\
\begin{itemize}
\item Pour toute famille holomorphe \ $\mathbb{G}$ \ de \ th\`emes \ $[\lambda]-$primitifs param\'etr\'ee par un espace complexe r\'eduit \ $Y$ \ telle que le th\`eme \ $G(y_0)$ \ soit isomorphe \`a \ $E(x_0)$, il existe un voisinage ouvert \ $U$ \ de \ $y_0$ \ dans \ $Y$, un voisinage ouvert \ $V$ \ de \ $x_0$ \ dans \ $X$, une application holomorphe \ $f : U \to V$ \ telle que les faisceaux de \ $\mathcal{O}_Y-(a,b)-$modules \ $f^*\mathbb{E}_{\vert U}$ \ et \ $\mathbb{G}_{\vert U}$ \ soit isomorphes.
\end{itemize}
Quand l'application \ $f$ \ est unique sur un voisinage ouvert assez petit de \ $x_0$, on dira que la famille est {\bf universelle} au voisinage de \ $x_0$.\\
Une famille verselle (resp. universelle)  au voisinage de chaque point de \ $X$ \ sera dite {\bf verselle} (resp. {\bf universelle} ).
\end{defn}

Voici le th\'eor\`eme principal de ce paragraphe.

\begin{thm}\label{standard est verselle}
Pour tout choix d'invariants fondamentaux  \ $\lambda_1, p_1, \dots, p_{k-1}$ \ la famille standard de th\`emes \ $[\lambda]-$primitifs param\'etr\'ee par \ $S(\lambda_1, p_1, \dots, p_{k-1})$  \ est verselle.
\end{thm}

\parag{D\'emonstration} Nous allons montrer ce r\'esultat par r\'ecurrence sur \ $k$. Les cas \ $k = 1$ \ et \ $k = 2$ \ ont d\'ej\`a \'et\'e trait\'es (voir 4.1.2). Supposons donc \ $k \geq 3$ \ et le cas \ $k-1$ \ \'etabli.Pr\'ecisons que l'assertion \'etant locale, il nous suffit de prouver l'assertion au voisinage d'un point donn\'e de \ $X$.\\
Notons \ $\mathbb{F}_1 \subset \mathbb{E}$ \ le sous-faisceau de \ $\mathcal{O}_X-(a,b)-$modules donnant la famille holomorphe des sous-th\`emes normaux de rang 1  de la famille \ $E$. Alors le faisceau \ $\mathbb{E}\big/\mathbb{F}_1$ \ est une famille holomorphe de th\`eme \ $[\lambda]-$primitifs de rang \ $k-1$  \ et d'invariants fondamentaux \ $\lambda_2, p_2, \dots, p_{k-1}$, o\`u \ $\lambda_2 = \lambda_1 + p_1 -1$.\\
L'hypoth\`ese de r\'ecurrence nous fournit alors, localement sur \ $X$ \ une application holomorphe \ $f : X \to S(\lambda_2, p_2, \dots, p_{k-1})$ \ telle que l'image r\'eciproque par \ $f$ \ de la famille standard associ\'ee soit isomorphe \`a la famille \ $\mathbb{E}\big/\mathbb{F}_1$.\\
Comme tout ceci est local au voisinage d'un point \ $x_0$ \ de \ $X$ \ que l'on suppose fix\'e, on peut supposer que la famille holomorphe \ $\mathbb{E}$ \ est donn\'ee par une application holomorphe \ $k-$th\'ematique \ $\varphi : X \to \Xi^{(k-1)}_{\lambda} $ \ v\'erifiant
$$\varphi(x) = s^{\lambda_k-1}.\frac{(Log\,s)^{k-1}}{(k-1)!} + \psi(x) $$
o\`u \ $\psi$ \ est holomorphe \`a valeurs dans \ $ \Xi^{(k-2)}_{\lambda} $. 

L'application holomorphe \ $f$ \ nous fournit en fait des applications holomorphes \ $S_2, \dots, S_{k-1} : X \to \C[[b]] $ \ v\'erifiant \ $S_j(0) \equiv 1$ \ et telles que, si l'on pose
 $$P_1 : = (a - \lambda_2.b).S_2^{-1}\dots S_{k-1}^{-1}.(a - \lambda_k.b)$$
  on ait \ $P_1.e = 0$ \ pour le g\'en\'erateur standard \ $e$ \ de la famille standard param\'etr\'ee par \ $S(\lambda_2, p_2, \dots, p_{k-1})$. Donc le g\'en\'erateur \ $f^*e$ \ de \ $\mathbb{E}\big/\mathbb{F}_1$ \ v\'erifie \'egalement \ $P_1.f^*e =0$ \ dans \ $\mathbb{E}\big/\mathbb{F}_1$. En identifiant \ $\mathbb{E}$ \ et son image par \ $\varphi$ \ dans \ $\mathcal{O}_X \ptc \Xi^{(k-1)}_{\lambda}$ \ et \ $\mathbb{F}_1$ \ avec \\ 
  $\mathbb{E} \cap \mathcal{O}_X \ptc \Xi^{(0)}_{\lambda}$, on identifie alors \ $\mathbb{E}\big/\mathbb{F}_1$ \ \`a un sous faisceau du faisceau quotient  
   $$\mathcal{O}_X \ptc \Xi^{(k-1)}_{\lambda} \big/ \mathcal{O}_X \ptc \Xi^{(0)}_{\lambda} \simeq \mathcal{O}_X \ptc \Xi^{(k-2)}_{\lambda}.$$
    On peut donc trouver \ $T_0, \dots, T_{k-1}$ \ des sections locales de \ $\mathcal{O}_X[[b]]$ \ telles que l'image de la section
 $$\sigma : =  \sum_{j=0}^{k-1} \ T_j.a^j.\varphi $$
 dans \ $\mathbb{E}\big/\mathbb{F}_1 $ \ soit \ $f^*e$. De plus comme \ $f^*e$ \ engendre \ $\mathbb{E}\big/\mathbb{F}_1$ \ la section \ $T_0$ \ devra \^etre un inversible de \ $\mathcal{O}_X[[b]]$ \ au voisinage de \ $x_0$, sinon la valeur en \ $x_0$ \ est dans \ $a.E(x_0) + b.E(x_0)$ \ et son image ne peut engendrer \ $E(x_0)\big/F_1(x_0)$. Alors la section \ $\sigma$ \ engendre localement \ $\mathbb{E}$ \ et v\'erifiera
 $$ P_1.\sigma \in \mathbb{F}_1 .$$
 Mais on sait que \ $\mathbb{F}_1 = \mathcal{O}_X[[b]] \otimes s^{\lambda_1 -1}$, ce qui permet d'\'ecrire
 $$ P_1.\sigma  = \Theta_1.s^{\lambda_1 -1} $$
 o\`u \ $\Theta_1$ \ est une section locale de \ $\mathcal{O}_X[[b]] $. La d\'ecomposition
 $$ \mathcal{O}_X\ptc E_{\lambda_1} = P_1.(\mathcal{O}_X \ptc E_{\lambda_1}) \oplus (\mathcal{O}_X \otimes V_1) $$
 permet alors d'\'ecrire \ $\Theta_1.s^{\lambda_1-1} = P_1.\alpha + S_1.s^{\lambda_1-1} $ \ o\`u \ $\alpha$ \ est une section locale de \ $\mathbb{F}_1$ \ et \ $S_1$ \ une section locale de \ $\mathcal{O}_X\otimes V_1$. De plus l'inversibilit\'e de \ $\Theta_1$ \ assure l'inversibilit\'e de \ $S_1$ \ dans \ $\mathcal{O}_X[[b]]$, c'est \`a dire l'inversibilit\'e de son terme constant en \ $b$ \ dans \ $\mathcal{O}_X$. Donc quitte \`a multiplier \ $\sigma$ \ et \ $\alpha$ \ par un inversible  \ $I$ \ de \ $\mathcal{O}_X$, on pourra supposer que le terme constant en \ $b$ \ de \ $S_1$ \ est identiquement \'egal \`a \ $1$.\\
 Alors \ $\tau : = I.(\sigma - \alpha)$ \ est encore un g\'en\'erateur local  de \ $\mathbb{E}$ \ et il v\'erifie
 $$ P_1.\tau = S_1.s^{\lambda_1-1} \quad {\rm avec} \quad S_1 \in \mathcal{O}_X\otimes V_1, S_1(0) \equiv 1 $$
 ce qui donne \ $ (a - \lambda_1.b).S_1^{-1}.P_1.\tau = 0 $. On constate alors que \ $\mathbb{E}$ \ est isomorphe \`a l'image r\'eciproque de la famille standard par l'application \ $g$ \ donn\'ee au voisinage de \ $x_0$ \ par \ $S_1, \dots, S_{k-1}$, en envoyant le g\'en\'erateur local \ $\tau$ \ sur l'image r\'eciproque \ $g^*e$ \ du  g\'en\'erateur standard \ $e$ \  de la famille param\'etr\'ee par \ $S(\lambda_1, p_1, \dots, p_{k-1})$. En effet si \ $P_0 : = (a -\lambda_1.b).S_1^{-1}.P_1$ \ on aura \ $P_0.\tau = 0$ \ ainsi que \ $P_0.g^*e = 0$, puisque \ $P_0.e = 0$. $\hfill \blacksquare$
 
 \bigskip
 
 \begin{cor}\label{Univers. cas stable.}
 Soient \ $\lambda_1, p_1, \dots, p_{k-1}$ \ les invariants fondamentaux d'un th\`eme \ $[\lambda]-$primitif. Si tout th\`eme admettant ces invariants fondamentaux est stable, la famille standard associ\'ee \`a ces invariants fondamentaux est universelle.
\end{cor}

\parag{Preuve} Ceci r\'esulte imm\'ediatement du fait que sous notre hypoth\`ese, deux param\`etres distincts donnent deux th\`emes non isomorphes gr\^ace au th\'eor\`eme \ref{unicite}. La versalit\'e de la famille standard montr\'ee ci-dessus au th\'eor\`eme \ref{standard est verselle} permet alors imm\'ediatement de conclure. $\hfill \blacksquare$

\subsection{Un contre-exemple.}

Nous allons donner un exemple  de th\`emes de rang 3  pour lesquels il n'existe pas de famille universelle.

On fixe les invariants \ $\lambda_1, p_1 = p_2 = 1$ \ pour les th\`emes de rang 3 que nous allons consid\'erer maintenant. On a donc \ $q_1 = 2 $ \ et  \ $q_2 = p_2 = 1$.

Notre objectif est de montrer la proposition suivante.

\begin{prop}\label{contre-exemple}
Il n'existe pas de famille universelle pour les th\`emes \ $[\lambda]-$primitifs de rang 3 d'invariants fondamentaux  \ $\lambda_1, p_1 = p_2 = 1$, au voisinage de chacun des th\`emes de param\`etres \ $(\alpha,\alpha,\gamma)$, avec \ $\alpha \not= 0$, c'est-\`a-dire au voisinage de chacun des th\`emes stables (sp\'eciaux) de la famille verselle standard.
\end{prop}

\bigskip

La preuve de cette proposition utilisera les trois lemmes suivants.

\bigskip

\begin{lemma}\label{isom.}
On consid\`ere, pour \ $\alpha,\beta, \gamma \in \C, \alpha.\beta \not= 0 $ \ les th\`emes de rang 3 \ $E_{\alpha,\beta,\gamma}$ \ d\'efinis de la fa{\c c}on suivante :
\begin{align*}
&  (a - \lambda.b).e_3 = (1 + \alpha.b).e_2 \\
&  (a - \lambda.b).e_2 = (1 + \beta.b + \gamma.b^2).e_1 \\
& (a - \lambda.b).e_1 = 0 .
\end{align*}
Pour \ $\beta \not= \alpha$, \ $E_{\alpha,\beta,\gamma}$ \ est isomorphe \`a \ $E_{\alpha,\beta, 0}$ \ quelque soit \ $\gamma$.\\
Pour \ $\beta = \alpha$, les th\`emes \ $E_{\alpha,\alpha,\gamma}$ \ et \ $E_{\alpha,\alpha,\gamma'}$ \ sont isomorphes si et seulement si \ $\gamma = \gamma'$.
\end{lemma}

\parag{Preuve} Cherchons une \ $\C[[b]]-$base \ $\varepsilon_3, \varepsilon_2, \varepsilon_1$ \ de \ $E_{\alpha,\beta,\gamma}$ \ v\'erifiant les conditions suivantes :
\begin{align*}
& \varepsilon_3 = e_3 + U.e_2 + V.e_1, \quad {\rm avec} \quad U,V \in \C[[b]] \tag{0}\\
& (a - \lambda.b).\varepsilon_3 = (1 + \alpha.b).\varepsilon_2 \tag{1} \\
& (a - \lambda.b).\varepsilon_2 = (1 + \beta.b + \gamma'.b^2).\varepsilon_1 \tag{2}\\
& (a - \lambda.b).\varepsilon_1 = 0. \tag{3}
\end{align*}
On sait en effet que \ $\alpha$ \ et \ $\beta$ \ sont d\'etermin\'es par la classe d'isomorphisme du th\`eme \ $E(\alpha,\beta,\gamma)$ \ puisque l'on a \ $p_1 = p_2 = 1$; on notera que \ $q_1 = p_1 + p_2 = 2$. La derni\`ere \'egalit\'e \ $(3)$ \  impose \ $\varepsilon_1 = \rho.e_1$ \ avec \ $\rho \in \C^*$.\\
Calculons les conditions impos\'ees \`a \ $U$ \ et \ $V$ :
\begin{align*}
& (a - \lambda.b).\varepsilon_3 =   (1 + \alpha.b).e_2 + b^2.U'.e_2 + U. (1 + \beta.b + \gamma.b^2).e_1 + b^2.V'.e_1\\
& \qquad \qquad  = (1+\alpha.b).\varepsilon_2 \qquad \qquad  {\rm et \ donc} \\
& \varepsilon_2 = Z.e_2 + T.e_1 \qquad \qquad {\rm avec} \\
& Z = (1 +\alpha.b)^{-1}.(1 + \alpha.b + b^2.U') \quad {\rm et}\\
& (1 + \alpha.b).T = U.(1 + \beta.b + \gamma.b^2) + b^2.V' \tag{4}
\end{align*}
On aura alors 
\begin{align*}
& (a - \lambda.b).\varepsilon_2 = Z.(1 + \beta.b + \gamma.b^2).e_1 + b^2.Z'.e_2 + b^2.T'.e_1 \\
& \qquad \qquad = (1 + \beta.b + \gamma'.b^2).\rho.e_1
\end{align*}
ce qui implique d\'ej\`a \ $Z' = 0 $ \ et comme \ $Z = 1 + (1 + \beta.b)^{-1}.b^2.U'$ \ on doit avoir \ $U \in \mathbb{C}$, et \ $Z = 1$. La relation \ $(2)$ \ donne maintenant, puisque \ $\varepsilon_2 = e_2 + T.e_1$
$$  (1 + \beta.b + \gamma.b^2).e_1 + b^2.T'.e_1 =  (1 + \beta.b + \gamma'.b^2).\rho.e_1 $$
ce qui impose \ $\rho = 1$ \ et \ $ T' = \gamma' - \gamma $. On aura donc \ $T = U + (\gamma' - \gamma).b$ \ en identifiant les termes constants de \ $(4)$. Cette \'egalit\'e \ $(4)$ \ impose de plus
\begin{align*}
& \alpha.U + \gamma' - \gamma = U.\beta \qquad \qquad {\rm et}  \tag{5} \\
& U.\gamma + V' = \alpha.(\gamma - \gamma') 
\end{align*}
On en d\'eduit que pour \ $\alpha \not= \beta $ \ on aura
$$U = \frac{\gamma - \gamma'}{\alpha - \beta} \quad {\rm et} \quad  V = V_0 + \frac{\gamma'- \gamma}{\beta - \alpha}.\big(\alpha.(\beta-\alpha) - \gamma\big).b .$$

Si \ $\beta = \alpha$, la relation \ $(5)$ \ impose \ $\gamma = \gamma'$. $\hfill \blacksquare$

\bigskip

Pour \ $(\alpha,\beta) \in X: = \{ (\alpha,\beta) \in (\C^*)^2, \alpha \not= \beta\}$ \ notons \ $E(\alpha, \beta)$ \ le th\`eme de rang \ $3$ \ d\'efini par \ $E(\alpha,\beta) : = \A\big/\A.(a -\lambda.b)(1+\beta.b)^{-1}(a - \lambda.b)(1 + \alpha.b)^{-1}(a - \lambda.b) $.

\bigskip 

\begin{lemma}\label{non stables} Il n'existe pas d'endomorphisme de rang \  $2$ \  de \ $E(\alpha,\beta)$ \ pour \ $\alpha \not= \beta$.
\end{lemma}

\parag{Preuve} Il nous suffit de montrer qu'il n'existe pas d'\'el\'ement \ $x : = e_2 + U.e_1$ \ dans \ $E(\alpha,\beta)$ \ v\'erifiant \ $(a-\lambda.b)(1+\alpha.b)^{-1}(a -\lambda.b).x = 0 $, o\`u \ $U \in \C[[b]]$.
Comme les \'el\'ements de \ $E(\alpha,\beta)$ \ annul\'es par \ $(a - \lambda.b)$ sont de la forme \ $\rho.e_1$ \ avec \ $\rho \in \C$, un tel \ $x$ \ doit v\'erifier
$$ (a - \lambda.b)x = \rho.(1 + \alpha.b).e_1 $$
ce qui impose \`a \ $U$ \ de v\'erifier la relation
$$ (1 + \beta.b) + b^2.U' = \rho.(1 + \alpha.b) .$$
On en conclut que l'on doit avoir \ $\rho = 1$ \ et donc \ $\alpha = \beta$. $\hfill \blacksquare$

\bigskip

Par contre, pour \ $\alpha = \beta \not= 0$ \ et \ $\gamma$ \ arbitraire on a stabilit\'e .

\begin{lemma}\label{Stab. alpha = beta}
Pour \ $\alpha \not= 0$ \ le (a,b)-module \  $E_{\alpha,\alpha,\gamma}$ \ est un th\`eme stable de rang 3.
\end{lemma}

\parag{Preuve} Il nous suffit de montrer qu'il existe \ $x : = e_2 + U.e_1$ \ v\'erifiant  
$$(a-\lambda.b)(1+\alpha.b)^{-1}(a -\lambda.b).x = 0 , $$
 o\`u \ $U \in \C[[b]]$. Comme \ $F_2 $ \ est un th\`eme, les \'el\'ements de \ $F_2$ \ annul\'es par \ $(a - \lambda.b)$ \ sont de la forme \ $\rho.e_1, \rho \in \C$. Donc \ $x$ \ doit v\'erifier
$$ (a - \lambda.b)x = \rho.(1 + \alpha.b).e_1 $$
ce qui impose \`a \ $U$ de v\'erifier la relation
$$ (1 + \beta.b+ \gamma.b^2) + b^2.U' = \rho.(1 + \alpha.b) .$$
On en conclut que l'on doit avoir \ $\rho = 1$ \ et \ $U = -\gamma.b + cste$. On a donc une solution
\ $ x : = e_2 -\gamma.b.e_1$. $\hfill \blacksquare$

\parag{Preuve de la proposition \ref{contre-exemple}} Le fait que les th\`emes stables de la famille standard consid\'er\'ee sont exactement  les \ $E_{\alpha,\alpha,\gamma}$ \ est d\'emontr\'e dans les lemmes \ref{non stables} et \ref{Stab. alpha = beta}.\\
La famille \ $(E_{\alpha,\beta,\gamma})_{(\alpha,\beta,\gamma) \in S(\lambda_1,  p_1 = p_2 =1)}$ \ est une famille holomorphe et m\^eme verselle en chaque point d'apr\`es le th\'eor\`eme \ref{standard est verselle}. Supposons trouv\'ee une famille universelle \ $(E_y)_{y \in Y}$ \ au voisinage du th\`eme \ $E(\alpha_0,\alpha_0,\gamma_0) \simeq E_{y_0}$, o\`u \ $Y$ \ est un espace complexe r\'eduit que l'on peut supposer plong\'e dans \ $\C^n$ \ au voisinage de \ $y_0$. Consid\'erons alors l'application holomorphe \ $\varphi : \Omega \to Y \hookrightarrow \C^N$ \ classifiant la famille standard sur un voisinage ouvert \ $\Omega$ \ de \ $(\alpha_0,\alpha_0,\gamma_0) \in (\C^*)^2\times \C $. Comme pour \ $\alpha \not= \beta$ \ la classe d'isomorphisme de \ $E_{\alpha,\beta,\gamma}$ \ ne d\'epend pas de \ $\gamma$ \ d'apr\`es le lemme \ref{isom.}, on aura \ $\frac{\partial \varphi}{\partial \gamma} \equiv 0$ \ sur l'ouvert \ $\{\alpha \not= \beta\}$ \ de \ $\Omega$. Ceci impose \`a \ $\varphi$ \ d'\^etre ind\'ependante de \ $\gamma$ \ ce qui donnerait l'isomorphisme entre \ $E_{\alpha,\alpha,\gamma}$ \ et \ $E_{\alpha,\alpha, \gamma'}$ \ pour tout \ $\alpha$ \ assez voisin de \ $\alpha_0$ \ et tout \ $\gamma,\gamma'$ \ assez voisins de \ $\gamma_0$. Ceci contredit le lemme \ref{isom.}. $\hfill \blacksquare$

\begin{cor}\label{Univers.}
La famille \ $E(\alpha,\beta)_{(\alpha,\beta)\in X}$ \ est universelle en chaque point de  \ $X : = (\C)^*\setminus \{\alpha = \beta\} $.
\end{cor}

\parag{Preuve} Notons \ $\mathbb{E}$ \ le faisceau sur \ $X$ \ de \ $\mathcal{O}_X-(a,b)-$module associ\'e \`a la famille holomorphe des \ $E(\alpha,\beta)$. Il nous suffit en fait de montrer que l'application holomorphe
$$\pi :  X \times C \to X $$
d\'efinie par \ $\pi(\alpha,\beta, \gamma) = (\alpha,\beta)$ \ v\'erifie bien que \ $\pi^*(\mathbb{E})$ \ est un faisceau de \ $\mathcal{O}_{X\times\C}-(a,b)-$modules isomorphe au faisceau associ\'e \`a la famille standard param\'etr\'ee par \ $X \times \C$. Mais l'isomorphisme (inverse) de l'isomorphisme cherch\'e est donn\'e par le calcul du lemme \ref{isom.} qui nous fournit, dans le cas \ $\gamma' = 0$, o\`u \ $(\alpha,\beta,\gamma)$ \ est consid\'er\'e comme param\`etre holomorphe dans \ $X \times \C$,  des sections holomorphes \ $U, V$ \ de \ $\mathcal{O}_{X\times \C}[[b]]$. L'isomorphisme (inverse) de l'isomorphisme cherch\'e est obtenu en envoyant  le g\'en\'erateur \ $e_3$ \ de la famille standard sur \ $\varepsilon_3(\gamma' = 0) : = e_3 + U.e_2 + V.e_1$, qui est le g\'en\'erateur de la famille \ $\pi^*(\mathbb{E})$. $\hfill \blacksquare$

\newpage

\section{Appendices.}

\subsection{Un lemme.}

Le r\'esultat suivant jouant un r\^ole clef dans la construction des bases standards, et donc dans la construction des familles verselles de th\`emes primitifs, nous en donnons ici les grandes lignes de la preuve pour la commodit\'e du lecteur.

\begin{lemma}\label{Ext}
Soient \ $E$ \ et \ $F$ \ deux (a,b)-modules r\'eguliers. Alors on a
$$ \dim_{\C}(Ext^1_{\A}(E,F)) - \dim_{\C}(Ext^0_{\A}(E,F)) = rg(E).rg(F).$$
\end{lemma}

\parag{Preuve} Commen{\c c}ons par montrer le cas o\`u \ $rg(E) = 1$.\\
 On a alors \ $E \simeq \A\big/\A.(a - \lambda.b)$ \ et donc \ $Ext^0_{\A}(E,F)$ \ et \ $Ext^1_{\A}(E,F)$ \ sont respectivement les noyaux et conoyaux de \ $a - \lambda.b : F \to F $. Montrons alors la formule par r\'ecurrence sur le rang de \ $F$. En rang 1 on a \ $F \simeq \A\big/\A.(a - \mu.b)$, et le calcul est \'el\'ementaire :
\begin{enumerate}
\item Pour \ $\lambda \not\in \mu + \mathbb{N}$ \ on a \ $Ext^0_{\A}(E,F) = \{0\}$ \ et \ $Ext^1_{\A}(E,F) \simeq \C.e_{\mu}$.
\item On a \ $Ext^0_{\A}(E,F) = \C.b^n.e_{\mu}$ \ et \
$Ext^1_{\A}(E,F) \simeq \C.e_{\mu} \oplus \C.b^{n+1}.e_{\mu}$ \ pour \ $\lambda = \mu + n $.
\end{enumerate}
D'o\`u l'assertion dans ce cas.\\
Faisons une r\'ecurrence sur l'entier \ $rg(F)$. \\
Si  \ $rg(F) \geq 2$, on a une suite exacte
$$ 0 \to G \to F \to E_{\mu} \to 0 $$
avec \ $ rg(G) = rg(F) - 1$ \  qui donnera la suite exacte d'espaces vectoriels de dimensions finies (d'apr\`es le th\'eor\`eme 1 de [B.95])
\begin{align*}
& 0 \to Ext^0_{\A}(E,G) \to Ext^0_{\A}(E,F)\to Ext^0_{\A}(E,E_{\mu})\to \\
& \to Ext^1_{\A}(E,G) \to Ext^1_{\A}(E,F) \to Ext^1_{\A}(E,E_{\mu})\to 0 
\end{align*}
qui donne que la somme altern\'ee des dimensions est nulle, ou encore
\begin{align*}
&  \dim(Ext^1_{\A}(E,G)) -  \dim(Ext^0_{\A}(E,G)) + \dim( Ext^1_{\A}(E,E_{\mu})) - \dim( Ext^0_{\A}(E,E_{\mu})) = \\
& \dim(Ext^1_{\A}(E,F)) - \dim(Ext^1_{\A}(E,F)) = (rg(F) - 1) + 1 = rg(F)
\end{align*}
gr\^ace \`a l'hypoth\`ese de r\'ecurrence. \\
Le cas o\`u \ $E$ \ est arbitraire et \ $F$ \ est de rang 1 s'obtient de fa{\c c}on analogue.\\
Enfin une r\'ecurrence maintenant sur l'entier \ $rg(E) + rg(F)$ \ donne le cas g\'en\'eral, \`a nouveau par un raisonnement analogue. $\hfill \blacksquare$

 \subsection{Exemple}\label{App.1}

Les in\'egalit\'es du th\'eor\`eme \ref{inclusion} sont pr\'ecises puisque dans l'exemple ci-apr\`es le th\`eme   \ $E' : = E\big/F_1$ \ de rang \ $3$, ne s'injecte pas dans \ $F_3$ \ puisque \ $E$ \ n'est pas stable, alors que l'on a \ $\mu_i - \lambda_i \geq 3-2 = 1$ \ (et m\^eme \ $\mu_2-\lambda_2 = k-1 = 2$) pour ces deux th\`emes.

\bigskip

Je d\'etaille l'exemple suivant :  $ k = 4$ , $p_1 = p_3 = 2$ \ et \ $p_2 = 3$. On a donc \ $q_1 = 5, q_2 = p_2 = 3$ \ et \ $q_3 = p_3 = 2$. On pose \ $E : = \A\big/\A.P $ \ avec 
$$ P : = (a - \lambda_1.b)S_1^{-1}.(a - \lambda_2.b).S_2^{-1}.(a - \lambda_3.b).S_3^{-1}.(a -\lambda_4.b) $$
\begin{align*}
& S_1 : = 1+ \delta.b + \varepsilon.b^2 + \theta.b^5 \\
& S_2 : = 1 + \beta.b + \gamma.b^3 \\
& S_3 : = 1 + \alpha.b^2 \quad {\rm et} \quad \alpha.\gamma.\varepsilon \not= 0 
\end{align*}

\begin{lemma}
L'espace vectoriel \ $Hom_{\A}(E,E)$ \ est de dimension \ $3$, si l'on a  \\
 $\alpha+ \varepsilon \not= 0 $.
\end{lemma}

\parag{preuve} On a un homomorphisme (unique \`a un scalaire multiplicatif non nul pr\`es) de rang 1 de \ $E$ \ dans \ $E$ : il envoie le g\'en\'erateur \ $e$ \ d'annulateur \ $\A.P$ \ sur \ $z_1 : = b^{\lambda_4 - \lambda_1}.e_1$ \ o\`u \ $e_1$ \ est un g\'en\'erateur de \ $F_1$ \ v\'erifiant \ $ (a - \lambda_1.b).e_1 = 0 $.
En effet \ $z_1$ \ est annul\'e par \ $P$, l'\'el\'ement \ $e_1$ \ est unique \`a un scalaire multiplicatif pr\`es, et l'image d'un homomorphisme de rang 1  est de normalis\'e \'egal \`a \ $F_1$. Comme c'est un quotient de \ $E$ de rang 1 il est isomorphe \`a \ $E\big/F_3 \simeq E_{\lambda_4}$.

\bigskip

L'espace vectoriel des homomorphismes de rang 4 modulo ceux de rang 3 est de dimension 1 et engendr\'e par l'identit\'e.

\bigskip

Cherchons maintenant la dimension de l'espace vectoriel des homomorphismes de rang 2 modulo ceux de rang \ $\leq 1$. Un tel homomorphisme a son image dans \ $F_2$ \ puisque le normalis\'e de l'image est \ $F_2$, il est donc donn\'e par un \'el\'ement \ $z_2 \in F_2 \setminus F_1$ \ qui est annul\'e par \ $P$. D'apr\`es le lemme \ref{lemme manquant}, \ $z_2$ \ v\'erifie alors 
$$ (a - \lambda_4.b).z_2 \in F_1 \quad {\rm et} \quad (a - \lambda_3.b).S_3^{-1}.(a - \lambda_4.b).z_2 = 0  $$
et il est donc de la forme
$$ z_2 : = \rho.b^{\lambda_4-\lambda_2}.e_2 + U.e_1 $$
o\`u \ $U \in \C[[b]]$. On a alors
\begin{align*}
& (a - \lambda_4.b).z_2 = b^{\lambda_4-\lambda_2}S_1.e_1 + (\lambda_1-\lambda_4).b.U.e_1 + b^2.U'.e_1 \\
& (a - \lambda_4.b).z_2 = \sigma.b^{\lambda_3-\lambda_1}.S_3.e_1.
\end{align*}
On doit donc avoir
$$ b^2.U' - 4b.U = \sigma.S_3.b^3 - \rho.S_1.b^3 .$$
Apr\`es simplification par \ $b$ \ on obtient une \'equation qui n'a de solution \ $U \in \C[[b]]$ \ que si le coefficient de \ $b^4$ \ dans \ $\sigma.S_3.b^2 - \rho.S_1.b^2 $ \  est nul. Ceci impose la relation 
$$ \sigma.\alpha = \rho.\varepsilon .$$
Comme \ $\alpha, \rho, \varepsilon$ \ sont non nuls, on a un choix (unique pour \ $\rho \in \C^*$ \ donn\'e) qui est non nul pour \ $\sigma$. On a alors une solution unique pour \ $U$ \ modulo \ $\C.b^4$.\\
On en conclut que l'espace vectoriel des homomorphismes de rang \ $\leq 2$ \ est de dimension 2.

\bigskip

Pour achever la preuve du lemme, il suffit de montrer que sous nos hypoth\`eses, il n'existe pas d'homomorphisme de rang 3 de \ $E$ \ dans \ $E$. Ceci revient \`a montrer  qu'il n'existe pas d'\'el\'ement
$$ z_3 = b^{\lambda_4-\lambda_3}.e_3 + V.e_2 + W.e_1 $$
avec \ $V, W \in \C[[b]]$, v\'erifiant \ $$(a - \lambda_4.b).z_3 = S_3.x_2 \quad {\rm   avec} \quad  
 (a - \lambda_2.b).S_2^{-1}.(a - \lambda_3.b).x_2 = 0 , \ x_2 \in F_2 \setminus F_1.$$
 Posons \ $x_2 : = \tau.b^{\lambda_3-\lambda_2}.e_2 + Z.e_1$ \ avec \ $\tau \not= 0$ \ et \ $Z \in \C[[b]]$.
 Alors on a
 $$ (a - \lambda_3.b).x_2 = \tau.b^{\lambda_3-\lambda_2}.S_1.e_1 + b^2.Z'.e_1 - (\lambda_3-\lambda_1).Z.b.e_1 = \eta.S_2.b^{\lambda_2-\lambda_1}.e_1. $$
 On a donc
 $$ b.Z' - 3Z = \eta.S_2 - \tau.S_1.b .$$
 Cette \'equation n'aura de solution \ $Z \in \C[[b]]$ \ que si le coefficient de \ $b^3$ \ dans le membre de droite est nul. Ceci impose la condition \ $\eta.\gamma = \tau. \varepsilon$. On aura alors \ $- 3Z(0) = \eta $ \ pour chaque solution \ $Z$, puisqu'elle est unique modulo \ $\C.b^3$.
 
 \bigskip
 
 Calculons
 \begin{align*}
 & (a - \lambda_4.b).z_3 = b^{p_3-1}.S_2.e_2 + b^2.V'.e_2 + (\lambda_2-\lambda_4).V.b.e_2 + \\
 & \qquad \qquad + V.S_1.e_1 + b^2.W'.e_1 + (\lambda_1-\lambda_4).W.b.e_1 = S_3.x_2 
 \end{align*}
  ce qui donne les \'equations
  \begin{align*}
  & b.V' - 3V = \tau.S_3.b - S_2  \\
  & b^2.W' - 4b.W = S_3.Z - V.S_1
  \end{align*}
  La premi\`ere \'equation n'a de solution que si \ $ \tau.\alpha = \gamma$ \ et on aura alors \ $ -3V(0) = -1$. La seconde ne peut avoir de solution que si  \ $Z(0) = V(0)$ \ puisque le membre de gauche est dans \ $b.\C[[b]]$. Ceci impose \ $\eta = -1$ \ et donc \ $\gamma + \tau. \varepsilon = 0 $. On doit donc avoir \ $\tau = \gamma/\alpha = -\gamma/\varepsilon $ \ ce qui est impossible pour \ $\alpha + \varepsilon \not= 0 $, puisque \ $\alpha.\gamma.\varepsilon \not= 0$. $\hfill \blacksquare$
  
  \parag{Remarque}
  On constate que pour \ $\alpha + \varepsilon = 0$, on trouve une solution \ $W \in \C[[b]]$ \ car le coefficient de \ $b^5$ \ dans \ $S_3.Z - V.S_1$ \ peut \^etre supprim\'e puisque  \ $Z$ \ est d\'efini modulo \ $b^3$ \ et que  \ $\alpha$, le coefficient de \ $b^2$ \ dans \ $S_3$, est non nul. Donc pour \ $\alpha + \varepsilon = 0$ \ le th\`eme \ $E$ \ est stable. $\hfill \square$
  
  \bigskip
  
  Il est facile de d\'eduire de ce qui pr\'ec\`ede que la classe d'isomorphisme de \ $E$ \ est ind\'ependant de \ $\theta \in \C$ \ pour \ $\alpha + \varepsilon \not= 0$\footnote{Noter que comme \ $P_1.F_3 \cap F_1 \subset b^3.F_1$ \ seul \ $\theta$ \ peut changer. Et on voit qu' il change effectivement gr\^ace au lemme \ref{instable 1}.}. On a donc une situation analogue \`a celle d\'ecrite au paragraphe 4.4, c'est-\`a-dire que l'on peut construire une famille universelle pour les th\`emes de rang 4 d'invariants fondamentaux \ $p_1 = p_3 = 2, p_2 = 3$ \ tels que \ $\alpha + \varepsilon \not= 0$, param\'etr\'ee par \ $\{(\alpha,\gamma,\varepsilon),(\beta,\delta) \in (\C^*)^3\times \C^2, \alpha + \varepsilon \not= 0 \}$. \\
  Une preuve analogue permet de montrer que pr\`es des th\`emes (stables, non sp\'eciaux) v\'erifiant \ $\alpha + \varepsilon = 0$, il n'existe pas de famille universelle.

\subsection{Existence d'applications k-th\'ematiques.}\label{App.2}
 
 D'abord un lemme de g\'eom\'etrie alg\'ebrique sur l'alg\`ebre \ $Z : = \C[[b]]$.
 
 \begin{lemma}\label{rg. sem. cont.1}
 Soit \ $E$ \ un (a,b)-module r\'egulier de rang \ $k$. Fixons une \ $\C[[b]]-$base  \ $e_1, \dots, e_k$ \ de \ $E$ \ et consid\'erons \ $E$ \ comme l'espace affine \ $Z^k$ \ sur la \ $\C-$alg\`ebre  \ $Z : =  \C[[b]]$.
 Pour chaque entier \ $p$ \ le sous-ensemble \ $X_p \subset E = Z^k$ \ d\'efini par
 $$ X_p : = \{ x \in E \ / \ rg( \A.x) \leq p \} $$
 est un sous-ensemble alg\'ebrique de \ $E = Z^k$, c'est-\`a-dire qu'il existe un ensemble fini de polyn\^omes \ $P_1, \dots, P_N$ \ dans \ $Z[x_1, \dots, x_k]$ \ tel que l'on ait
 $$ X_p = \{ x \in Z^k \ / \ P_j(x) = 0 \quad \forall j \in [1,N] \}.$$
 \end{lemma}
 
 \parag{Preuve} Comme \ $E$ \ est r\'egulier de rang \ $k$, pour chaque \ $x \in E$, le sous-(a,b)-module \ $\A.x$ \ est monog\`ene r\'egulier de rang \ $\leq k$. Il est donc engendr\'e sur \ $\C[[b]]$ \ par \ $\{x, a.x, \dots, a^{k-1}.x \}$. Pour \'ecrire que le rang de \ $\A.x$ \ est \ $\leq p$, il suffit d'\'ecrire que tous les mineurs \ $(q,q)$ \ de la matrice de ces \ $k$ \ vecteurs dans la base \ $e_1, \dots, e_k$ \ sont nuls pour \ $p+1 \leq k $, ce qui fournit les polyn\^omes \ $P_1, \dots, P_N$ \ de l'\'enonc\'e. $\hfill \blacksquare$
 
 \bigskip
 
 Et une cons\'equence imm\'ediate :
 
 \begin{cor}\label{rg. sem. cont.2}
 Soit \ $X$ \ un espace complexe r\'eduit et soit \ $E$ \ un (a,b)-module r\'egulier de rang \ $k$. Soit \ $f : X \to E$ \ une application holomorphe\footnote{ En fixant une base  \ $e_1, \dots, e_k$ \ c'est une section globale du faisceau \ $\mathcal{O}_X[[b]]^k$.}. On a une stratification finie
 $$ X_0 \subset X_1 \subset \dots \subset X_k = X $$
 par des sous-ensembles analytiques ferm\'es telle que, pour chaque \ $q \in [1,k]$\ le sous-ensemble  \ $X_q \setminus X_{q-1}$ \ soit exactement l'ensemble des \ $x \in X$ \ tels que le rang de \ $\A.f(x)$ \ soit \'egal \`a \ $q$.
 \end{cor}
 
 \parag{Remarque} Le quotient de deux fonctions holomorphes \ $f : D \to \C[[b]]$ \ et \\
  $g : D \to \C[[b]]$ \ avec \ $g(0) \not= 0$ \ peut \^etre bien d\'efini pour chaque valeur de \ $z \in D$, sans pour autant que \ $f/g$ \ soit holomorphe sur \ $D$. Par exemple \ $z \to \frac{z+b^2}{z+b} $ \ est bien d\'efini pour chaque valeur de \ $z \in D$, mais elle n'est cependant pas holomorphe. En effet une relation
 $$ z + b^2 = (z+b).(\sum_{j=0}^{\infty} \ a_j(z).b^j) $$
 conduit imm\'ediatement \`a \ $a_0 \equiv 1$ \ et \ $a_1 = \frac{-1}{z}$ ! $\hfill \square$
 
 \begin{lemma}\label{mero}
 Soient \ $f , g : X \to \C[[b]]$ \ deux applications holomorphes d'un espace complexe r\'eduit \ $X$ \ \`a valeurs dans \ $\C[[b]]$. Supposons \ $X$ \ irr\'eductible et \ $g \not\equiv 0$. Supposons que pour chaque \ $x \in X$ \ le quotient \ $f(x)/g(x)$ \ soit dans \ $\C[[b]]$. Alors il existe  un ouvert de Zariski dense \ $X'$ \ de \ $X$ \ sur lequel l'application \\
  $x \to f(x)/g(x) \in \C[[b]]$ \ est holomorphe.
 \end{lemma}
 
 \parag{Preuve} On peut supposer que \ $f \not\equiv 0$ \ sur \ $X$. Il existe alors deux ouverts de Zariski denses \ $X_1$ \ et \ $X_2$ \ tels que sur \ $X_1$ (resp. sur \ $X_2$) la valuation en \ $b$ \  de \ $f(x)$ \ (resp. de \ $g(x)$) soit constante \'egale \`a \ $k$ \ (resp. \`a \ $l$). La condition impos\'ee montre que l'on a \ $k \geq l$, et sur \ $X_1\cap X_2$ \ on peut \'ecrire
 $$ f(x) = b^k.F(x) \quad g(x) = b^l.G(x) $$
 o\`u \ $F, G $ \ sont des fonctions holomorphes \`a valeurs inversibles dans \ $\C[[b]]$. Il ne reste plus qu'\`a se convaincre que la fonction \ $x \mapsto b^{k-l}.F(x)/G(x)$ \ est holomorphe sur \ $X_1 \cap X_2$, ce qui est \'el\'ementaire. $\hfill \blacksquare$
 
 \begin{cor}\label{exist. k-thm.}
 Soit \ $f : X \to E$ \ une application holomorphe d'un espace complexe r\'eduit et irr\'eductible dans un (a,b)-module r\'egulier. Il existe un ouvert dense \ $X'$ \ de \ $X$ \ sur lequel la restriction de \ $f$ \ d\'efinit une application \ $k-$th\'ematique via \ $x \mapsto \A.f(x)$, o\`u \ $k \leq rg(E)$.
 \end{cor}
 
 \parag{Preuve} Le point est que l'on trouve un  ouvert de Zariski \ $X'$ \  sur lequel le rang du (a,b)-module monog\`ene \ $\A.f(x)$ \ est constant gr\^ace au premier lemme. On r\'esoud  ensuite le syst\`eme de Cramer avec param\`etre sur cet ouvert dense, mais on trouve, pour les fonctions \ $x \mapsto S_j(x)\in \C[[b]]$ \ donnant la relation
 $$ a^k.f(x) = \sum_{j=0}^{k-1} \ S_j(x).a^j.f(x) $$
 des fonctions m\'eromorphes. Le second lemme donne alors un ouvert de Zariski \ $X''$ \ de \ $X'$ \ sur lequel ces fonctions sont holomorphes. $\hfill \blacksquare$
 
 \parag{Remarque} Dans le corollaire ci-dessus on prendra garde que l'ouvert dense trouv\'e est un ouvert de Zariski d'un ouvert de Zariski de \ $X$, qui n'est pas, en g\'en\'eral, un ouvert de Zariski de \ $X$. $\hfill \square$
 
\newpage

\section{R\'ef\'erences.}

\begin{itemize}

\item{[A-G-V]} Arnold,V. Goussein-Zad\'e, S. Varchenko,A. {\it Singularit\'es des applications diff\'erentiables}, volume 2,  \'edition MIR, Moscou 1985.

\item{[B. 93]}  Barlet, D. \textit{Th\'eorie des (a,b)-modules I}, Univ. Ser. Math. Plenum (1993), p. 1-43.

\item{[B. 97]} Barlet, D. \textit{Th\'eorie des (a,b)-modules II. Extensions}, Pitman Res. Notes Math. Ser. 366, Longman (1997), p. 19-59.

\item{[B. 05]}  Barlet, D. \textit{Module de Brieskorn et formes hermitiennes pour une singularit\'e isol\'ee d'hypersurface}, Revue Inst. E. Cartan (Nancy) vol. 18 (2005), p. 19-46.

\item{[B. 09]} Barlet, D. \textit{ P\'eriodes \'evanescentes et (a,b)-modules monog\`enes}, preprint Institut E. Cartan (Nancy) 2009 \ $n^0 1$, p. 1- 46. A para\^itre dans le Bulletino U.M.I. (9) II (2009).

 \item{[Br.70]} Brieskorn, E. {\it Die Monodromie der Isolierten Singularit{\"a}ten von Hyperfl{\"a}chen}, Manuscripta Math. 2 (1970), p. 103-161.
 
 \item{[M. 75]} Malgrange, B. {\it Le polyn\^ome de Bernstein d'une singularit\'e isol\'ee}, in Lect. Notes in Math. 459, Springer (1975), p.98-119.

\item{[S. 89]} Saito, M. {\it On the structure of Brieskorn lattices}, Ann. Inst. Fourier 39 (1989), p.27-72.

\end{itemize}

\end{document}